\documentclass[onefignum,onetabnum]{siamart190516}

\usepackage{adjustbox, algpseudocode, amsmath, amssymb, blindtext, bm, bbm, dblfloatfix, esint, fancyhdr, float, graphicx, letltxmacro, marginnote, mathtools, subcaption, xcolor, titlesec, esint}
\usepackage[shortlabels]{enumitem}
\usepackage[linesnumbered,algoruled,algo2e]{algorithm2e}
\usepackage{lipsum}
\usepackage{amsfonts}
\usepackage{epstopdf}
\usepackage{tablefootnote}

\usepackage{listings}
\usepackage{color}

\usepackage{cite}

\newsiamremark{assumption}{Assumption}
\newsiamremark{remark}{Remark}
\newsiamremark{example}{Example}

\headers{Extended Galerkin Neural Networks}{M. Ainsworth and J. Dong}

\DeclareMathOperator*{\argmax}{arg\,max}
\DeclareMathOperator*{\argmin}{arg\,min}

\title{Extended Galerkin neural network approximation of singular variational problems with error control}

\author{Mark Ainsworth\thanks{Division of Applied Mathematics, Brown University
  (\email{mark\_ainsworth@brown.edu}).}
\and Justin Dong\thanks{Center for Applied Scientific Computing, Lawrence Livermore National Laboratory
  (\email{dong9@llnl.gov}).}}

\begin{document}
	
\maketitle
	
	\begin{abstract}
            We present extended Galerkin neural networks (xGNN), a variational framework for approximating general boundary value problems (BVPs) with error control. The main contributions of this work are (1) a rigorous theory guiding the construction of new weighted least squares variational formulations suitable for use in neural network approximation of general BVPs (2) an ``extended'' feedforward network architecture which incorporates and is even capable of learning singular solution structures, thus greatly improving approximability of singular solutions. Numerical results are presented for several problems including steady Stokes flow around re-entrant corners and in convex corners with Moffatt eddies in order to demonstrate efficacy of the method.

	\end{abstract}
	
	\vspace{5mm}
	

	

\section{Introduction}

Neural network methods for approximating the solutions of partial differential
equations (PDEs) have seen a surge of interest in recent years
\cite{dissanayake, lagaris, pinns, wan} including extensions to inverse
problems \cite{pinns, pakravan2021solving}, learning solution operators
\cite{deeponet,fourieronet} of PDEs, and even complementing existing numerical
methods, e.g. by learning finite difference stencils \cite{bar2019learning} and
guiding adaptive mesh refinement \cite{yang2022multi}. Despite impressive
results that are often reported, many current techniques have not yet been set on a
firm theoretical foundation and the prospective user may be wary of adopting
approaches where critical decisions may be made on the basis of the numerical
approximation.

In our previous work~\cite{gnn1} we aimed to develop effective techniques for
using neural networks to approximate the solutions of PDEs which are,
crucially, supported by a rigorous theory of convergence and the provision of a
computable \emph{a posteriori} error estimator for the error in the resulting
approximation. The approach in \cite{gnn1} gives an iterative procedure in
which the networks are successively enriched until the error, estimated using
the \emph{a posteriori} error estimator, meets a desired tolerance. Examples presented
in~\cite{gnn1} showed the approach to be effective for second and fourth order
self-adjoint elliptic PDEs while the application to a singularly-perturbed
elliptic system with boundary layers was presented in \cite{gnn2}.

One purpose of the current work is to describe an extension of the Galerkin
neural network approach in previous work~\cite{gnn1,gnn2} to quite broad
classes of non-self-adjoint and/or indefinite problems of the general form
        \begin{align} \label{eq:general bvp intro}
                \begin{dcases}
                        \mathcal{L}[u] = f &\text{in}\;\Omega\\
                        \mathcal{B}[u] = g &\text{on}\;\partial\Omega,
                \end{dcases}
        \end{align}
where $\mathcal{L}$ and $\mathcal{B}$ are linear differential operators. The
only assumption needed is the the rather mild condition that the problem is
well-posed in the sense that it admits a unique solution that depends
continuously on the data using appropriate norms. Under this condition, we show
that our previous results extend to the more general scenario given
in~\eqref{eq:general bvp intro} including the provision of a computable \emph{a posteriori} error estimator based on using a least squares formulation. While least
squares formulations have featured in previous work using neural networks to
approximate PDEs in a somewhat ad hoc fashion -- e.g. mean-squared error (MSE) or $L^{2}$ norm as a loss metric -- we argue that the least squares
formulation is dictated by the properties of the underlying continuous
problem~\eqref{eq:general bvp intro} rather than the whim of the user if one
wishes to obtain a rigorous convergence theory. Ad hoc choices of formulation correspond to a BVP which may not be uniquely solvable or for which there is no continuous dependence on the data. We shall demonstrate in this work that the use of such formulations can result in nonphysical solutions (see Examples \ref{ex:poisson rlambda} and \ref{ex:triangular})


The universal approximation property of neural networks is often cited as the driving force behind adopting them for traditional scientific computing tasks. This means
that neural network based approaches are applicable to broad classes of problems
in which solutions may exhibit widely varying behaviour.  Nevertheless, in many
common applications such as computational fluid dynamics and solid mechanics,
solutions of \eqref{eq:general bvp intro} exhibit characteristic features such
as boundary layers and singularities even when the data is smooth. Such
features generally result in a degradation in the rate of convergence of
traditional methods, such as standard finite element methods. In principle,
neural network approaches are capable of approximating such features at least
as well as adaptive finite element methods~\cite{schwab1998p} provided that the
network is trained appropriately. In practice, the training of networks is
often a bottleneck in the approach with the net result that the overall rate of
convergence that is achieved falls short of what is theoretically possible.

The second purpose of the current work is to show how knowledge of the characteristic features of the solution to a given PDE can be naturally incorporated into our approach to improve rate of convergence. Some traditional methods such as extended finite element methods \cite{xfem} and generalized finite element methods \cite{genfem} utilize prior knowledge of these features. More recent work has explored the use of tailored functions to supplement deep neural networks in specific examples, such as \cite{arzani2023theory} which utilizes boundary layer functions in conjunction with feedforward networks to learn simple boundary layers in a singularly-perturbed problem. In the current work, we demonstrate how extended Galerkin neural networks can be used to both incorporate and \emph{learn} singular solution features. Unifying theory is provided which demonstrates that the size of the networks needed to resolve a particular solution with extended Galerkin neural networks depends only on the smooth, non-characteristic, part of the solution which is more easily approximated using standard feedforward architectures, thus resulting in greatly improved approximation rates. Numerical examples show that the proposed method can be highly effective for problems exhibiting singularities and features such as eddies. 


The rest of this work is structured as follows. In Section \ref{sec:gnn}, we
present the extended Galerkin neural network framework, including presentation
and analysis of a least squares variational formulation for a very general
class of PDE as well as the enriched neural network architecture of xGNN. We
describe in Sections \ref{sec:applications}-\ref{sec:eigenvalue training} how
the preceding theory can be applied to several interesting applications,
including Stokes flow in polygonal domains. Conclusions follow in Section
\ref{sec:conclusions}.

	\section{Extended Galerkin neural networks (xGNN)} \label{sec:xgnn}

        We first give a brief summary of the basic Galerkin neural network framework for symmetric, positive-definite (SPD) problems in Section \ref{sec:gnn} following \cite{gnn1}. We then demonstrate how to rigorously extend this approach to non-symmetric and/or indefinite and negative-definite problems in Section \ref{sec:lsq} as well as enrich the neural network approximation space with known solution structures in Section \ref{sec:knowledge-based fn}. Collectively, we refer to these advancements as extended Galerkin neural networks (xGNN). 
        
	\subsection{The basic Galerkin neural network framework} \label{sec:gnn}
	
	The basic Galerkin neural network framework \cite{gnn1} considers the prototypical variational problem
	\begin{align} \label{eq:variational}
		u \in X \;:\; a(u,v) = F(v) \;\;\;\forall v \in X,
	\end{align}
	
	\noindent where $a$ is assumed to be an SPD, continuous, and coercive bilinear form with respect to $(X, ||\cdot||_{X})$ which defines an energy norm $|||\cdot||| := a(\cdot,\cdot)^{1/2}$, and $F$ is assumed to be continuous. Given an initial approximation $u_{0} \in X$ to \eqref{eq:variational}, the goal of the Galerkin neural network is to iteratively construct a finite-dimensional subspace $S_{j} := \text{span}\{u_{0}, \varphi_{1}^{NN}, \dots, \varphi_{j}^{NN}\}$ such that the basis function $\varphi_{i}^{NN}$ satisfies
	\begin{align} \label{eq:basis fn}
		\varphi_{i}^{NN} = \argmax_{v \in V_{\mathbf{n}^{(i)},L_{i}}^{\sigma,C} \cap B} \langle r(u_{i-1}), v \rangle,
	\end{align}
	
	\noindent where 
    \begin{align*}
        \langle r(u_{i-1}),v \rangle = F(v) - a(u_{i-1},v)
    \end{align*}
    
    \noindent is the weak residual in the current approximation $u_{i-1}$, $B$ is the closed unit ball in $X$, and $V_{\mathbf{n},L}^{\sigma}$, $\mathbf{n} \in \mathbb{N}^{L}$ is the set
    \begin{align}
	V_{\mathbf{n},L}^{\sigma} := &\left\{ v \in X \;:\; v(x) = (\sigma \circ \mathbf{T}_{L} \circ \dots \sigma \circ \mathbf{T}_{1}(x)) \cdot \mathbf{c}, \;\mathbf{T}_{j}(t) := t\cdot \mathbf{W}^{(j)} + \mathbf{b}^{(j)} \right.\notag\\
        &\left. \;\mathbf{W}^{(j)} \in \mathbb{R}^{n_{j-1}\times n_{j}}, \;\mathbf{b}^{(j)} \in \mathbb{R}^{1\times n_{j}}, \;\mathbf{c}\in \mathbb{R}^{n_{L}\times 1}, \;x \in \mathbb{R}^{1\times d} \right\},
    \end{align}

    \noindent with the convention that the data $x$ has dimension $\mathbb{R}^{1\times d}$ and $n_{0} = d$. The set $V_{\mathbf{n},L}^{\sigma}$ describes the realizations of a multilayer feedforward neural network with depth $L$, widths $\mathbf{n}$, and suitable (non-polynomial) activation function $\sigma \in X$. The set $V_{\mathbf{n},L}^{\sigma,C}$ denotes functions in $V_{\mathbf{n},L}^{\sigma}$ with bounded parameters:
    \begin{align}
        V_{\mathbf{n},L}^{\sigma,C} := \left\{ v \in V_{\mathbf{n},L}^{\sigma} : ||(\mathbf{W},\mathbf{b},\mathbf{c})||_{\mathcal{NN}} := \sum_{j=1}^{L} ||\mathbf{W}^{(j)}||_{\infty} + ||\mathbf{b}^{(j)}||_{\infty} + ||\mathbf{c}||_{\infty} \leqslant C \right\}.
    \end{align}

        \noindent The notion of bounded network parameters is a technical necessity to ensure existence of $\varphi_{i}^{NN}$. However, in practice, we observe as in \cite{gnn1} that the parameters corresponding to $\varphi_{i}^{NN}$ are uniformly bounded and thus, no clipping of the parameters is necessary to ensure boundedness. Once the basis function $\varphi_{i}^{NN}$ is known, one generates a new, Galerkin approximation $u_{i}$ using the finite-dimensional space $S_{i}$ as follows:
	\begin{align} \label{eq:ui}
		u_{i} \in S_{i} \;:\; a(u,v) = F(v) \;\;\;\forall v \in S_{i}.
	\end{align}
	
	One can show \cite{gnn1} that the basis function $\varphi_{i}^{NN}$ approximates the normalized error $\varphi_{i} := (u-u_{i-1})/|||u-u_{i-1}|||$ in the current approximation and thus the subspaces $S_{j}$ may be viewed as augmenting the initial approximation $u_{0}$ with a sequence of increasingly finer-scale error correction terms. If the network used to approximate $\varphi_{i}$ is sufficiently rich, then the approximation $u_{i}$ is exponentially convergent in the number of basis functions generated, and one can also show that the objective function $\langle r(u_{i-1}), \varphi_{i}^{NN} \rangle$ (i.e. weak residual) is an \emph{a posteriori} estimator of the energy error $|||u-u_{i-1}|||$.
		
		
    
    One can thus generate an approximation $u_{i}$ such that $|||u-u_{i-1}||| \lessapprox \texttt{tol}$ for a given tolerance \texttt{tol} by checking whether $\langle r(u_{i-1}),\varphi_{i}^{NN} \rangle < \texttt{tol}$ after $\varphi_{i}^{NN}$ is generated. In the affirmative case, we assume that the energy error is within the desired tolerance and terminate the subspace generation procedure. Otherwise, we generate another basis function until the desired tolerance is reached.  
	
	The training procedure for learning $\varphi_{i}^{NN}$ consists of a standard gradient-based step to update the hidden parameters $\mathbf{W}$ and $\mathbf{b}$ in conjunction with the least squares solution of the linear system
	\begin{align} \label{eq:lstsq solve}
		\begin{dcases}
    			\mathbf{A}\mathbf{c} &= \mathbf{F}\\
    			\mathbf{A}_{k\ell} &= a(\sigma(x\cdot \mathbf{W}_{k}+\mathbf{b}_{k}), \sigma(x\cdot \mathbf{W}_{\ell}+\mathbf{b}_{\ell}))\\
    			\mathbf{F}_{k} &= L(\sigma(x\cdot \mathbf{W}_{k}+\mathbf{b}_{k})) - a(u_{i-1}, \sigma(x\cdot \mathbf{W}_{k}+\mathbf{b}_{k})) 
		\end{dcases}
	\end{align}
	
	\noindent in order to update the activation coefficients $\mathbf{c}$. For ease of notation, we consider only the case when $L=1$ in \eqref{eq:lstsq solve} while noting that the linear system for $L>1$ differs only in that $\sigma(x\cdot \mathbf{W}_{k} + \mathbf{b}_{k})$ will be replaced by the corresponding more complicated expression for the $j$th component of $\sigma \circ \mathbf{T}_{L} \circ \dots \sigma \circ \mathbf{T}_{1}(x)$. The notation $\mathbf{W}_{k}$ and $\mathbf{b}_{k}$ are shorthand for $\mathbf{W}[:,j]$ and $\mathbf{b}[1,j]$, respectively. The linear system in \eqref{eq:lstsq solve} corresponds to the orthogonal projection of the error $u-u_{i-1}$ onto the subspace $\Phi := \text{span}\{\sigma(x\cdot \mathbf{W}_{k}+\mathbf{b}_{k})\}_{k=1}^{n}$. 
 
    Altogether, given a parameter initialization $(\mathbf{W}, \mathbf{b}, \mathbf{c})$, we update the parameters by the rules
	\begin{align} \label{eq:GNN update}
		  \mathbf{W} &\leftarrow \mathbf{W} + \nabla_{\mathbf{W}}\left[ \frac{\langle r(u_{i-1}),v \rangle}{|||v|||} \right]\\
		\mathbf{b} &\leftarrow \mathbf{b} + \nabla_{\mathbf{b}}\left[ \frac{\langle r(u_{i-1}),v \rangle}{|||v|||} \right]\\
		\mathbf{c} &\leftarrow \argmin_{\mathbf{c} \in \mathbb{R}^{n}} ||\mathbf{A}\mathbf{c} - \mathbf{F}||_{\ell^{2}}.
	\end{align}

	\subsection{Extension to non-self-adjoint and non-positive-definite problems} \label{sec:lsq}

    The framework described in Section \ref{sec:gnn} considers self-adjoint, positive-definite boundary value problems whose variational formulation naturally gives rise to a bilinear form $a(\cdot,\cdot)$ that is SPD. Here, we extend the approach to a more general boundary value problem with the strong form 
	\begin{align} \label{eq:general bvp}
		\begin{dcases}
			\mathcal{L}[\mathbf{u}] = \mathbf{f} &\text{in}\;\Omega\\
			\mathcal{B}[\mathbf{u}] = \mathbf{g} &\text{on}\;\partial\Omega,
		\end{dcases}
	\end{align}
	where $\mathcal{L} : \mathcal{X} \to \mathcal{V}$ and $\mathcal{B} : \mathcal{X} \to \mathcal{W}$ are linear differential operators, $\mathbf{f} \in \mathcal{V}$, and $\mathbf{g} \in \mathcal{W}$ for appropriate Hilbert spaces $\mathcal{V}, \mathcal{W}$, and $\mathcal{X}$. We shall defer discussion of how to interpret \eqref{eq:general bvp} as a variational problem until \eqref{eq:general lsq}.

    The operator $\mathcal{L}$ corresponds to the differential equation posed over the domain $\Omega$ while $\mathcal{B}$ corresponds to the boundary conditions on the domain boundary $\partial\Omega$. In order to maintain generality, we aim to impose as few requirements on these operators as possible beyond what is necessary for the problem \eqref{eq:general bvp} to be well-posed. A minimal requirement is that the operators should be bounded in the sense that there exists a positive constant $C$ such that 
    $||\mathcal{L}[\mathbf{v}]||_{\mathcal{V}} \leqslant C ||\mathbf{v}||_{\mathcal{X}}$
    and
	$||\mathcal{B}[\mathbf{v}]||_{\mathcal{W}} \leqslant C ||\mathbf{v}||_{\mathcal{X}}$.
In addition, we shall assume that the problem is well-posed in the sense that it admits a unique solution $\mathbf{u} \in \mathcal{X}$ that depends continuously on the data $\mathbf{f}$ and $\mathbf{g}$:
	\begin{align} \label{eq:a priori}
		||\mathbf{u}||_{\mathcal{X}} \leqslant C(||\mathbf{f}||_{\mathcal{V}} + ||\mathbf{g}||_{\mathcal{W}})
	\end{align}
 for some constant $C>0$.

    These assumptions are rather mild and encompass a broad range of boundary value problems that typically arise in physical applications including problems that are non-self-adjoint or which are not positive-definite. Examples will be given later. In order to apply the Galerkin neural network framework described in Section \ref{sec:gnn}, we formulate \eqref{eq:general bvp} in the form \eqref{eq:variational} where the bilinear and linear forms correspond to the following least squares formulation of \eqref{eq:general bvp}:
	\begin{align} \label{eq:general lsq}
		\mathbf{u} \in \mathcal{X} \;:\; a_{LS}(\mathbf{u},\mathbf{v}) = F_{LS}(\mathbf{v}) \;\;\;\forall \mathbf{v} \in \mathcal{X}.
	\end{align}
 where, for a given $\delta \geqslant 1$, we define
\begin{align}
		a_{LS}(\mathbf{u},\mathbf{v}) = (\mathcal{L}[\mathbf{u}], \mathcal{L}[\mathbf{v}])_{\mathcal{V}} 
                                      + \delta(\mathcal{B}[\mathbf{u}], \mathcal{B}[\mathbf{v}])_{\mathcal{W}}
\end{align}
and 
\begin{align}
		F_{LS}(\mathbf{v}) = (\mathbf{f},\mathcal{L}[\mathbf{v}])_{\mathcal{V}} + \delta(\mathbf{g}, \mathcal{B}[\mathbf{v}])_{\mathcal{W}} .
\end{align}

\noindent Variational formulations such as \eqref{eq:general lsq} are employed in least squares finite element methods \cite{bochev} often with a reduction of the PDE to first order. In the current work, we make no such conversion to first-order systems but note that this conversion is fully-compatible with the proposed method. 

One advantage of least squares formulations is that the bilinear form $a_{LS} : \mathcal{X} \times \mathcal{X} \to \mathbb{R}$ is SPD. This follows since $a_{LS}(\mathbf{v}, \mathbf{v}) \geqslant 0$ with equality if and only if $\mathcal{L}[\mathbf{v}] = \mathbf{0}$ and $\mathcal{B}[\mathbf{v}] = \mathbf{0}$ or, by the assumption of the uniqueness of the solution of \eqref{eq:general bvp},  $\mathbf{v} = \mathbf{0}.$ Moreover, the variational problem \eqref{eq:general lsq} is well-posed:
\begin{theorem} \label{thm:general weighted cont coercive}
		Suppose there exist constants $C_{1}, C_{2}>0$ such that
		\begin{align} \label{eq:ctyofLandB}
		||\mathcal{L}[\mathbf{v}]||_{\mathcal{V}} \leqslant C_{1} ||\mathbf{v}||_{\mathcal{X}}
   \mbox{ and }			
			||\mathcal{B}[\mathbf{v}]||_{\mathcal{W}} \leqslant C_{2} ||\mathbf{v}||_{\mathcal{X}}
		\end{align}
  for all $\mathbf{v} \in \mathcal{X}$, and that for all $\mathbf{f} \in \mathcal{V}$ and $\mathbf{g} \in \mathcal{W}$, 
  problem \eqref{eq:general bvp} admits a unique solution $\mathbf{u} \in \mathcal{X}$ that depends continuously on the data
		\begin{align} \label{eq:general a priori}
			||\mathbf{u}||_{\mathcal{X}} \leqslant C_{3}(||\mathbf{f}||_{\mathcal{V}} + ||\mathbf{g}||_{\mathcal{W}})
		\end{align}
  for some constant $C_{3} > 0$. Then the bilinear form $a_{LS}$ is continuous and coercive, and the linear form $F_{LS}$ is continuous in the sense that there exist constants $\mathcal{C}_{1}, \mathcal{C}_{2}, \mathcal{C}_{3} > 0$ such that
		\begin{align*}
			a_{LS}(\mathbf{u}, \mathbf{v}) \leqslant \mathcal{C}_{1} ||\mathbf{u}||_{\mathcal{X}} ||\mathbf{v}||_{\mathcal{X}}, \;\;\;\mathcal{C}_{2}||\mathbf{v}||_{\mathcal{X}}^{2} \leqslant a_{LS}(\mathbf{v}, \mathbf{v})
  \mbox{ and }F_{LS}(\mathbf{v})\leq \mathcal{C}_3 ||\mathbf{v}||_\mathcal{X}
		\end{align*}
  for all $\mathbf{u}, \mathbf{v} \in \mathcal{X}$. Consequently, \eqref{eq:general lsq} is uniquely solvable.
	\end{theorem}
 \begin{proof} 
 Using the Cauchy-Schwarz inequality and continuity of $\mathcal{L}$ and $\mathcal{B}$ immediately shows that the bilinear form will be continuous,
		\begin{align*}
			a_{LS}(\mathbf{u}, \mathbf{v}) &\leqslant ||\mathcal{L}[\mathbf{u}]||_{\mathcal{V}} ||\mathcal{L}[\mathbf{v}]||_{\mathcal{V}} + \delta||\mathcal{B}[\mathbf{u}]||_{\mathcal{W}} ||\mathcal{B}[\mathbf{v}]||_{\mathcal{W}}\\
			&\leqslant (C_1^2+\delta C_2^2)||\mathbf{u}||_{\mathcal{X}} ||\mathbf{v}||_{\mathcal{X}}.
		\end{align*}
Likewise,    
\begin{align}
		F_{LS}(\mathbf{v})  
        \le \left(C_1 ||\mathbf{f}]||_{\mathcal{V}} + \delta C_2||\mathbf{g}||_\mathcal{W} \right) 
			||\mathbf{v}||_{\mathcal{X}} .
\end{align}
Similarly, from \eqref{eq:general a priori} we obtain $||\mathbf{v}||_{\mathcal{X}} \leqslant C(||\mathcal{L}[\mathbf{v}]||_{\mathcal{V}} + ||\mathcal{B}[\mathbf{v}]||_{\mathcal{W}})$ for $\mathbf{v}\in\mathcal{X}$ and, thanks to Young's inequality for products, deduce that the bilinear form is coercive: 
		\begin{align*}
			||\mathbf{v}||_{\mathcal{X}}^{2} \leqslant C_3^2(
                    ||\mathcal{L}[\mathbf{v}]||_{\mathcal{V}}  + ||\mathcal{B}[\mathbf{v}]||_{\mathcal{W}}
            )^{2} \leqslant 2C_{3}^{2}(
                    ||\mathcal{L}[\mathbf{v}]||_{\mathcal{V}}^{2}  + ||\mathcal{B}[\mathbf{v}]||_{\mathcal{W}}^{2}
            ) \leqslant 2C_3^2 a_{LS}(\mathbf{v}, \mathbf{v}).
  \end{align*}
\end{proof}

	We conclude the preceding discussion with a simple example that demonstrates the necessity of choosing $\mathcal{V}$ and $\mathcal{W}$ so that the $a_{LS}$-energy norm is norm equivalent to $||\cdot||_{\mathcal{X}}$. 
	\begin{figure}[t!]
		\centering
		\includegraphics[width=1.7in]{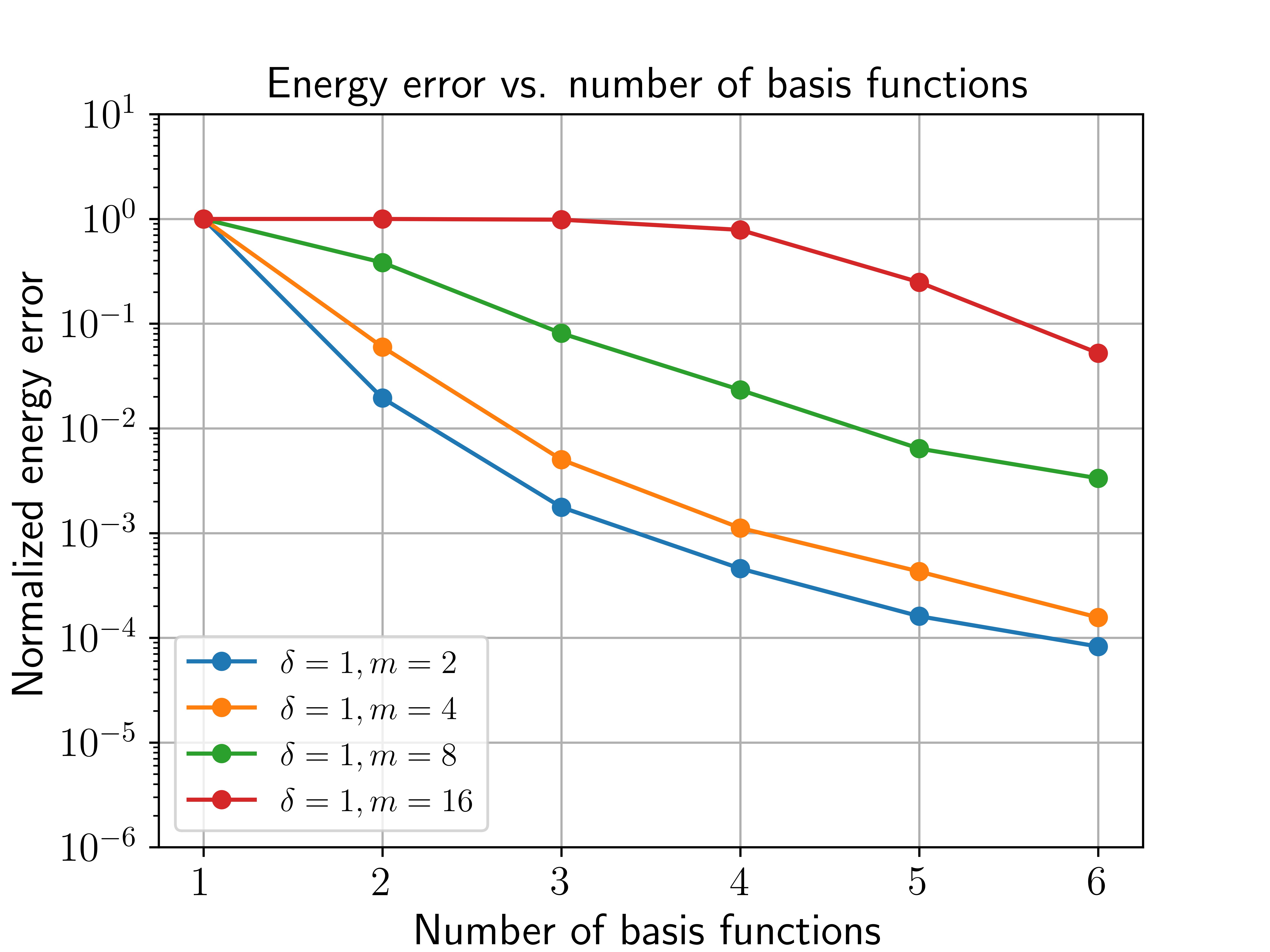}
		\quad
		\includegraphics[width=1.7in]{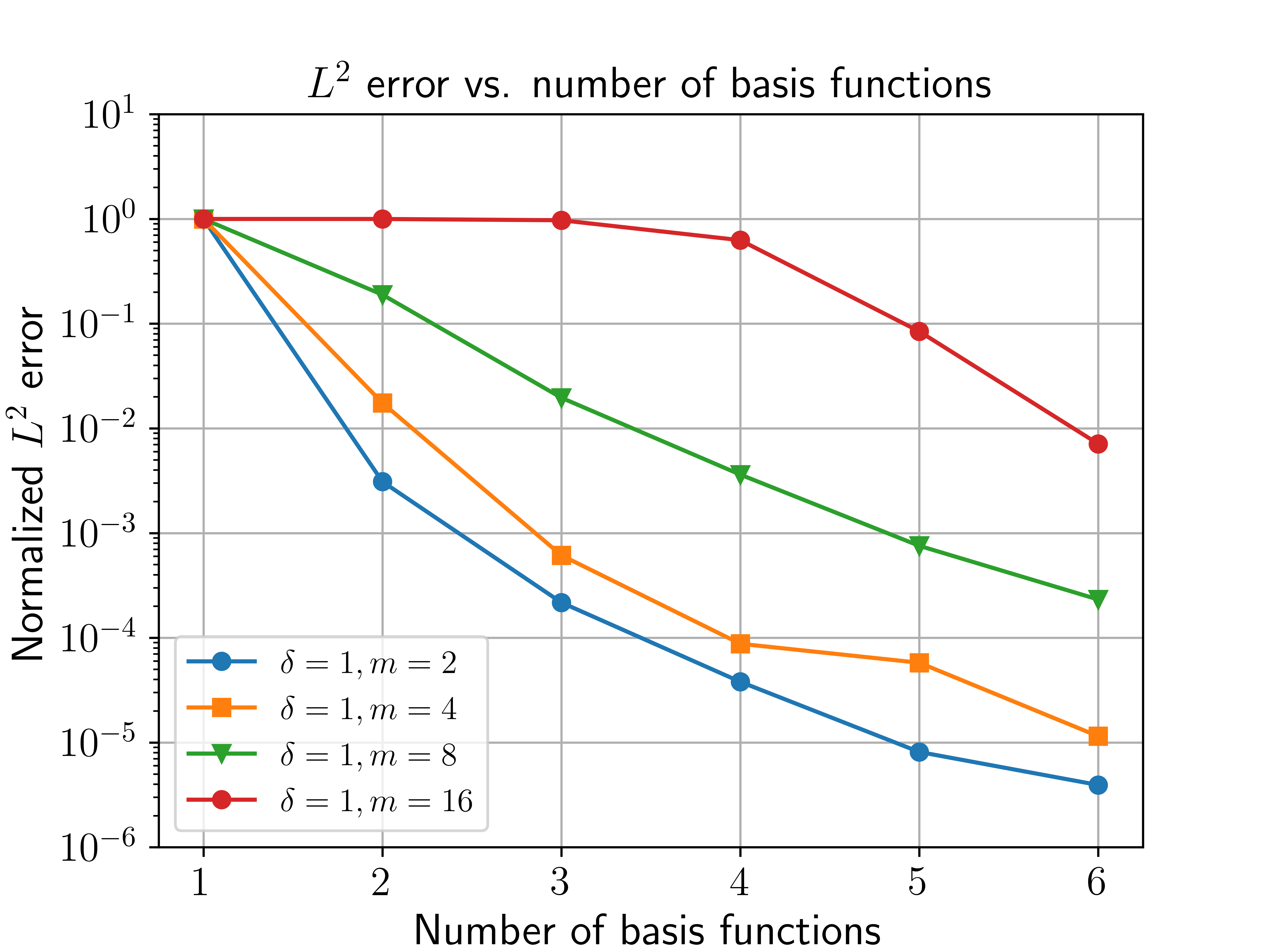}
		\quad
		
		\includegraphics[width=1.7in]{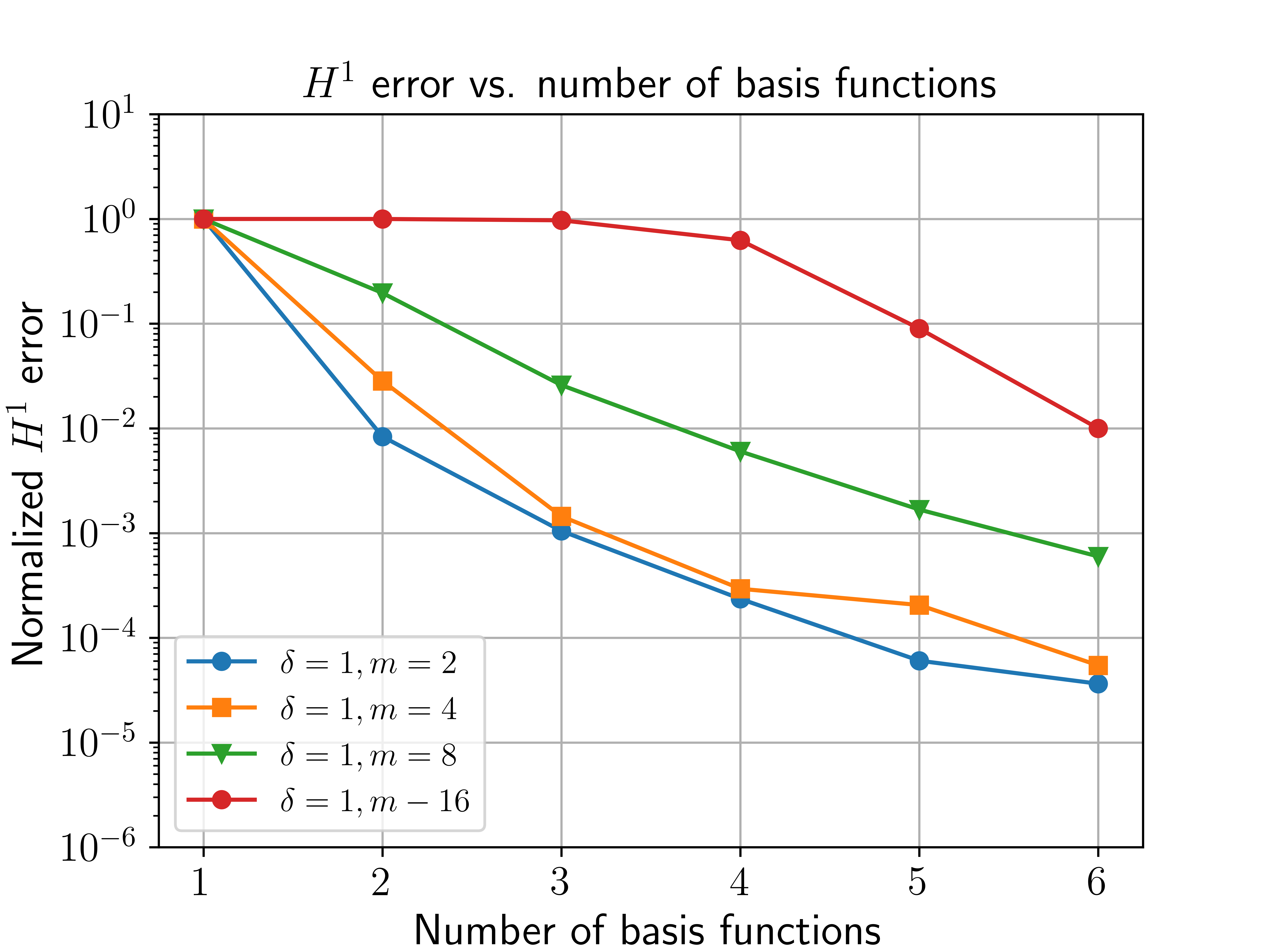}
		\quad
		\includegraphics[width=1.7in]{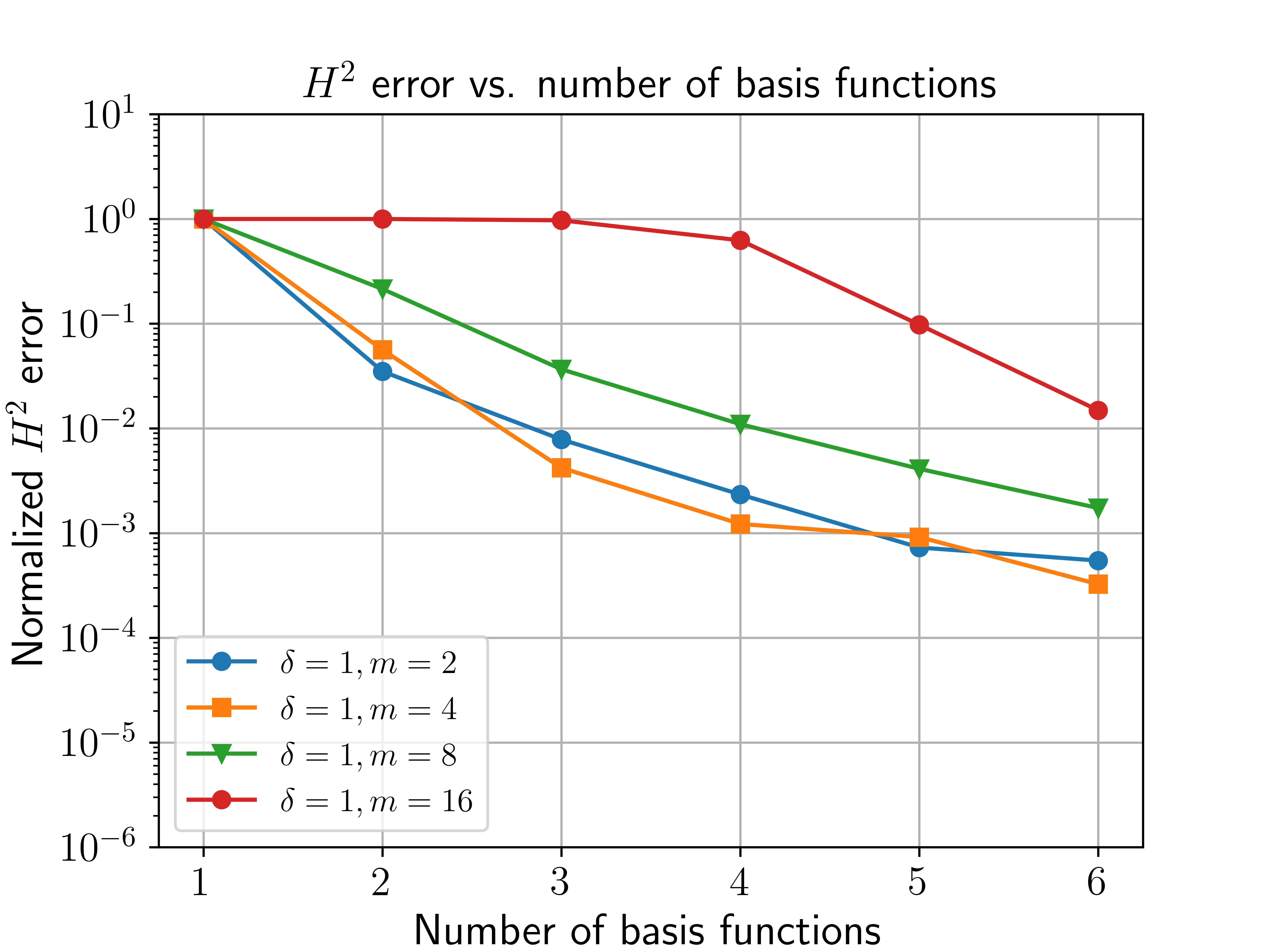}
		\caption{Errors with respect to the number of basis functions learning using the Galerkin neural network approach applied to the variational problem \eqref{eq:coercivity example} for $\delta =1$, $\mathcal{W} = L^{2}(\partial\Omega)$, and the sequence of problems $m=2,4,8,16$.}
		\label{fig:poisson boundaryweights1}
	\end{figure}
 
	\begin{example} \label{ex:cont coercivity}
		Consider a Poisson equation with Dirichlet data so that the corresponding least squares variational problem is given by
          \begin{align} \label{eq:coercivity example}
             u \in \mathcal{X} \;:\; a_{LS}(u,v) = (-\Delta u, -\Delta v)_{\mathcal{V}} + \delta(u,v)_{\mathcal{W}} = F_{LS}(v)\;\;\;\forall v \in \mathcal{X},
        \end{align}

        \noindent where the linear form $F_{LS}$ is chosen so that the true solution is given by the harmonic function $u_{m}(r,\theta) = r^{m}\sin(m\theta)$ for fixed $m \in \mathbb{N}$ in the unit circle.

        Here, $\mathcal{X}=H^{2}(\Omega)$, $\mathcal{V} = L^{2}(\Omega)$, and $\mathcal{W}$ should be taken as $H^{3/2}(\partial\Omega)$ as suggested by the following estimate for the Poisson equation which reads \cite{kondratev}
        \begin{align}
		||u||_{H^{2}(\Omega)} \leqslant C(||f||_{L^{2}(\Omega)} + ||g||_{H^{3/2}(\partial\Omega)}).
	\end{align}
        
        \noindent In practice, physics-informed neural networks (PINNs) are typically applied in conjunction with the choice $\mathcal{W} = H^{s}(\partial\Omega)$, $s=0$, i.e. $\mathcal{W} = L^{2}(\partial\Omega)$. However, choosing $s\neq 3/2$ means that the variational problem \eqref{eq:general lsq} is no longer well-posed.
        
        In order to see this, we argue as follows. A straightforward computation reveals that
	\begin{align} \label{eq:norm calcs}
		||u_{m}||_{H^{2}(\Omega)} \propto m^{3/2}, \;\;\;||u_{m}||_{H^{s}(\partial\Omega)} \propto m^{s}.
	\end{align}

        \noindent Consequently, the constant $\mathcal{C}_{2}$ in Theorem \ref{thm:general weighted cont coercive} satisfies
        \begin{align}
            \mathcal{C}_{2} \leqslant ||\mathcal{B}[u_{m}]||_{\mathcal{W}}^{2} / ||u_{m}||_{\mathcal{X}}^{2} \propto m^{2(s-3/2)}, \;\;\;\;\;\;\forall m \in \mathbb{N},
        \end{align}

        \noindent and hence, since $m \in \mathbb{N}$ can be chosen arbitrarily large, if $s<3/2$, then $\mathcal{C}_{2} = 0$. Conversely, if $s > 3/2$, then a similar argument shows that the constant $\mathcal{C}_{1}$ in Theorem \ref{thm:general weighted cont coercive} satisfies
        \begin{align}
            \mathcal{C}_{1} \geqslant ||\mathcal{B}[u_{m}]||_{\mathcal{W}}^{2} / ||u_{m}||_{\mathcal{X}}^{2} \propto m^{2(s-3/2)}, \;\;\;\;\;\;\forall m \in \mathbb{N},
        \end{align}

        \noindent in which case $\mathcal{C}_{1}$ must be unbounded.
        
        These arguments show that the well-posedness result of Theorem \ref{thm:general weighted cont coercive} holds for the Poisson equation only if $s=3/2$, which means that the space $\mathcal{W}$ is dictated by the variational problem and cannot be chosen arbitrarily by the practitioner as is often implied in the literature.

        In order to see the practical implications of using the $H^{s}(\partial\Omega)$, $s\neq 3/2$ norm, rather than the $H^{3/2}(\partial\Omega)$ norm, we return to the case where the true solution has the form $r^{m}\sin(m\theta)$ and we take advantage of \eqref{eq:norm calcs} to replace the $H^{s}$-norm on the boundary by a weighted $L^{2}$-norm, i.e. $||u_{m}||_{H^{s}(\partial\Omega)} = m^{s}||u_{m}||_{L^{2}(\partial\Omega)}$. This is equivalent to taking $\mathcal{W} = L^{2}(\partial\Omega)$ and penalty parameter $\delta = (m^{s})^{2} = m^{2s}$. Of course, this approach will not be possible for computing the $H^{s}(\partial\Omega)$ for general problems but will suffice for the current example.
	
    We first illustrate the scenario when $s=0$ and consider \eqref{eq:coercivity example} with $m=2,4,8,16$. Hyperparameters for this example are provided in Appendix \ref{app:examples}. Figure \ref{fig:poisson boundaryweights1} shows the error in various metrics for the Galerkin neural network method with formulation \eqref{eq:coercivity example}. For the $L^{2}(\Omega)$, $H^{1}(\Omega)$, and $H^{2}(\Omega)$ norms, we observe a delay in the onset of convergence as $m$ is increased. As $m$ increases, the boundary data becomes more oscillatory; the underpenalization of the boundary data for $s=0$ means that the resolution of the interior terms in \eqref{eq:coercivity example} is prioritized in the earlier Galerkin neural network iterations at the expense of satisfying the boundary condition, resulting in the behavior observed in Figure \ref{fig:poisson boundaryweights1}.
	\begin{figure}[t!]
		\centering
		\includegraphics[width=1.7in]{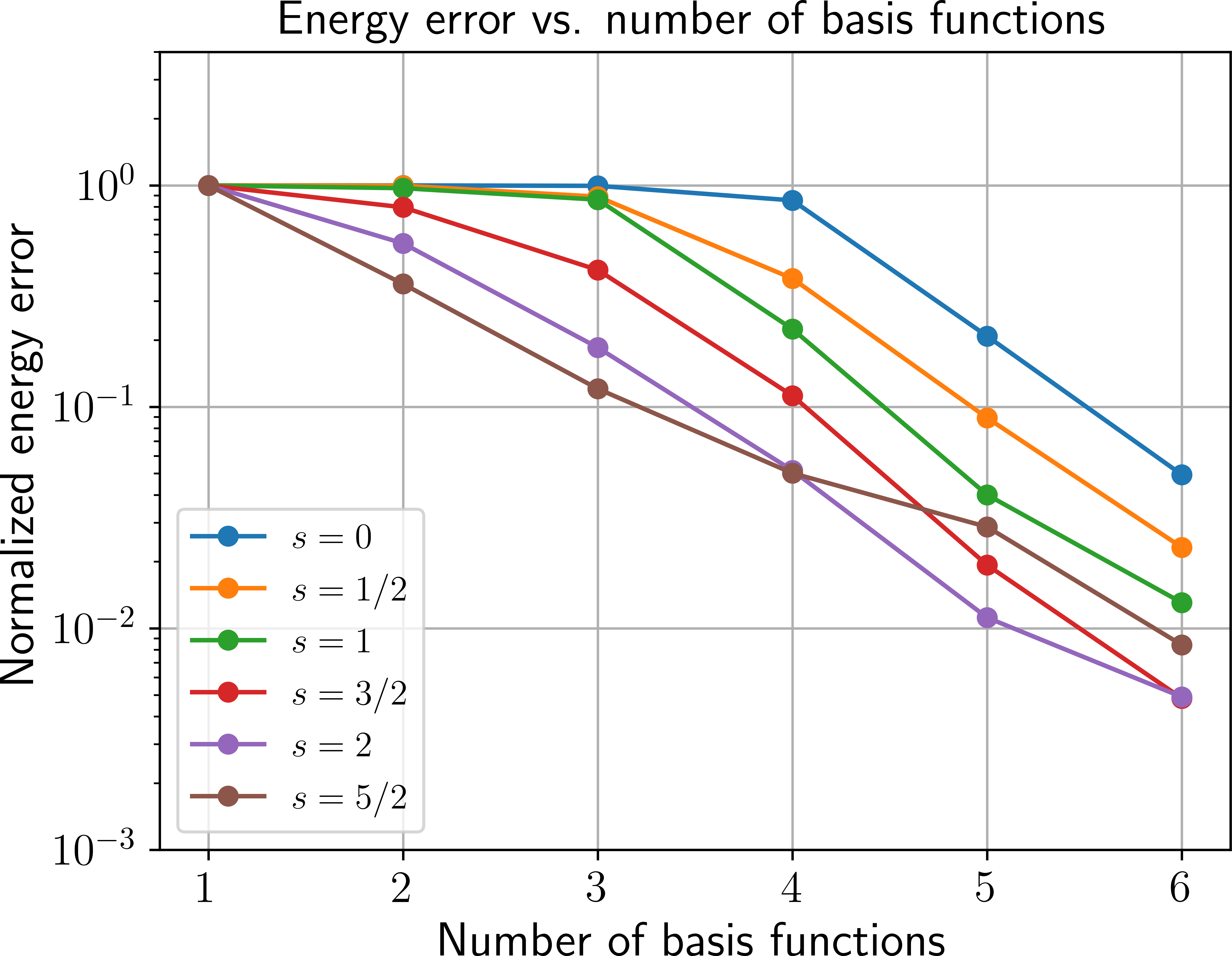}
		\quad
		\includegraphics[width=1.7in]{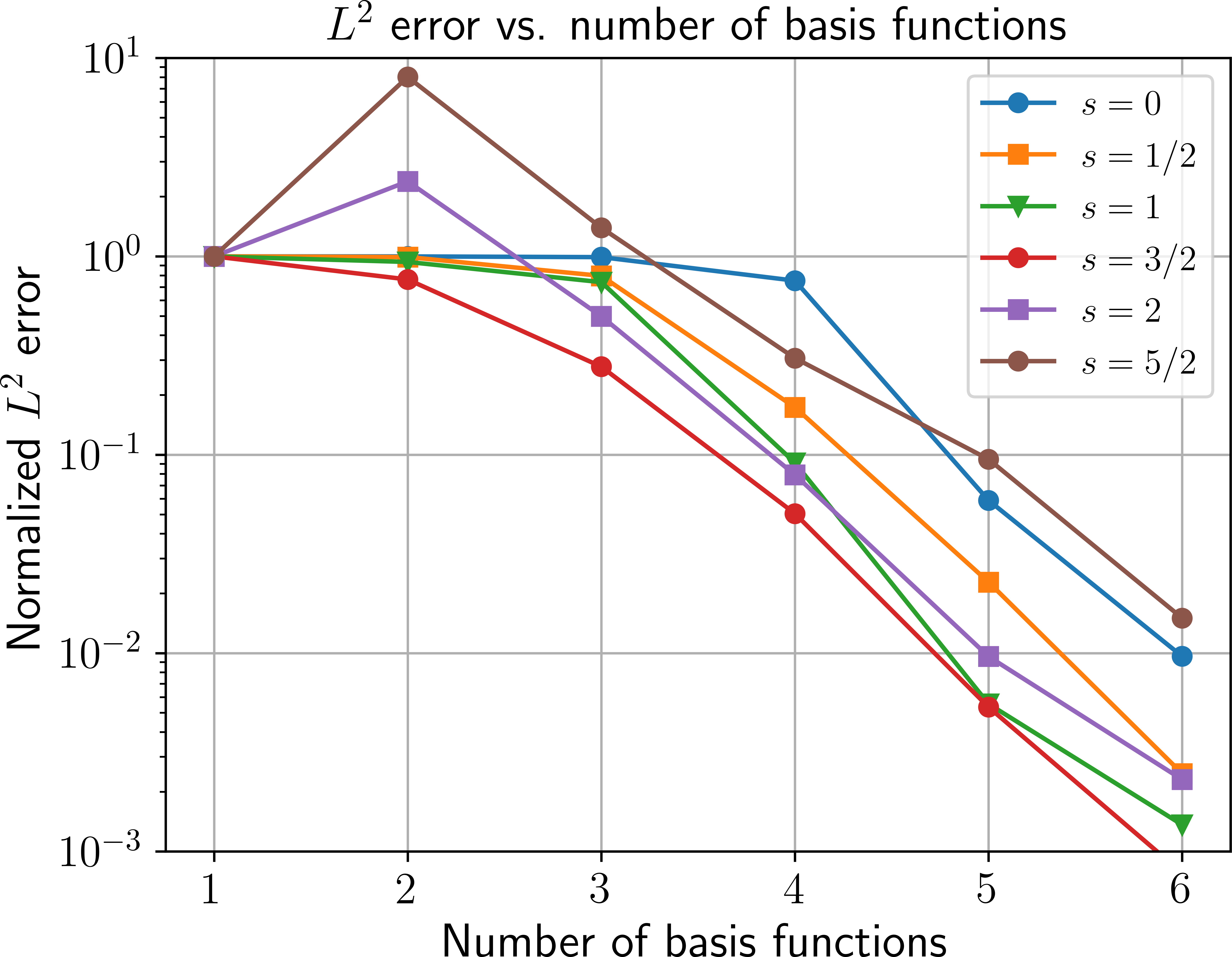}
		\quad
		
		\includegraphics[width=1.7in]{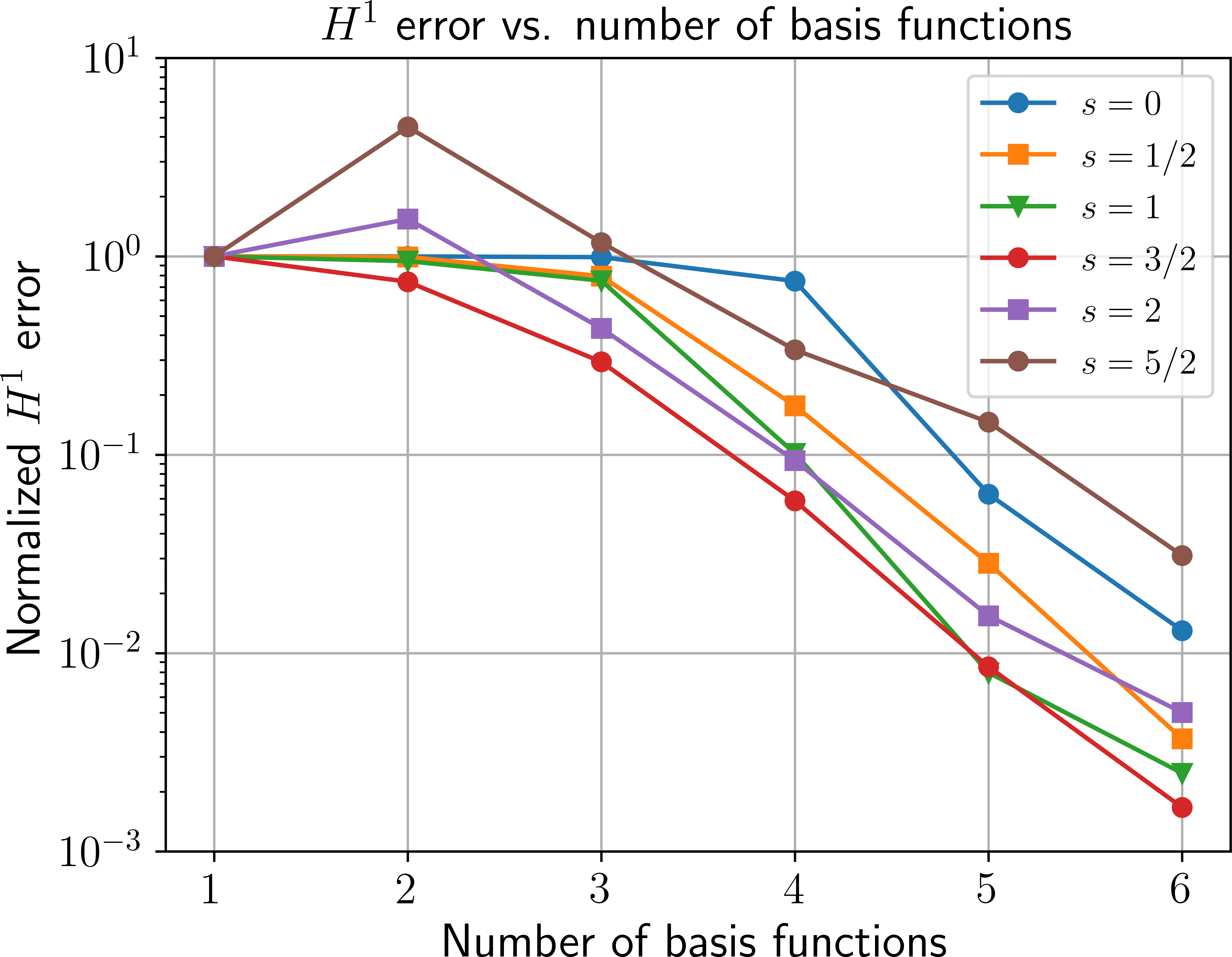}
		\quad
		\includegraphics[width=1.7in]{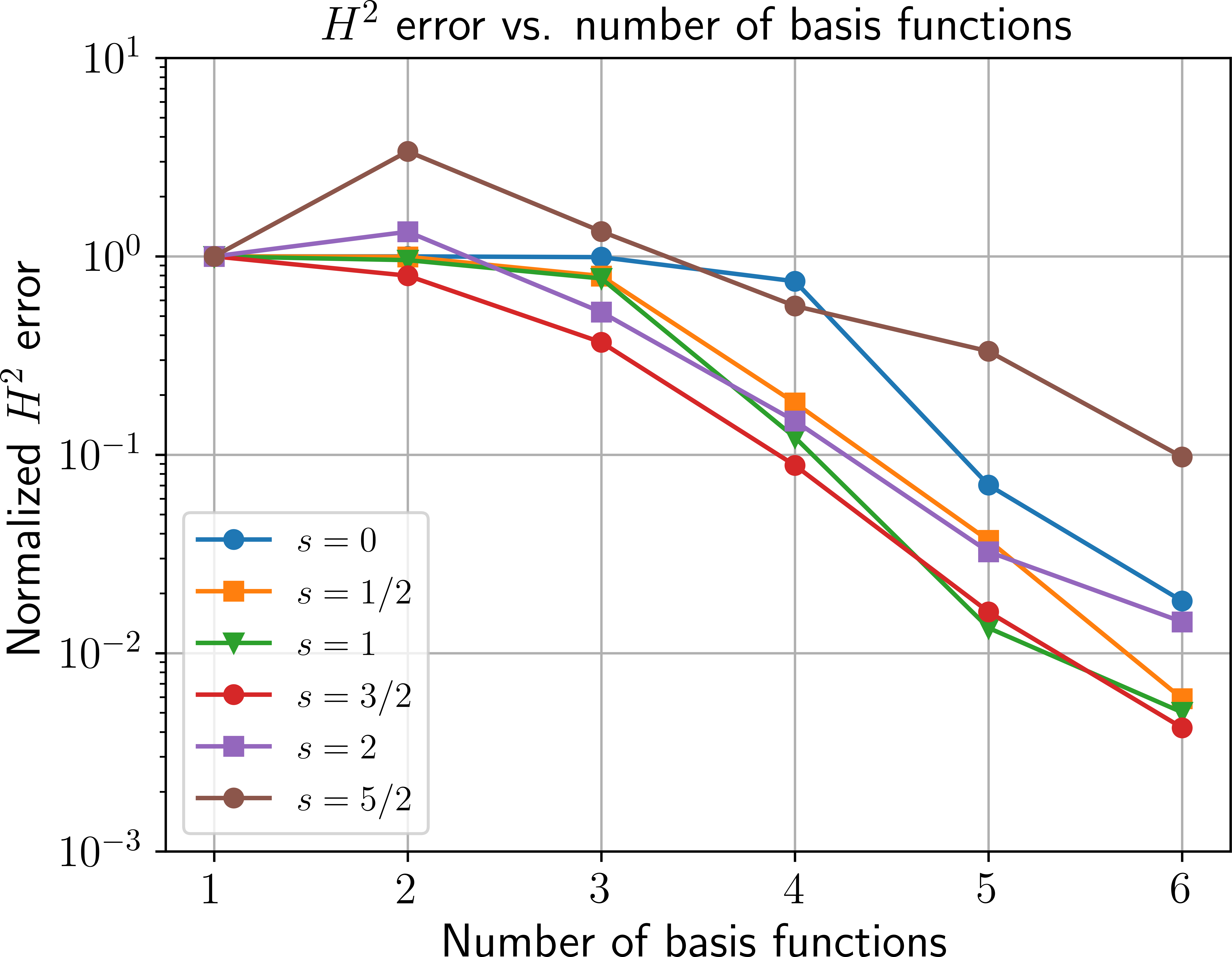}
		\caption{Errors with respect to the number of basis functions learned using the Galerkin neural network applied to the variational problem \eqref{eq:coercivity example} for several choices of penalty parameter $\delta = m^{2s}$ with $m=16$ and $s = 0,1/2,1,3/2,2,5/2$.}
		\label{fig:poisson boundaryweights2}
            \vspace{-4mm}
	\end{figure}
	
    In order to further amplify this point, Figure \ref{fig:poisson boundaryweights2} shows the results when we fix $m=16$ and consider the choice $s = 0, 1/2, 1, 3/2, 2, 5/2$. Regardless of whether the error is measured in the $L^{2}(\Omega)$, $H^{1}(\Omega)$, or $H^{2}(\Omega)$ norms, the greatest accuracy is obtained when $s=3/2$.
    

    The above example illustrates the importance of using the correct norms that arise in the well-posedness of the problem. In practice, the presence of fractional norms presents difficulties but is possible using techniques such as those in \cite{boundarypenalty}. However, generally speaking their computation is impractical and we instead use the Sobolev norm whose order is the integer part of the order of the desired Sobolev norm. For instance, the $H^{1}(\partial\Omega)$ penalization in Figure \ref{fig:poisson boundaryweights2} still performs quite well when compared to the optimal $H^{3/2}(\partial\Omega)$ penalization.
	\end{example}

	\subsection{Extended neural network architectures} \label{sec:knowledge-based fn}
            We return to the question of how to approximate the solution of variational problems of the form
        \begin{align} \label{eq:scalar singular lsq}
            u \in \mathcal{X} \;:\; a_{LS}(u,v) = F_{LS}(v) \;\;\;\forall v \in \mathcal{X},
        \end{align}

        \noindent where $a_{LS}$ and $F_{LS}$ are defined analogously to \eqref{eq:general lsq}. The focus of this section, in contrast to that of Section \ref{sec:lsq}, is on problems with known solution structures.
        
        Let $\Phi : \Omega \times \mathbb{R} \to \mathbb{R}$ be a given function and define $\mathcal{I}_{\Phi} := \{ \mu \in \mathbb{R} \;:\; \Phi(\cdot,\mu) \in \mathcal{X} \}$. For a given $m \in \mathbb{N}$, define
        \begin{align} \label{eq:singular NN}
            V_{m}^{\Phi} := \left\{ v \;:\; \Omega \to \mathbb{R} \;\bigg|\; v(x) = \sum_{j=1}^{m} d_{j}\Phi(x;\mu_{j}), \;\mathbf{d}\in\mathbb{R}^{m}, \;\mu_{j} \in \mathcal{I}_{\Phi}\right\},
        \end{align}
    
        \noindent and, for fixed $\mu \in \mathcal{I}_{\Phi}^{m}$, we also define the subset $V_{m}^{\Phi}(\mu)$ by
        \begin{align}
            V_{m}^{\Phi}(\mu) := \left\{ v \;:\; \Omega \to \mathbb{R} \;\bigg|\; v(x) = \sum_{j=1}^{m} d_{j}\Phi(x;\mu_{j}), \;\mathbf{d}\in\mathbb{R}^{m}\right\}.
        \end{align}
        
        \noindent The main difference between the sets $V_{m}^{\Phi}$ and $V_{m}^{\Phi}(\mu)$ is that in the former, the nonlinear parameters $\mu$ are allowed to vary over the course of training while in the latter, they are fixed. We shall first present the extended Galerkin neural network approach on the simpler space $V_{m}^{\Phi}(\mu)$ in which the parameter $\mu$ are assumed to be 
        known in advance before presenting the approach on the more complicated set $V_{m}^{\Phi}$ where $\mu$ is treated as a trainable parameter. As in the basic Galerkin neural network in Section \ref{sec:gnn}, we define the set of realizations whose parameters are bounded:
        \begin{align}
            \begin{aligned}
                V_{m}^{\Phi,C} :&= \{ v \in V_{m}^{\Phi} \;:\; ||(\mu,\mathbf{d})||_{\Lambda} := ||\mathbf{d}||_{\infty} + ||\mu||_{\infty} \leqslant C \},\\
                V_{m}^{\Phi,C}(\mu) :&= \{v \in V_{m}^{\Phi}(\mu) \;:\; ||(\mu,\mathbf{d})||_{\Lambda} := ||\mathbf{d}||_{\infty} + ||\mu||_{\infty} \leqslant C \}.
            \end{aligned}
        \end{align}
 
        The solutions of many classical PDEs exhibit well-known structures $\Psi(x;\lambda_{i})$, where $\lambda$ is a parameter that depends on the specific form of the domain, boundary conditions, and coefficients of the PDE. This means that for each fixed $\mathcal{M} \in \mathbb{N}$, the solution $u \in \mathcal{X}$ of such PDEs admits a representation of the form
        \begin{align} \label{eq:singular solution}
            u(x) = u_{\infty}(x) + u_{\Psi}(x;\lambda),
        \end{align}

        \noindent where $u_{\infty}$ is the portion of the solution with high Sobolev regularity and $u_{\Psi} \in V_{\mathcal{M}}^{\Psi}(\lambda)$ the portion with low regularity. For instance, for the Poisson equation, the singular part of the solution near corners of a polygonal domain is a sum of terms of the form $\Psi(r,\theta;\lambda_{j}) \propto r^{\lambda_{j}}\sin(\lambda_{j}\theta)$, where $\lambda_{j}$ depends on the interior angle at the corner and $(r,\theta)$ are polar coordinates centered at the corner. 
         \begin{figure}
            \centering
            \includegraphics[width=2.7in]{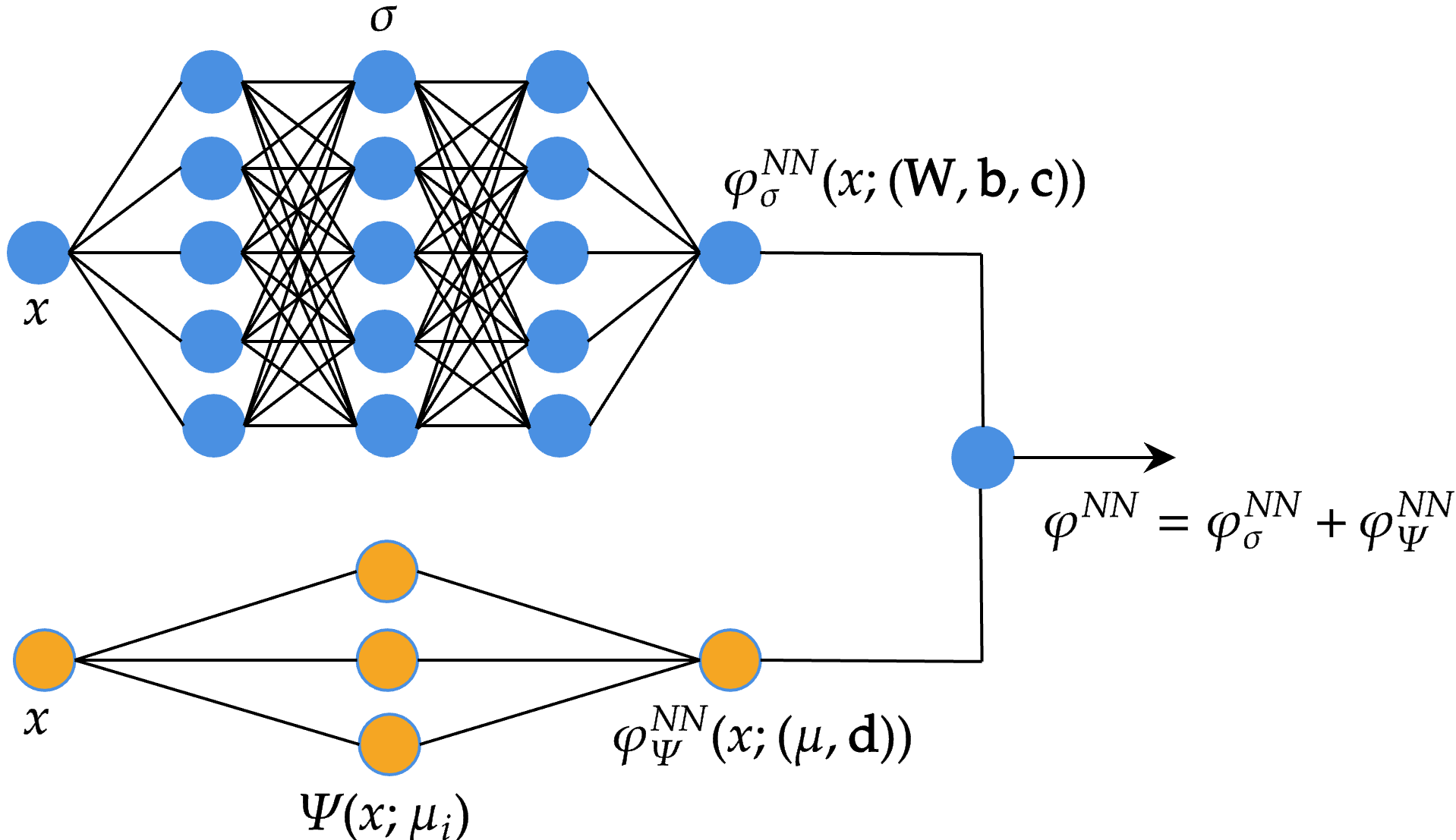}
            \caption{Extended neural network architecture for extended Galerkin neural networks. Each basis function $\varphi^{NN}$ consists of a high-regularity component $\varphi_{\sigma}^{NN}$ which is the realization of a standard feedforward neural network and a knowledge-based function $\varphi_{\Psi}^{NN}$ which utilizes known solution structures to augment the solution.}
            \label{fig:enriched GNN}
        \end{figure}   
        
        If such knowledge-based functions are available, then it makes sense for them to be incorporated into the Galerkin neural network framework and seek
    \begin{align} \label{eq:singular basis fn}
        \varphi_{i}^{NN} = \argmax_{\substack{v \in V_{\mathbf{n},L}^{\sigma, C_{1}} \bigoplus V_{\mathcal{M}}^{\Psi, C_{2}}(\lambda)\\ |||v||| = 1}} \langle r(u_{i-1}), v \rangle.
    \end{align}
    
    \noindent This means that the function $\varphi_{i}^{NN}$ has a structure that reflects that of the true solution in \eqref{eq:singular solution}, namely it is the sum of a neural network function $\varphi_{\sigma,i}^{NN} \in V_{\mathbf{n},L}^{\sigma, C_{1}}$ and a function $\varphi_{\Psi,i}^{NN} \in V_{\mathcal{M}}^{\Psi,C_{2}}(\lambda)$ belonging to the ``knowledge-based'' space. Given an initial approximation $u_{0} = u_{\sigma,0} + u_{\Psi,0}$ to \eqref{eq:scalar singular lsq} with smooth part $u_{\sigma,0}$ and low-regularity part $u_{\Psi,0} \in V_{\mathcal{M}}^{\Psi}$, we define the space
    \begin{align}
        S_{i}^{\Psi} := \text{span}\{u_{\sigma,0}/|||u_{0}|||, u_{\Psi,0}/|||u_{0}|||, \varphi_{\sigma,1}^{NN}, \varphi_{\Psi,1}^{NN}, \dots, \varphi_{\sigma,i}^{NN}, \varphi_{\Psi,i}^{NN}\},
    \end{align}

    \noindent and the approximation $u_{i}$ to the solution of \eqref{eq:scalar singular lsq} is defined to be the solution of the problem
    \begin{align} \label{eq:singular galerkin}
        u_{i} \in S_{i}^{\Psi} \;:\; a(u,v) = F(v) \;\;\;\forall v \in S_{i}^{\Psi}.
    \end{align}
    
    \noindent The solution may be written as a sum of smooth and singular parts $u_{i} = u_{\sigma,i} + u_{\Psi,i}$, where 
    \begin{align} \label{eq:u_sigma and u_Psi}
        u_{\sigma,i} \in \text{span}\{u_{\sigma,0}/|||u_{0}|||, \varphi_{\sigma,1}^{NN}, \dots, \varphi_{\sigma,i}^{NN}\}, \;\;\;u_{\Psi,i} \in \text{span}\{u_{\Psi,0}/|||u_{0}|||, \varphi_{\Psi,1}^{NN}, \dots, \varphi_{\Psi,i}^{NN}\}.
    \end{align}

    We briefly remark on the existence of the basis function $\varphi_{i}^{NN}$ defined as the maximizer of the weak residual in \eqref{eq:singular basis fn}. The existence of $\varphi_{i}^{NN}$ is guaranteed if $V_{\mathbf{n},L}^{\sigma,C_{1}} \oplus V_{\mathcal{M}}^{\Psi,C_{2}}(\lambda)$ is a compact subset of $\mathcal{X}$. The compactness of $V_{\mathbf{n},L}^{\sigma,C_{1}}$ is addressed in \cite[Proposition 3.5]{petersen2021topological} and discussed in the context of the basic Galerkin neural network approach in \cite[Remark 2.4]{gnn1}. To establish compactness of $V_{\mathcal{M}}^{\Psi}(\lambda)$, we utilize the following lemma.
    \begin{lemma} \label{lemma:compactness}
        Let
        \begin{align*}
            \Lambda := \left\{ \mu \in \mathcal{I}_{\Phi}^{m}, \;\mathbf{d} \in \mathbb{R}^{m} \;:\; ||(\mu,\mathbf{d})||_{\Lambda} := ||\mathbf{d}||_{\infty} + ||\mu||_{\infty} \leqslant C \right\}.
        \end{align*}
        
        \noindent Suppose the mapping $\mathcal{N} : (\Lambda, ||\cdot||_{\Lambda}) \to (V_{m}^{\Phi,C}, |||\cdot|||)$ defined by
        \begin{align*}
            \mathcal{N}(\mu, \mathbf{d}) = \sum_{i=1}^{m} d_{i}\Phi(\cdot;\mu_{i})
        \end{align*}
        is continuous. Then $V_{m}^{\Phi,C}$ is compact in $(\mathcal{X}, |||\cdot|||)$.
    \end{lemma}

    The set $V_{\mathcal{M}}^{\Psi,C_{2}}(\lambda)$ is always compact in $\mathcal{X}$ since the corresponding mapping $\mathcal{N}$ in Lemma \ref{lemma:compactness} is always continuous whenever the nonlinear parameter $\mu$ is fixed, as is the case for the space $V_{\mathcal{M}}^{\Psi,C_{2}}$. Thus, the maximization problem in \eqref{eq:singular basis fn} is well-defined. We postpone discussions of compactness of the more complicated set $V_{m}^{\Phi}$ until its use in Proposition \ref{thm:xgnn learning lambda} (see Remark \ref{remark:compactness}). 

    The effectiveness of the extended Galerkin neural network approach depends on the universal approximation properties of neural networks. The following result, proved in \cite{gnn1} for the basic Galerkin neural network in the case of single hidden layer networks, is directly applicable to problems whose solutions take the form \eqref{eq:singular solution}.
    \begin{proposition}[{\cite[Proposition 2.5]{gnn1}}] \label{prop:phi approximation}
        Let $0 < \varepsilon < 1$ be given and consider networks consisting of a single hidden layer, i.e. $\mathbf{n}^{(i)} = n$. Then there exist $n(\varepsilon, u_{\infty}+u_{\Psi}-u_{i-1}) \in \mathbb{N}$ and $C(\varepsilon, u_{\infty}+u_{\Psi}-u_{i-1}) > 0$ such that if $n \geqslant n(\varepsilon, u_{\infty}+u_{\Psi}-u_{i-1})$ and $C \geqslant C(\varepsilon, u_{\infty}+u_{\Psi}-u_{i-1})$, then $\varphi_{i}^{NN}$ given by \eqref{eq:basis fn} satisfies $|||\varphi_{i} - \varphi_{i}^{NN}||| \leqslant 2\varepsilon/(1-\varepsilon)$.
    \end{proposition}

    \noindent It will not have escaped the reader's attention that the width $n$ appearing in Proposition \ref{prop:phi approximation} depends on both the smooth and singular parts of the solution. 

    By way of contrast, the next result, Proposition \ref{thm:xgnn theory} applies to multilayer networks and shows that if the extended Galerkin neural network approach is used, then the width $\mathbf{n}$ now only depends on the smooth part of the solution, whereas if the knowledge-based functions are omitted, then the width $\mathbf{n}$ would also depend on the singular part of the solution (which we expect to be larger). This expectation is borne out in Example \ref{ex:poisson Lshaped}.
    \begin{proposition} \label{thm:xgnn theory}
        \noindent Suppose $\varepsilon>0$. Then given $u_{0} \in (\mathcal{X}, |||\cdot|||)$ and $L_{i} \in \mathbb{N}$, there exists $\mathbf{n}^{(i)} \in \mathbb{N}^{L_{i}}$ and $C_{1} > 0$ depending only on $\varepsilon$ and the error $u_{\infty} - u_{\sigma,i-1}$ in the approximation of the high-regularity solution component as well as $C_{2}>0$ depending only on $\varepsilon$ and the error $u_{\Psi} - u_{\Psi,i-1}$ in the approximation of the low-regularity solution component such that
        \begin{align} \label{eq:singular approximation result}
            |||\varphi_{i} - \varphi_{i}^{NN}||| \leqslant \varepsilon,
        \end{align}

        \noindent where $\varphi_{i}^{NN} := \varphi_{\sigma,i}^{NN} + \varphi_{\Psi,i}^{NN} \in V_{\mathbf{n}^{(i)},L_{i}}^{\sigma,C_{1}} \oplus V_{\mathcal{M}}^{\Psi}(\lambda)$ is defined by \eqref{eq:singular basis fn}.
    \end{proposition}

    \begin{proof}
        We partition the normalized approximation error $\varphi_{i} = (u - u_{i-1}) / |||u-u_{i-1}||| = \varphi_{\sigma,i} + \tilde{\varphi}_{\Psi}$, where $\varphi_{\sigma,i} := (u_{\infty} - u_{\sigma,i-1})/|||u-u_{i-1}|||$ and $\varphi_{\Psi} := (u_{\Psi} - u_{\Psi,i-1})/|||u-u_{i-1}||| \in V_{\mathcal{M}}^{\Psi}(\lambda)$.
        
        Thanks to \cite{hornik1990universal, kidger2020universal}, we may approximate $\varphi_{\sigma,i}$ by $\tilde{\varphi}_{\sigma} \in V_{\mathbf{n}, L_{i}}^{\sigma}$ such that
        \begin{align}
            |||\varphi_{\sigma,i} - \tilde{\varphi}_{\sigma}||| \leqslant \varepsilon/2
        \end{align} 
        provided that $\mathbf{n}$ is sufficiently large depending on $\varepsilon$ and $u_{\infty} - u_{\sigma,i-1}$. Thus, setting $\tilde{\varphi} := \tilde{\varphi}_{\sigma} + \varphi_{\Psi}$, we have the estimate
        \begin{align}
            |||\varphi_{i} - \tilde{\varphi}||| \leqslant \varepsilon/2.
        \end{align}
        \noindent Moreover, since $|||\varphi_{i}||| = 1$, we have $|||\tilde{\varphi}||| - 1 \in (-\varepsilon/2, \varepsilon/2)$.

        We next demonstrate that the normalization of $\tilde{\varphi}$, which we shall define by $\hat{\varphi} := \tilde{\varphi} / |||\tilde{\varphi}|||$, is also an adequate approximation to $\varphi_{i}$. By the triangle inequality, we have
        \begin{align} \label{eq:intermediate estimate}
            |||\varphi_{i} - \hat{\varphi}||| &\leqslant |||\varphi_{i} - \tilde{\varphi}|||  + |||\tilde{\varphi} - \frac{\tilde{\varphi}}{|||\tilde{\varphi}|||}|||\notag\\
            &\leqslant \varepsilon/2 + \frac{1}{|||\tilde{\varphi}|||} |||(|||\tilde{\varphi}||| - 1)\tilde{\varphi}||| \leqslant \varepsilon.
        \end{align}
        
        \noindent Choose $\mathbf{n}^{(i)}$ so that $\mathbf{n}^{(i)}_{j} \geqslant \mathbf{n}_{j}$ for $1 \leqslant j \leqslant L_{i}$, $C_{1} \geqslant ||(\mathbf{\hat{W}},\mathbf{\hat{b}},\mathbf{\hat{c}})||_{\mathcal{NN}}$, and $C_{2} \geqslant ||(\hat{\mu}, \mathbf{\hat{d}})||_{\Lambda}$ where $(\mathbf{\hat{W}}, \mathbf{\hat{b}}, \mathbf{\hat{c}})$ and $(\hat{\mu}, \mathbf{\hat{d}})$ are the parameters corresponding to the realizations $\tilde{\varphi}_{\sigma}/|||\tilde{\varphi}|||$ and $\varphi_{\Psi}$, respectively. 
        
        It remains to show that the maximizer $\varphi_{i}^{NN}$ is sufficiently close in norm $\hat{\varphi}$ and hence, $\varphi_{i}$. The argument does not differ from the corresponding proof for the basic Galerkin neural network framework in \cite{gnn1}. For completeness, we summarize the argument here. Since $\varphi_{i}^{NN}$ is the maximizer of $\langle r(u_{i-1}),\cdot \rangle = a(u-u_{i-1}, \cdot)$, we have $a(\varphi_{i},\varphi_{i}^{NN}) \geqslant a(\varphi_{i},v)$ for all $v \in V_{\mathbf{n}^{(i)},L_{i}}^{\sigma, C_{1}} \oplus V_{\mathcal{M}}^{\Psi, C_{2}}(\lambda)$ with $|||v|||=1$, and
        \begin{align*}
            |||\varphi_{i} - \varphi_{i}^{NN}|||^{2} &= |||\varphi_{i}|||^{2} - 2a(\varphi_{i}, \varphi_{i}^{NN}) + |||\varphi_{i}^{NN}|||^{2}\\
            &\leqslant |||\varphi_{i}|||^{2}  - 2a(\varphi_{i},\hat{\varphi}) + |||\hat{\varphi}|||^{2} = |||\varphi_{i} - \hat{\varphi}|||^{2} \leqslant \varepsilon^{2}.
        \end{align*}
    \end{proof}

    At first sight, one might expect that the knowledge-based functions need only be utilized once and for all in the first iteration, i.e. $i=1$. However, as the Galerkin iteration proceeds and the smooth part of the approximation encapsulated by the feedforward neural network becomes more refined, it may be necessary to adjust the linear coefficients $\mathbf{d}^{(i)}$ of the knowledge-based functions accordingly. For this reason, the singular functions are incorporated at every subsequent Galerkin iteration.

    Given Proposition \ref{thm:xgnn theory}, we may also obtain approximation results for $u_{i} \approx u$ and the \emph{a posteriori} estimator $\langle r(u_{i-1}),\varphi_{i}^{NN} \rangle \approx |||u-u_{i-1}|||$.
    \begin{corollary} \label{cor:approximation}
        The approximation $u_{i}$ given in \eqref{eq:singular galerkin} satisfies the estimate
        \begin{align} \label{eq:error estimate xgnn}
            |||u-u_{i}||| \leqslant |||u-u_{0}||| \cdot \min\left\{1, \varepsilon \right\}^{i}.
        \end{align}

        \noindent Moreover, if $\varepsilon < 1$, then each $\varphi_{i}^{NN}$ satisfies
         \begin{align} \label{eq:a posteriori}
            \langle r(u_{i-1}),\varphi_{i}^{NN} \rangle \leqslant |||u-u_{i-1}||| \leqslant \frac{1}{1-\varepsilon} \langle r(u_{i-1}),\varphi_{i}^{NN} \rangle.
        \end{align}
    \end{corollary}

    \begin{proof}
        The proof of \eqref{eq:error estimate xgnn} follows from \cite[Proposition 2.6]{gnn1}. The upper bound of \eqref{eq:a posteriori} follows from the fact that
        \begin{align*}
            \left|\langle r(u_{i-1}),\varphi_{i}^{NN}\rangle - |||u-u_{i-1}|||\right| \leqslant |||u-u_{i-1}|||\cdot |||\varphi_{i}^{NN} - \varphi_{i}||| = \varepsilon|||u-u_{i-1}|||
        \end{align*}

        \noindent while the lower bound follows from the observation that $\langle r(u_{i-1}), \varphi_{i}^{NN} \rangle := a(u-u_{i-1}, \varphi_{i}^{NN}) \leqslant |||u-u_{i-1}|||$. 
    \end{proof}
    
    In some cases, one might have knowledge of the structure of $\Psi$ while not having access to the exact values of $\lambda$. The parameters $\mu^{(i)}$ corresponding to $\varphi_{\Psi,i}^{NN}$ may then be treated as trainable nonlinear parameters and learned during the training step. One might also not know the structure of $\Psi$ in its entirety but instead have a suitable function $\Phi(\cdot;\lambda)$ which approximates $\Psi(\cdot,\lambda)$. In this case, a result analogous to Proposition \ref{thm:xgnn theory} may still be obtained. In what follows, we shall assume that $\Phi(\cdot,\mu) \in \mathcal{X}$, $\mu \in \mathcal{I}_{\Phi}$ is a family of functions such that for each $\varepsilon>0$, there exists $m_{*}(\varepsilon, u_{\Psi}) \in \mathbb{N}$ and $\tilde{\varphi}_{\Phi} \in V_{m_{*}}^{\Phi}$ such that
    \begin{align} \label{eq:singular estimate}
        |||u_{\Psi} - \tilde{\varphi}_{\Phi}||| \leqslant \varepsilon/2.
    \end{align}

    We propose seeking $\varphi_{\Psi,i}^{NN} \in V_{m_{*}}^{\Phi,C_{2}} \oplus (V_{m_{*}}^{\Phi,C_{2}}(\mu^{(i-1)}) \oplus \dots \oplus V_{m_{*}}^{\Phi,C_{2}}(\mu^{(0)}))$, where $\mu^{(1)}, \dots, \mu^{(i-1)}$ are the learned nonlinear hyperparameters corresponding to $\varphi_{\Psi,1}^{NN}, \dots, \varphi_{\Psi,i-1}^{NN}$ from the prior Galerkin iterations. The basis function $\varphi_{i}^{NN}$ is then obtained according to
    \begin{align} \label{eq:singular basis fn trained}
        \varphi_{i}^{NN} = \argmax_{\substack{v \in V_{\mathbf{n}^{(i)}}^{\sigma,C_{1}} \oplus \left( V_{m_{*}}^{\Phi,C_{2}} \oplus V_{m_{*}}^{\Phi,C_{2}}(\mu^{(i-1)}) \oplus \dots \oplus V_{m_{*}}^{\Phi,C_{2}}(\mu^{(0)}) \right)\\ |||v|||=1}} \langle r(u_{i-1}),v \rangle.
    \end{align}
    Note that the set functions $V_{m_{*}}^{\Phi, C_{2}} \oplus (V_{m_{*}}^{\Phi, C_{2}}(\mu^{(i-1)}) \oplus \dots \oplus V_{m_{*}}^{\Phi,C_{2}}(\mu^{(0)}))$ from which the knowledge-based part of the basis function is obtained has width $i\cdot m_{*}$. This is in contrast to \eqref{eq:singular basis fn} where the width of the set of knowledge-based functions was fixed and did not depend on $i$. The necessity of this growth is due to the nonlinear nature of the set $V_{m_{*}}^{\Phi,C}$. In particular, to correct the linear coefficients $\mathbf{d}^{(i)}$ corresponding to the singular part of the previous basis functions as the solution is refined, the knowledge-based functions $\Phi(\cdot;\mu^{(i)})$ must be utilized at every Galerkin iteration. The following result demonstrates that it is possible to obtain an analogous approximation result to Proposition \ref{thm:xgnn theory} when using approximate knowledge-based functions and inexact singular hyperparameters.
    \begin{remark} \label{remark:compactness}
        Compactness of the set $V_{m}^{\Phi,C}$ depends on the precise structure of $\Phi$ and cannot be proven for general $\Phi$. We briefly note, however, that the special choice $\Phi(x;\mu) = r(x)^{\mu}\sin(\mu\theta(x))$, where $(r,\theta)$ are the standard polar coordinate, results in $V_{m}^{\Phi,C}$ being a compact subset of $\mathcal{X}$. This function commonly arises in the solution of PDEs with singular behavior in non-convex corners. It is not difficult to show that $||\Phi(\cdot;\mu_{1}) - \Phi(\cdot;\mu_{2})||_{H^{s}(\Omega)} \to 0$ as $\mu_{1} \to \mu_{2}$ provided that $\mu_{1},\mu_{2}>0$ are chosen appropriately.



    \end{remark}
    
    \begin{proposition} \label{thm:xgnn learning lambda}
        \noindent Suppose $\varepsilon>0$ and there exists $\tilde{\varphi}_{\Phi} \in V_{m_{*}}^{\Phi}$ for some sufficiently large $m_{*} \in \mathbb{N}$ such that \eqref{eq:singular estimate} holds. Then given $u_{0} \in (\mathcal{X}, |||\cdot|||)$ and $L_{i} \in \mathbb{N}$, there exists $\mathbf{n}^{(i)} \in \mathbb{N}^{L_{i}}$ and $C_{1} > 0$ depending only on $\varepsilon$ and the error $u_{\infty} - u_{\sigma,i-1}$ in the approximation of the high-regularity component as well as $C_{2}>0$ depending only on $\varepsilon$ and the error $u_{\Psi} - u_{\Psi,i-1}$ low-regularity component such that
        \begin{align} 
            |||\varphi_{i} - \varphi_{i}^{NN}||| \leqslant 2\varepsilon,
        \end{align}
        
        \noindent where $\varphi_{i}^{NN}$ is defined by \eqref{eq:singular basis fn trained}.
    \end{proposition}
    
    \begin{proof}
        The two key differences between this scenario and the one presented in Proposition \ref{thm:xgnn theory} are (1) the set $V_{m}^{\Phi}$ from which we seek the basis functions is nonlinear since $\lambda$ is a trainable parameter, meaning that the widths $m^{(i)}$ must grow with each Galerkin iteration in order to allow for adjustments to the linear coefficients $\mathbf{d}^{(i)}$ as the smooth part of the approximation is refined (2) the approximation error includes additional error introduced by the inexact function $\Phi$. 
    
        We need only prove the existence of $\mathbf{n}^{(i)}$, $C_{1}$, $C_{2}$, and $\hat{\varphi} \in V_{\mathbf{n}^{(i)}, L_{i}}^{\sigma, C_{1}} \oplus V_{m_{*}}^{\Phi,C_{2}} \oplus V_{m_{*}}^{\Phi, C_{2}}(\mu^{(i-1)}) \oplus \dots \oplus V_{m_{*}}^{\Phi,C_{2}}(\mu^{(0)}))$ with $|||\hat{\varphi}|||=1$ such that $|||\hat{\varphi} - \varphi_{i}||| \leqslant 2\varepsilon$. The remainder of the proof does not differ from that of Proposition \ref{thm:xgnn theory}. We again set $\hat{\varphi} = \tilde{\varphi}/|||\tilde{\varphi}|||$ where $\tilde{\varphi} := \tilde{\varphi}_{\sigma} + \tilde{\varphi}_{\Psi}$, $\tilde{\varphi}_{\sigma} \in V_{\mathbf{n}, L_{i}}^{\sigma}$ is defined the same way as in the proof of Proposition \ref{thm:xgnn theory}, and the existence of $\tilde{\varphi}_{\Psi}$ is established as follows. Recall that $\varphi_{\Psi} \propto u_{\Psi} - u_{\Psi,i-1}$. Due to \eqref{eq:singular estimate}, a function from the set $V_{m_{*}}^{\Phi}$ is sufficient to approximate $u_{\Psi}$ to within $\varepsilon/2$. Likewise, due to \eqref{eq:u_sigma and u_Psi} we have $u_{\Psi,i-1} \in V_{m_{*}}^{\Phi}(\mu^{(0)}) \oplus \dots \oplus V_{m_{*}}^{\Phi}(\mu^{(i-1)})$. Thus, there exists $\tilde{\varphi}_{\Psi} \in V_{m_{*}}^{\Phi} \oplus (V_{m_{*}}^{\Phi}(\mu^{(i-1)}) \oplus \dots \oplus V_{m_{*}}^{\Phi}(\mu^{(0)}))$ such that
        \begin{align}
            |||\varphi_{\Psi} - \tilde{\varphi}_{\Psi}||| \leqslant \varepsilon/2
        \end{align}
        and we have
        \begin{align}
            \begin{aligned}
                |||\varphi_{i} - \hat{\varphi}||| 
                &\leqslant |||\varphi_{\sigma,i} - \tilde{\varphi}_{\sigma}||| + |||\varphi_{\Psi} - \tilde{\varphi}_{\Psi}||| + |||\tilde{\varphi} - \frac{\tilde{\varphi}}{|||\tilde{\varphi}|||}|||\\
                &\leqslant \varepsilon/2 + \varepsilon/2 + \varepsilon = 2\varepsilon.
            \end{aligned}
        \end{align}

        We choose $\mathbf{n}^{(i)}$ such that $\mathbf{n}^{(i)}_{j} \geqslant \mathbf{n}_{j}$, $1 \leqslant j \leqslant L_{i}$; $C_{1}$ such that $C_{1} \geqslant ||(\mathbf{\hat{W}}, \mathbf{\hat{b}}, \mathbf{\hat{c}})||_{\mathcal{NN}}$; and $C_{2}$ such that $C_{2} \geqslant ||(\hat{\mu}, \mathbf{\hat{d}})||_{\Lambda}$, where $(\mathbf{\hat{W}}, \mathbf{\hat{b}}, \mathbf{\hat{c}})$ and $(\hat{\mu}, \mathbf{\hat{d}})$ are the parameters corresponding to the realizations $\tilde{\varphi}_{\sigma}/|||\tilde{\varphi}|||$ and $\tilde{\varphi}_{\Psi}/|||\tilde{\varphi}|||$, respectively.
    \end{proof}

    An analogous result to Corollary \ref{cor:approximation} may be derived for the case described by \eqref{eq:singular basis fn trained}. The proof does not differ from that of Corollary \ref{cor:approximation}.

    Finally, we briefly discuss the training process for computing the maximizer in \eqref{eq:singular basis fn trained}. The nonlinear hyperparameters corresponding to the knowledge-based part of the basis function, $\varphi_{\Psi,i}^{NN}$, are given by $[\mu^{(i)}; \mu^{(i-1)}; \dots; \mu^{(0)}]$, where $\mu^{(i)}$ corresponding to the part of $\varphi_{\Psi,i}^{NN}$ in $V_{m_{*}}^{\Phi,C_{2}}$. Only $\mu^{(i)}$ should be trained since $\mu^{(i-1)}, \dots, \mu^{(0)}$ are fixed according to the previous Galerkin iterations. We update $\mu^{(i)}$ using a gradient-based optimizer; for convenience, we show a gradient descent update here:
    \begin{align} 
        \mu^{(i)} \leftarrow \mu^{(i)} + \nabla_{\mu^{(i)}}\left[ \frac{\langle r(u_{i-1}),v(\cdot;\mu^{(i)}) \rangle}{|||v(\cdot;\mu^{(i)})|||} \right].
    \end{align}
		
    \noindent The corresponding linear coefficients $\mathbf{c}$ and $\mathbf{d}$ are updated simultaneously by solving the following least squares system
    \begin{align} \label{eq:knowledge learn2}
        [\mathbf{c};\mathbf{d}] = \argmin_{\mathbf{c} \in \mathbb{R}^{n}, \mathbf{d} \in \mathbb{R}^{i\cdot m_{*}}} ||\mathbf{A}[\mathbf{c}; \mathbf{d}] - \mathbf{F}||_{\ell^{2}}
    \end{align}
		
    \noindent with the matrix $\mathbf{A}$ and vector $\mathbf{F}$ describing the projection of the error $u-u_{i-1}$ onto the space spanned by the extended neural network activation functions. For ease of notation, we only state the linear system for the case when $L^{(i)}=1$ and note that extension to the case of an arbitrary depth network is straightforward (see Section \ref{sec:gnn}):
    \begin{align}
        \begin{dcases}
            \mathbf{A} = \begin{bmatrix}
                \mathbf{A}^{\sigma\sigma} & \mathbf{A}^{\sigma\Psi}\\
                (\mathbf{A}^{\sigma\Psi})^{T} & \mathbf{A}^{\Psi\Psi}
            \end{bmatrix}, \;\;\;\mathbf{F} = \begin{bmatrix}
                \mathbf{F}^{\sigma}\\
                \mathbf{F}^{\Psi}
            \end{bmatrix}\\
            \mathbf{A}_{k\ell}^{\sigma\sigma} = a_{LS}(\sigma((\cdot)\cdot \mathbf{W}_{k} + \mathbf{b}_{k}), \sigma((\cdot)\cdot\mathbf{W}_{\ell} + \mathbf{b}_{\ell})), &1\leqslant k,\ell \leqslant n^{(i)}\\
            \mathbf{A}_{k\ell}^{\sigma\Psi} = a_{LS}(\sigma((\cdot)\cdot\mathbf{W}_{k} + \mathbf{b}_{k}), \Phi(\cdot; \mu^{(i)}_{\ell}), &1\leqslant k \leqslant n^{(i)}, \;\;\;1\leqslant \ell \leqslant m_{*}\\
            \mathbf{A}_{k\ell}^{\Psi\Psi} = a_{LS}(\Phi(\cdot; \mu^{(i)}_{k}), \Phi(\cdot; \mu^{(i)}_{\ell})), &1\leqslant k,\ell\leqslant m_{*}\\
            \mathbf{F}_{k}^{\sigma} = F_{LS}(\sigma((\cdot)\cdot\mathbf{W}_{k} + \mathbf{b}_{k})), &1\leqslant k\leqslant n^{(i)}\\
            \mathbf{F}_{k}^{\Psi} = F_{LS}(\Psi(\cdot; \mu^{(i)}_{k}), &1\leqslant k\leqslant m_{*}.
        \end{dcases}
    \end{align}

	\section{Applications} \label{sec:applications}
	
	We consider several indefinite and non-self adjoint boundary value problems and present (1) well-posed least squares variational formulations and (2) numerical results demonstrating the efficacy of the extended Galerkin neural network approach for each example. All hyperparameters for each example are provided in Appendix \ref{app:examples}. Of particular interest to us are problems whose solutions contain singular features, for instance due to low-regularity data or non-convex domain corners.
	
	As a motivating example, consider the incompressible stationary Stokes flow in velocity-pressure form given by the boundary value problem
	\begin{align} \label{eq:stokes strong}
		\begin{dcases}
			-\Delta\mathbf{u} + \nabla p = \mathbf{f} &\text{in}\;\Omega\\
			\text{div}\;\mathbf{u} = g &\text{in}\;\Omega\\
			\mathbf{u} = \mathbf{u}_{D} &\text{on}\;\partial\Omega,
		\end{dcases}
	\end{align}
	
	\noindent where $\Omega \subset \mathbb{R}^{2}$ is a bounded and connected polygonal domain, $\mathbf{u}$ is the fluid velocity, and $p$ the fluid pressure. To ensure solvability of \eqref{eq:stokes strong}, we require the compatibility condition $\int_{\Omega} g\;d\mathbf{x} = -\int_{\partial\Omega} \mathbf{u}_{D}\cdot\mathbf{n}\;ds$.
 
    The most natural variational formulation associated with \eqref{eq:stokes strong} is to seek $(\mathbf{u},p) \in \mathbf{H^{1}_{0}}(\Omega) \times (L^{2}(\Omega)/\mathbb{R})$ such that
	\begin{align} \label{eq:stokes natural}
		\begin{dcases}
			(\nabla\mathbf{u}, \nabla\mathbf{v})_{\Omega} - (p, \text{div}\;\mathbf{v})_{\Omega} &= (\mathbf{f},\mathbf{v})_{\Omega} + (\mathbf{u}_{D},\mathbf{v})_{\partial\Omega}\\
			(\text{div}\;\mathbf{u},q)_{\Omega} &= (g,q)_{\Omega}
		\end{dcases}
	\end{align}
	
	\noindent for all $(\mathbf{v},q) \in \mathbf{H^{1}_{0}}(\Omega) \times (L^{2}(\Omega)/\mathbb{R})$. Formulation \eqref{eq:stokes natural} is not SPD and therefore the basic Galerkin neural network formulation outlined in Section \ref{sec:gnn} is not applicable. Instead, one could attempt to use a least squares formulation based directly on \eqref{eq:stokes strong}. However, if $\Omega$ contains non-convex corners, then $\mathbf{u} \notin \mathbf{H^{2}}(\Omega)$ and $p \notin H^{1}(\Omega)$, and a pure $\mathbf{L}^{2}(\Omega)$-based least squares formulation would not be valid. The spaces $\mathcal{X}$, $\mathcal{V}$, and $\mathcal{W}$ must be chosen carefully to account for the reduced regularity. The same objection is valid for the Poisson equation. Fortunately, this issue can be resolved by choosing the spaces $\mathcal{X}$, $\mathcal{V}$, and $\mathcal{W}$ to be appropriately weighted Sobolev spaces. This approach is explained in the remainder of this section. For simplicity, we begin by considering the Poisson equation.
	
 	\subsection{Poisson equation} \label{sec:lsq poisson}
	
	We demonstrate the construction of a well-posed weighted least squares variational formulation for the Poisson equation:
	\begin{align} \label{eq:poisson2}
		\begin{dcases}
			-\Delta u = f &\text{in}\;\Omega\\
			u = g &\text{on}\;\partial\Omega.
		\end{dcases}
	\end{align}
	
	 \noindent We consider the case where $\Omega$ is a polygon with vertices $\{x^{(i)}\}_{i=1}^{M}$.
  
  	\subsubsection{Regularity and variational formulation}
    As mentioned in Section \ref{sec:xgnn}, the singular solution structure for \eqref{eq:poisson2} takes the form $\Psi(r,\theta) = \sum_{\lambda} c_{\lambda}r^{\lambda}\sin(\lambda\theta)$, where $\lambda > 1/2$. Similar to the Stokes problem, a pure $L^{2}$-based least squares formulation is not valid. However, the solution does belong to a weighted space $H^{2}_{\beta}(\Omega)$ provided that the weight $\beta>0$ is sufficiently large. 

    For $\beta = (\beta_{1},\dots,\beta_{M}) \geqslant \mathbf{0}$ and $s>0$, the spaces $H^{s}_{\beta}(\Omega)$ are weighted Sobolev spaces equipped with the norm
	\begin{align}
		||u||_{H^{s}_{\beta}(\Omega)}^{2} = ||u||_{H^{s-1}(\Omega)}^{2} + \sum_{|\alpha|=s} ||\prod_{i=1}^{M}r_{i}^{\beta_{i}}D^{\alpha}u||_{L^{2}(\Omega)}^{2},
	\end{align}
	
	\noindent where $r_{i}$ is the local radial coordinate centered at vertex $x^{(i)}$. If $s=0$, then the term involving the $H^{s-1}(\Omega)$ norm is dropped. In what follows, we will commonly use $(\cdot,\cdot)_{\beta,\Omega}$ to denote the $L^{2}_{\beta}(\Omega)$ inner product. The corresponding spaces $\mathbf{H}^{s}_{\beta}$ and $\mathbf{L}^{2}_{\beta}$ of vector-valued functions are defined in the usual way.

    A simple calculation shows that $r^{\lambda}\sin(\lambda\theta) \in H^{2}_{\beta}(\Omega)$ if and only if $\beta > 1-\lambda$. Thus, one might expect a continuous dependence result of the form \eqref{eq:a priori} to hold in the space $\mathcal{X} = H^{2}_{\beta}(\Omega)$, $\beta > 1-\lambda$. The following result, proved in \cite{kondratev}, confirms this expectation.
	\begin{figure}[t!]
		\centering
		\includegraphics[width=1.3in]{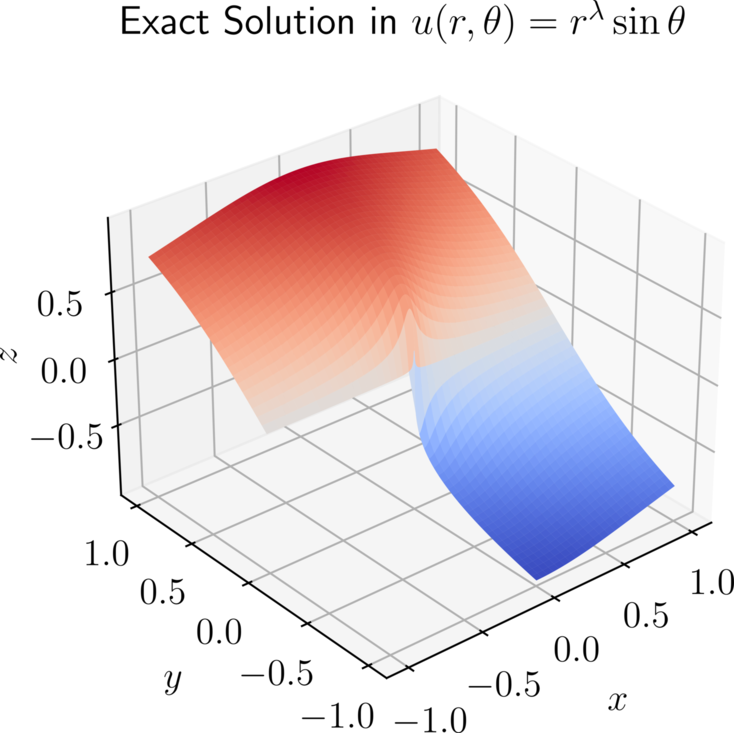}
		\quad
		\includegraphics[width=1.3in]{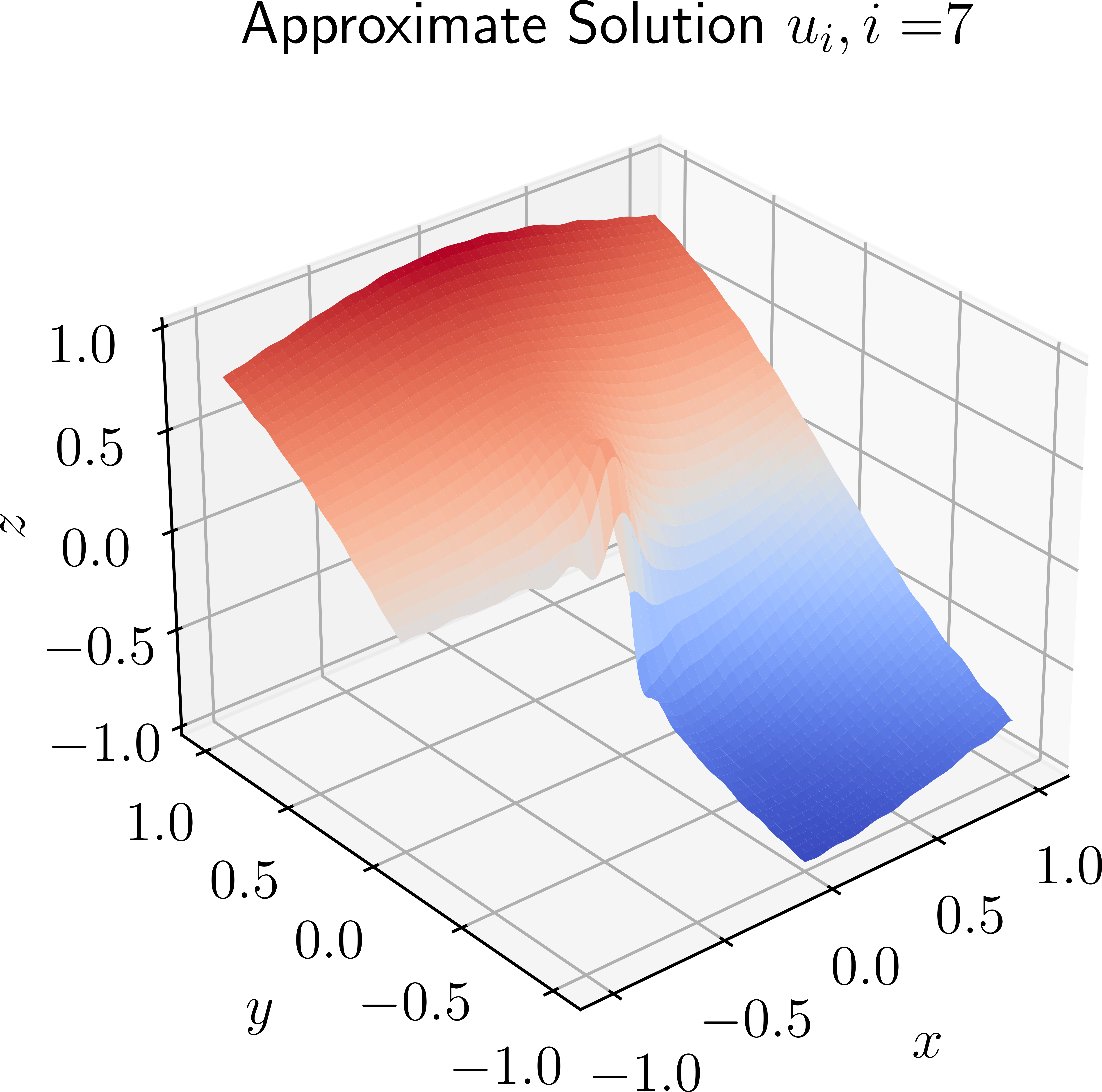}
		\quad
		\includegraphics[width=1.3in]{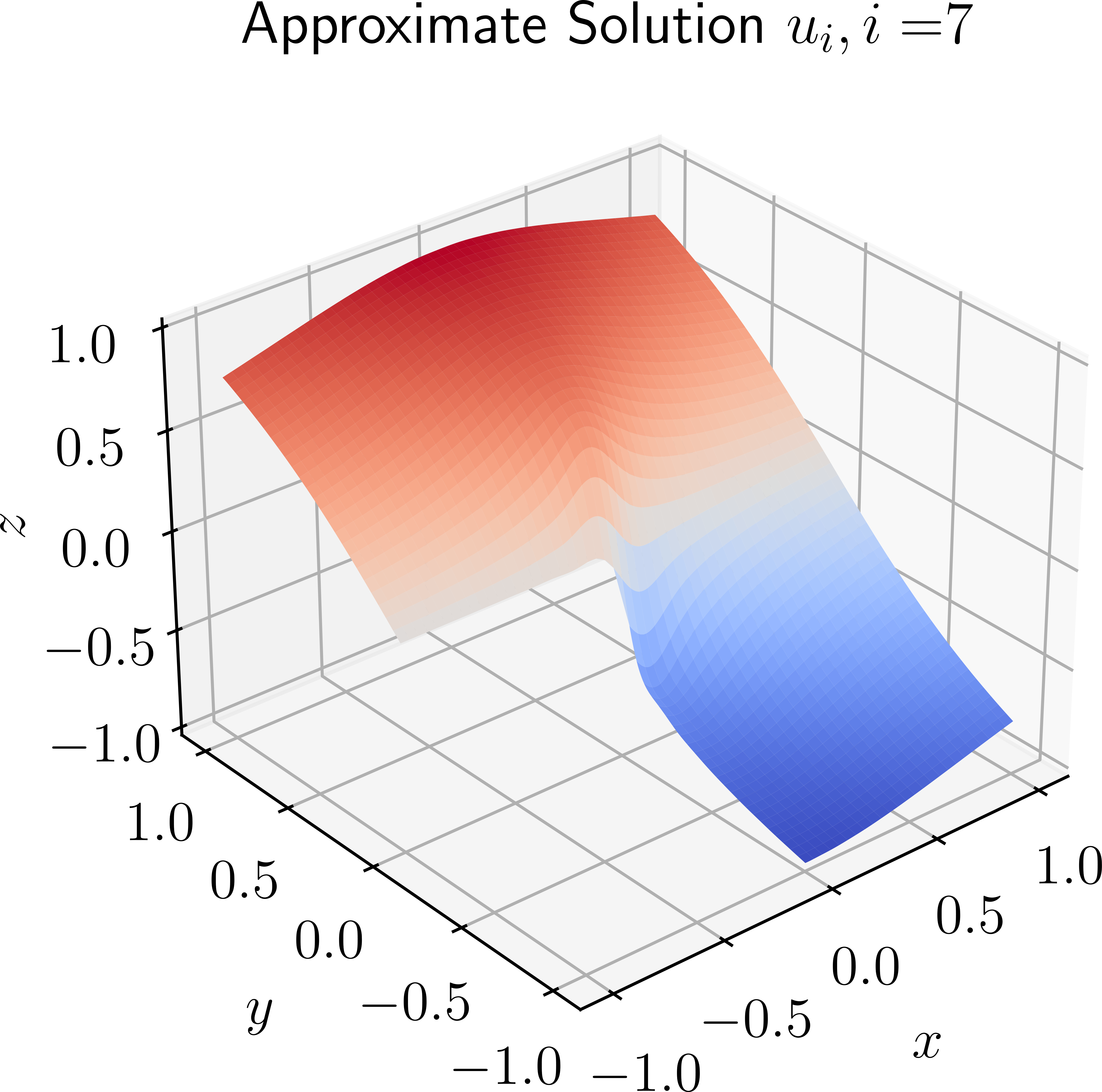}
		\caption{The true solution $u(r,\theta) = r^{\lambda}\sin{\theta}$ (left) and approximate solutions obtain using Galerkin neural networks applied to the variational problem \eqref{eq:poisson lsq} with $\beta = 0$ (middle) and $\beta=1$ (right).}
		\label{fig:poisson rlambda1}
	\end{figure}
	\begin{figure}[b!]
            \captionsetup{font={color=red}}
		\centering
		\includegraphics[width=1.3in]{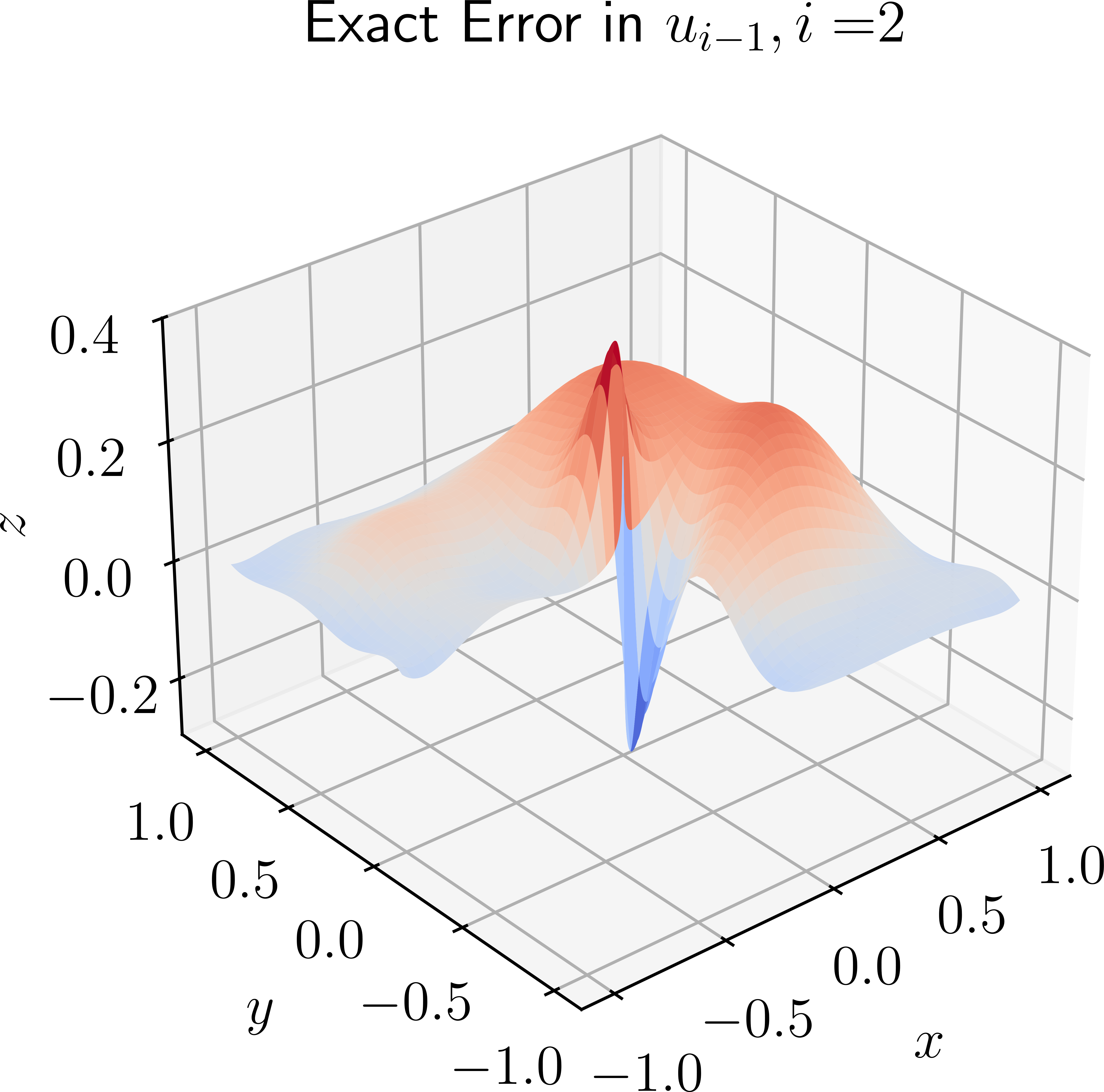}
		\quad
		\includegraphics[width=1.3in]{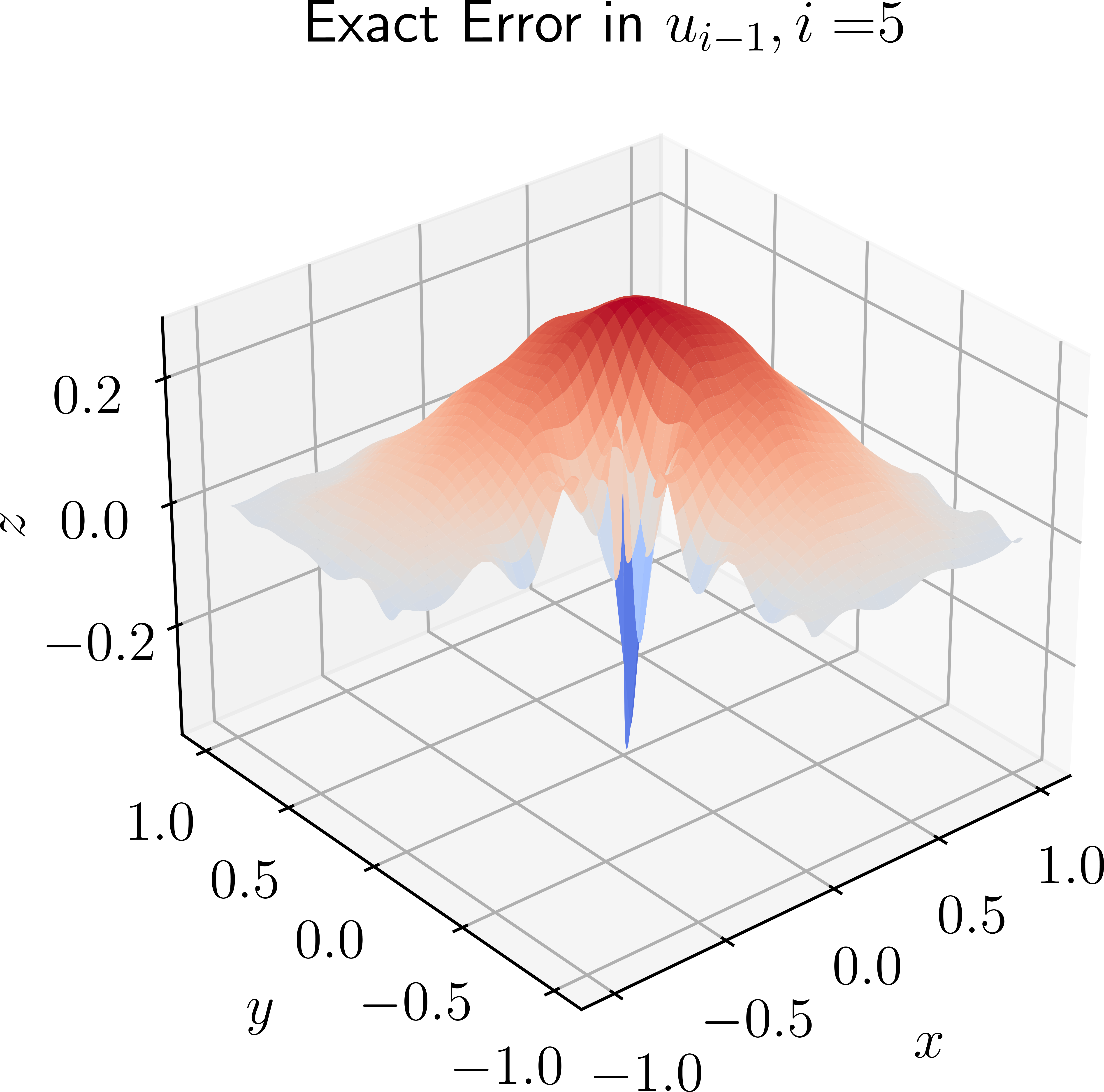}
		\quad
		\includegraphics[width=1.3in]{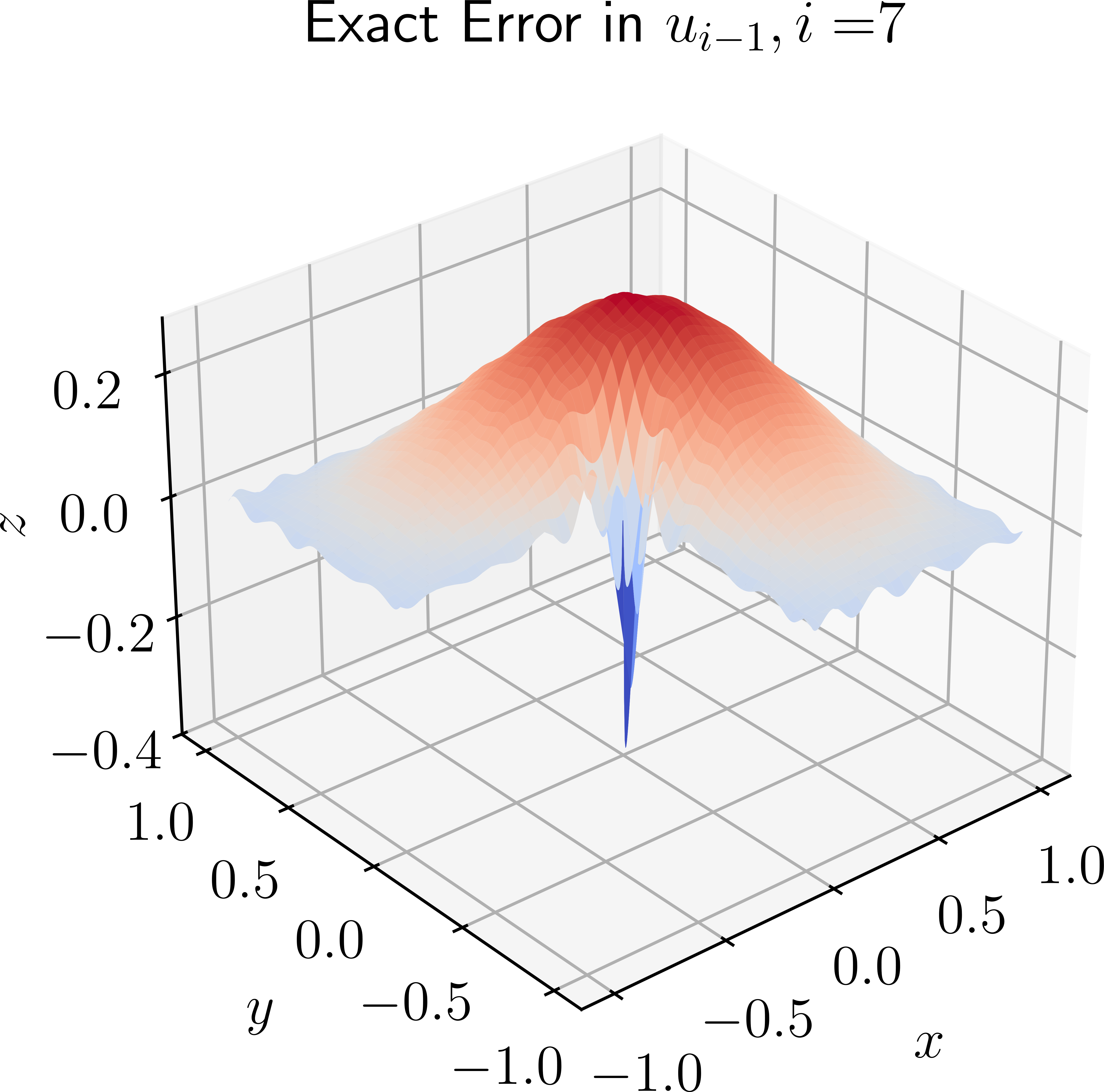}
		
		\includegraphics[width=1.3in]{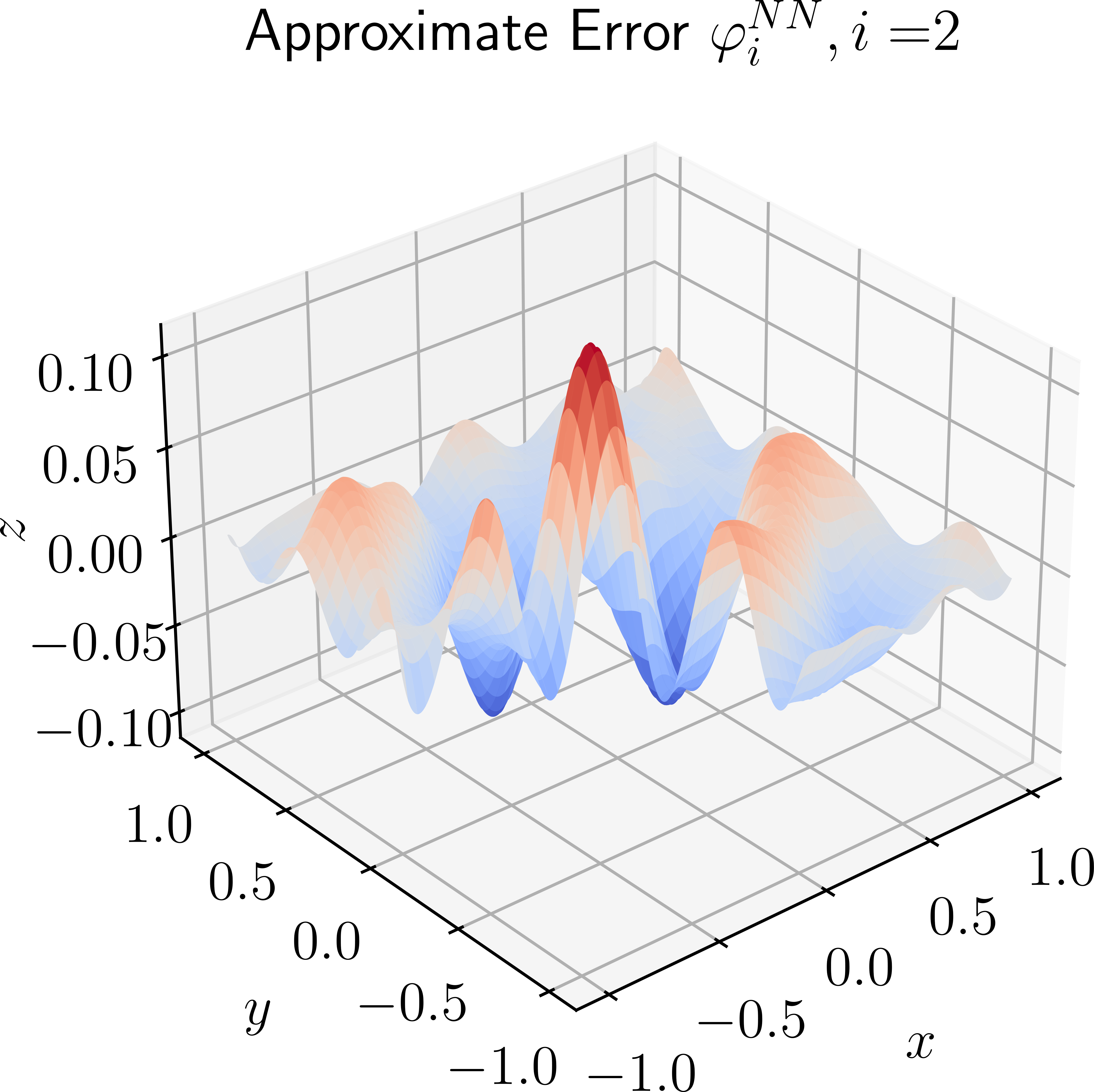}
		\quad
		\includegraphics[width=1.3in]{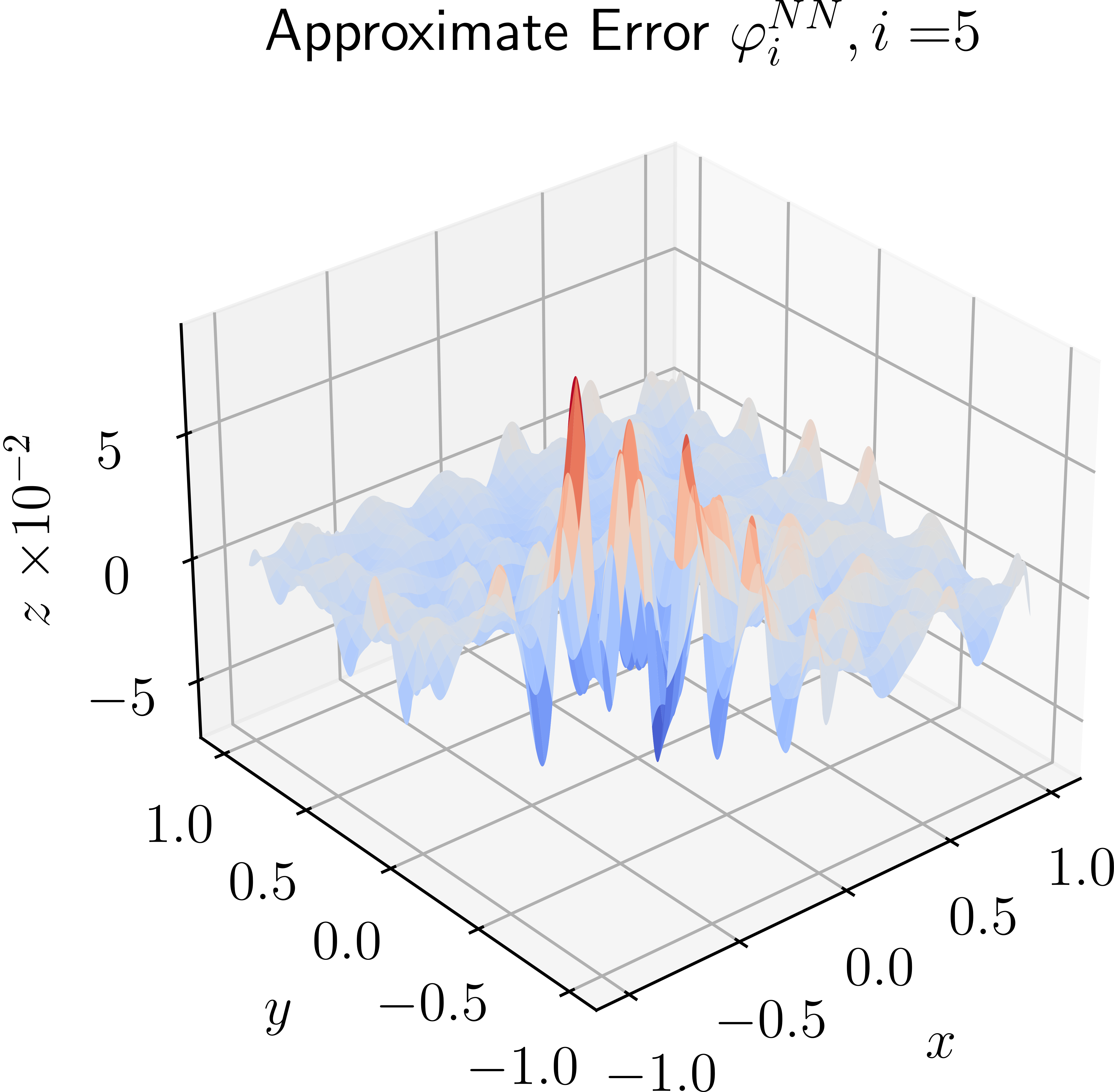}
		\quad
		\includegraphics[width=1.3in]{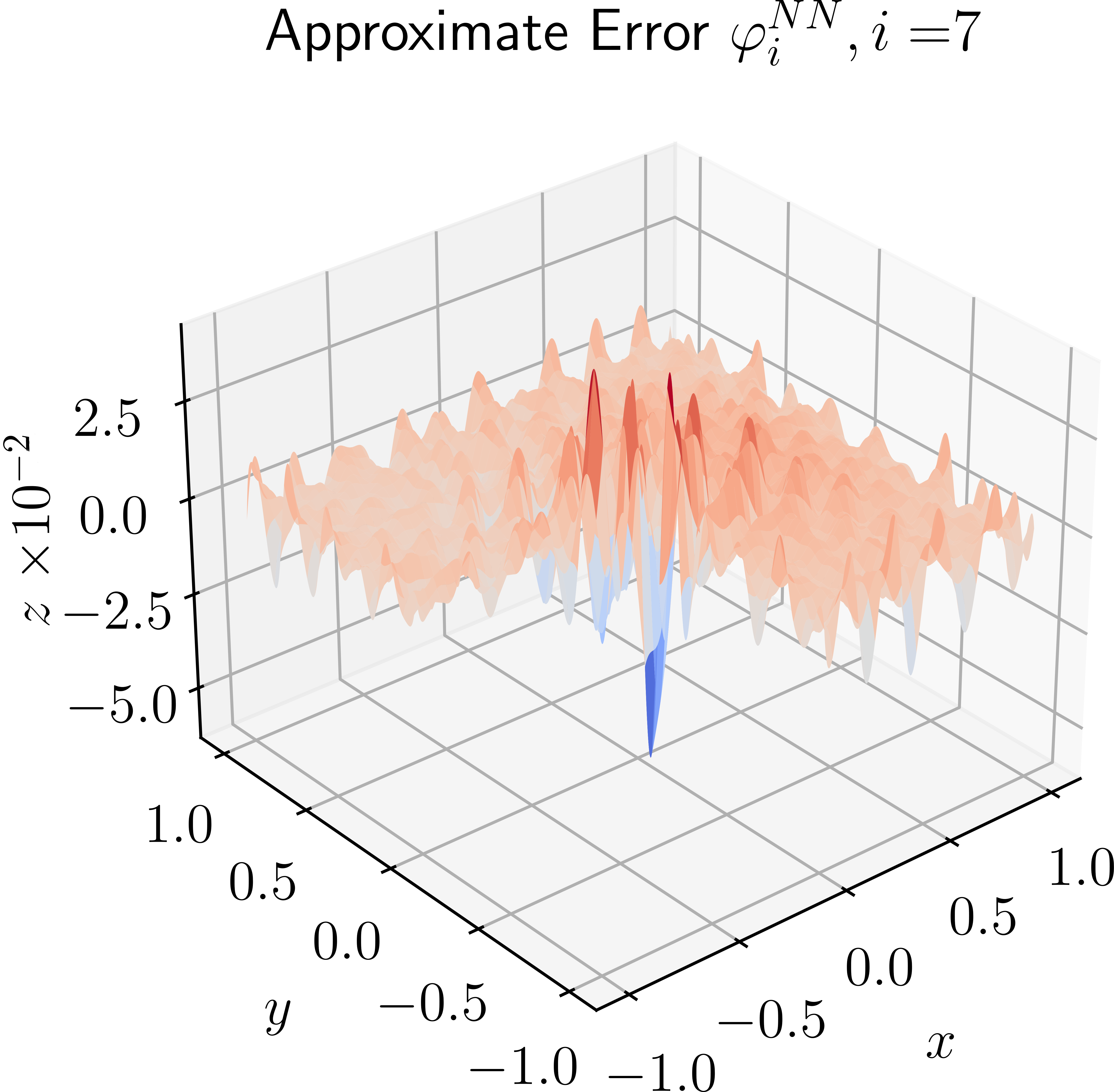}
		\caption{Example \ref{ex:poisson rlambda}: for $\beta=0$, the exact errors $u-u_{i-1}$ (top row) and approximate errors which correspond to basis functions $\varphi_{i}^{NN} \approx u-u_{i-1}$ learned using the Galerkin neural network approach (bottom row).}
		\label{fig:poisson rlambda2}
	\end{figure}

	\begin{theorem}
		\cite{kondratev} Suppose $f \in L^{2}_{\beta}(\Omega)$ and $g \in H^{3/2}_{\beta}(\partial\Omega)$. Then \eqref{eq:poisson2} admits a unique solution $u \in H^{2}_{\beta}(\Omega)$ and
		\begin{align*}
			||u||_{H^{2}_{\beta}(\Omega)} \leqslant C(||f||_{L^{2}_{\beta}(\Omega)} + ||g||_{H^{3/2}_{\beta}(\partial\Omega)}),
		\end{align*}
		
		\noindent for some $C>0$ independent of $u$, $f$, and $g$, where
		\begin{align*}
			||u||_{H^{3/2}_{\beta}(\partial\Omega)} := \inf_{U|_{\partial\Omega} = u} ||U||_{H^{2}_{\beta}(\Omega)}.
		\end{align*}
	\end{theorem}
		
	 We propose the weighted least squares variational formulation given by
	\begin{align}
		u \in H^{2}_{\beta}(\Omega) \;:\; (\Delta u, \Delta v)_{\beta,\Omega} &+ \delta(u,v)_{H_{\beta}^{3/2}(\partial\Omega)} = -(f,\Delta v)_{\beta,\Omega} + \delta(g,v)_{H_{\beta}^{3/2}(\partial\Omega)}\notag\\
		&\hspace{60mm}\forall v \in H^{2}_{\beta}(\Omega).\label{eq:poisson lsq}
	\end{align}
    
	\begin{figure}[t!]
            \captionsetup{font={color=red}}
		\centering
		\includegraphics[width=1.3in]{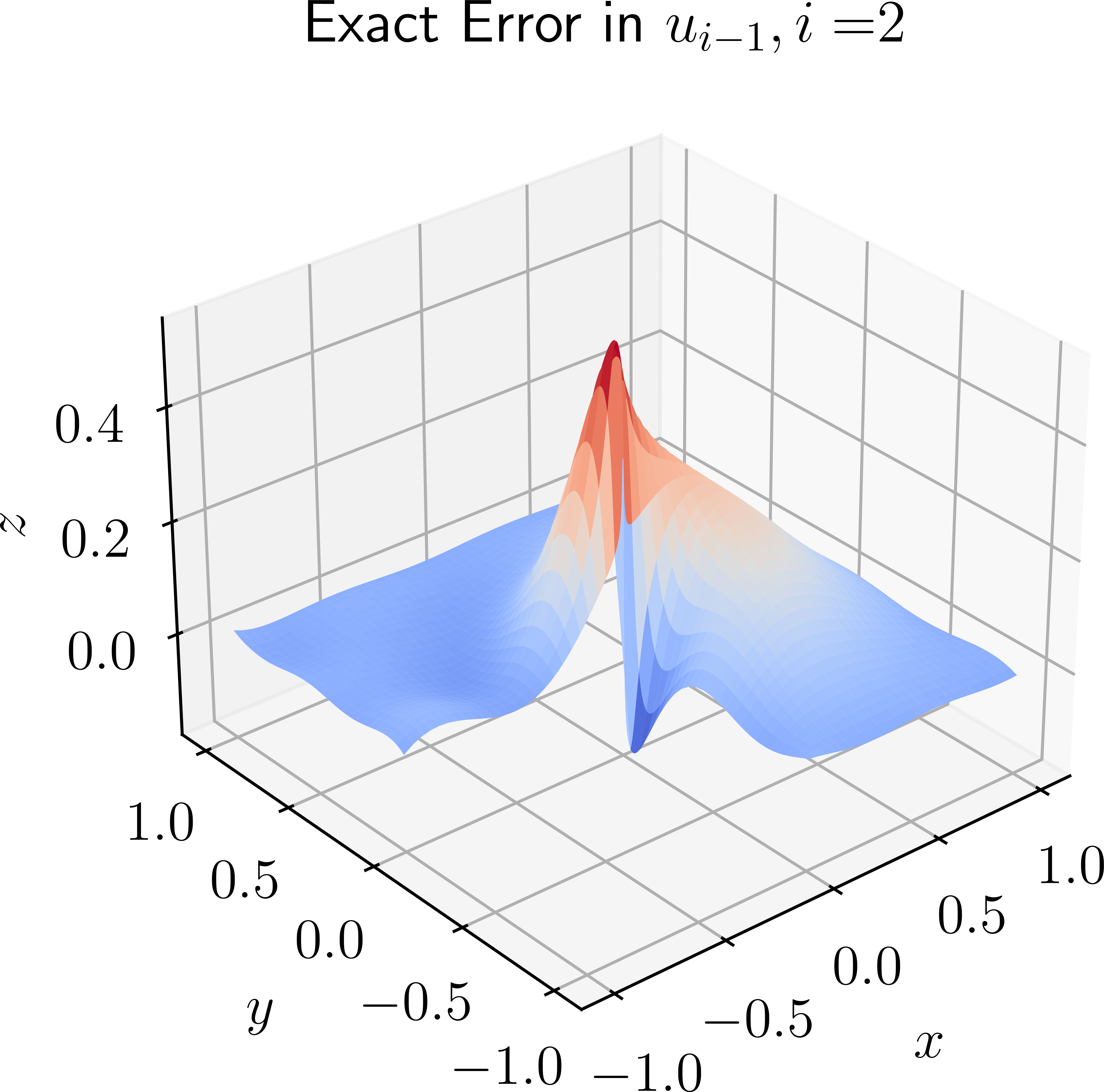}
		\quad
		\includegraphics[width=1.3in]{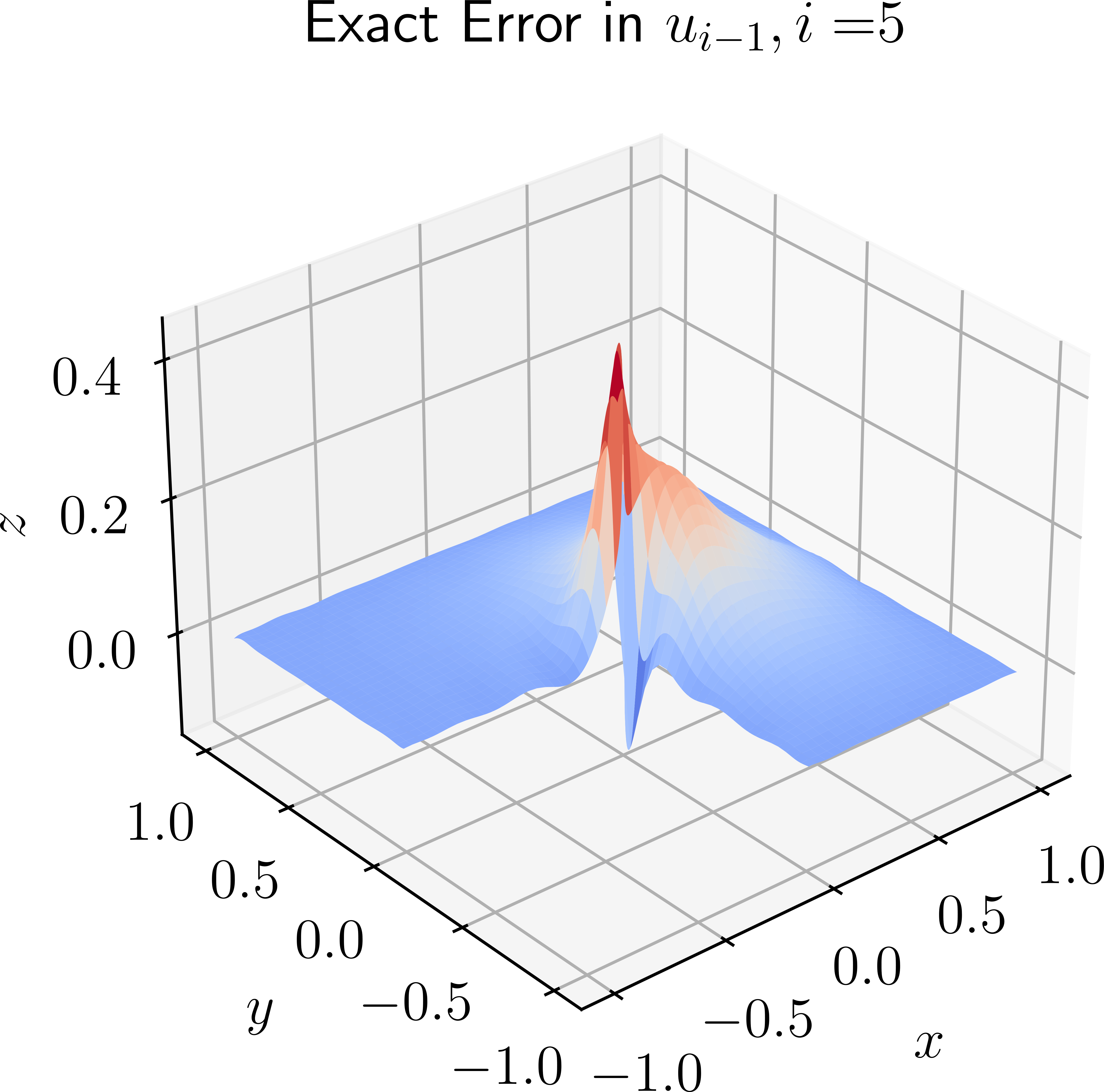}
		\quad
		\includegraphics[width=1.3in]{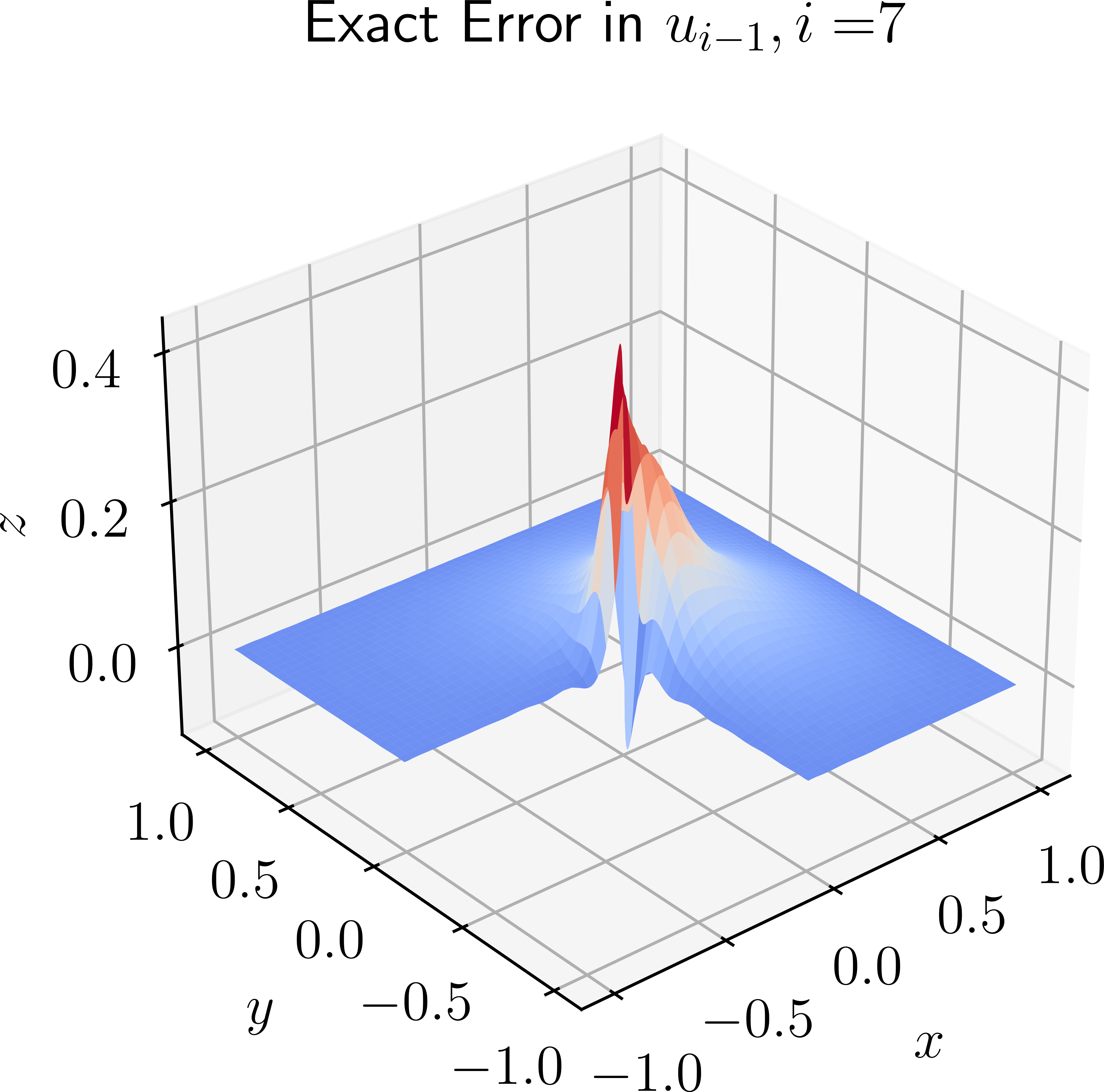}
		
		\includegraphics[width=1.3in]{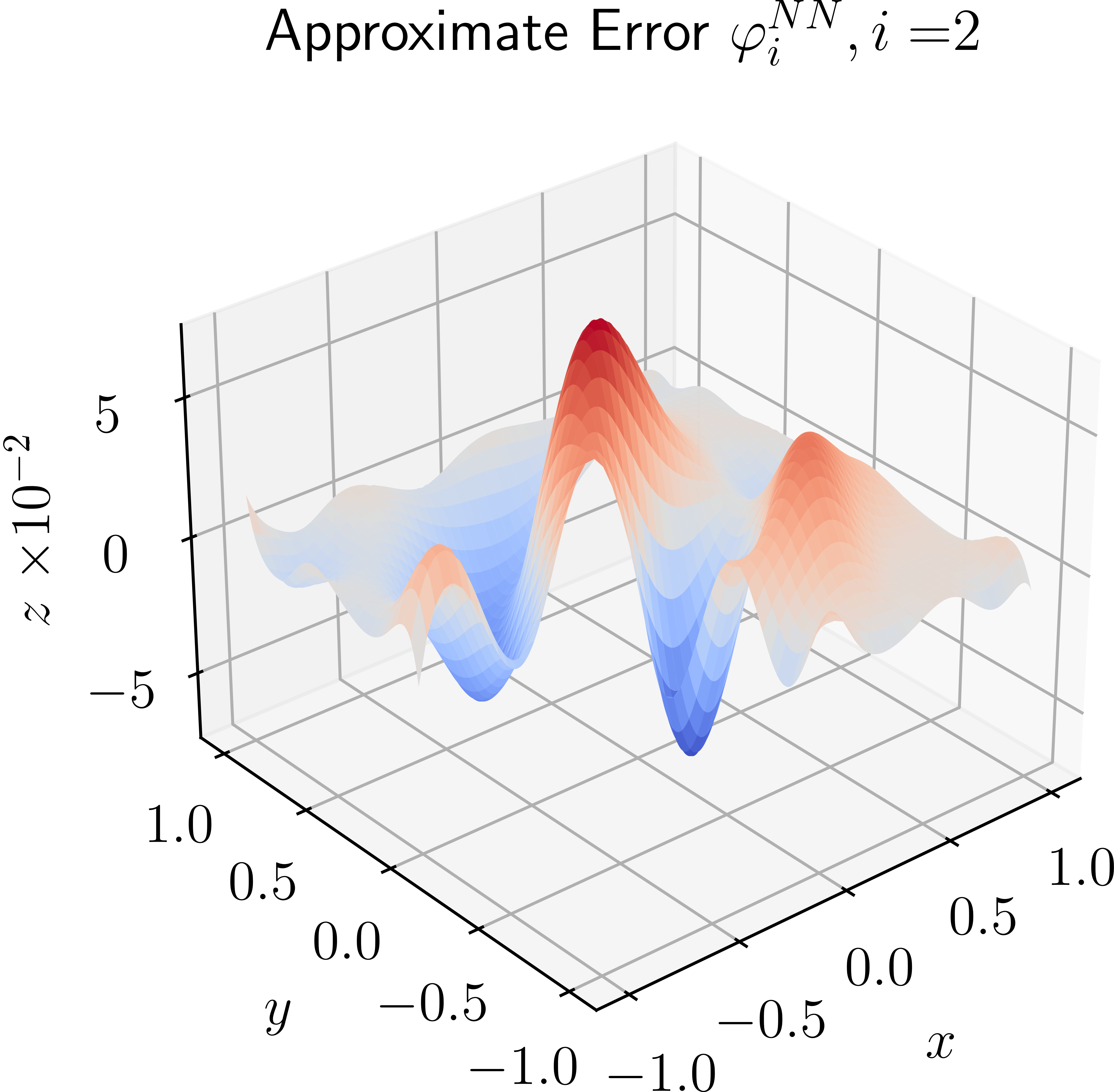}
		\quad
		\includegraphics[width=1.3in]{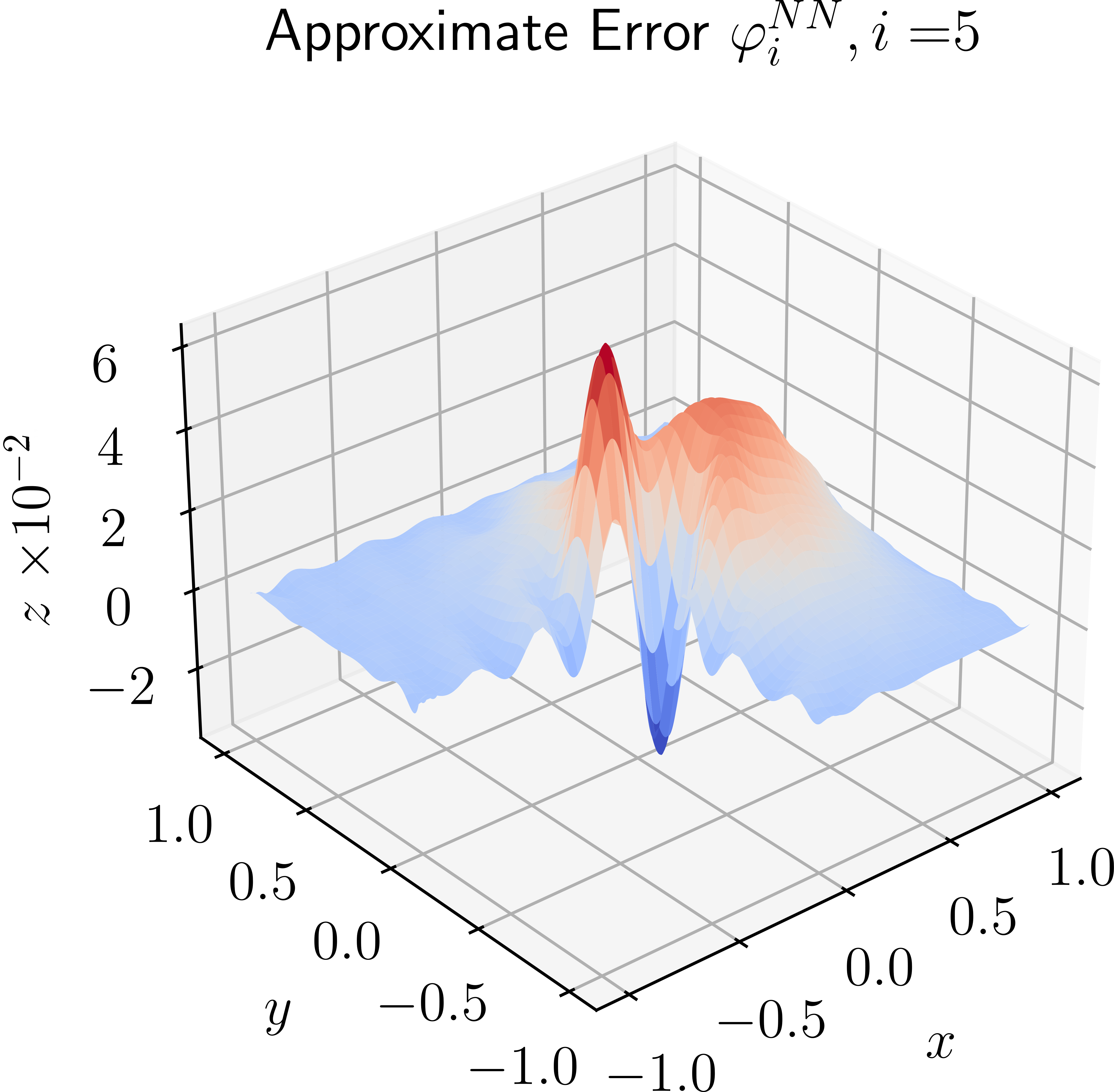}
		\quad
		\includegraphics[width=1.3in]{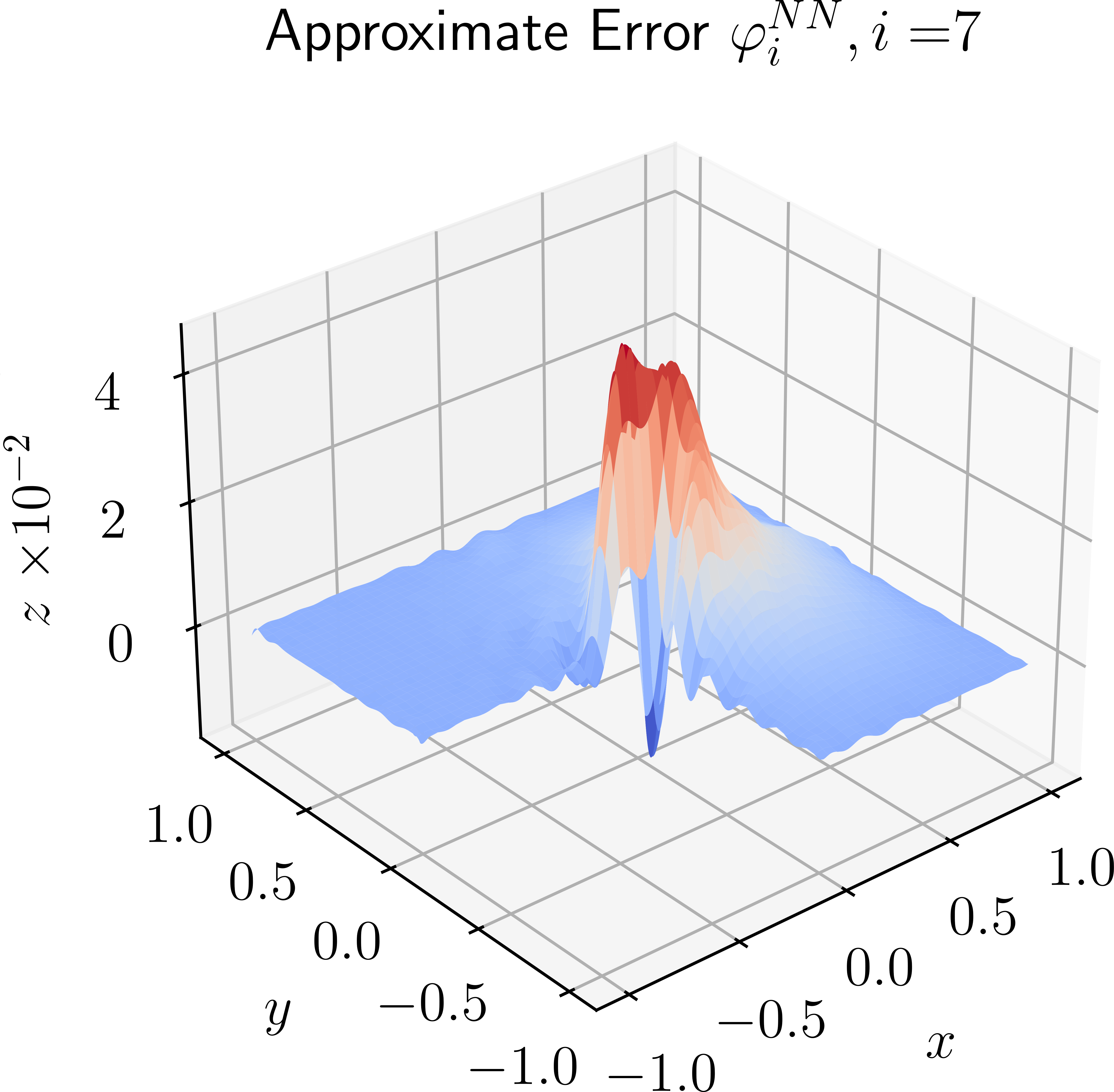}
            \quad
            \includegraphics[width=1.8in]{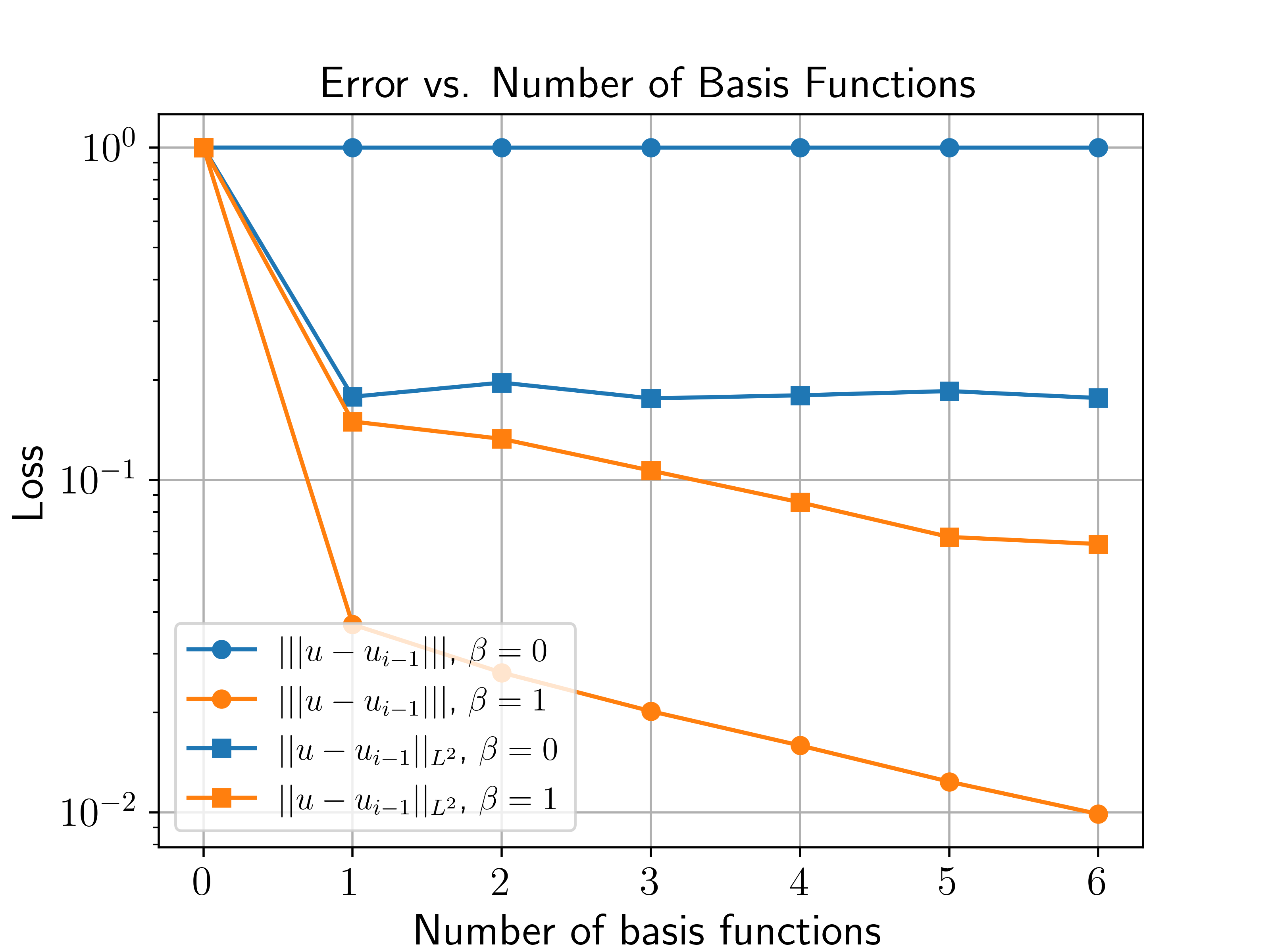}
		\caption{Example \ref{ex:poisson rlambda}: for $\beta=1$, the exact errors $u-u_{i-1}$ (top row) and approximate errors which correspond to learned basis functions $\varphi_{i}^{NN} \approx u-u_{i-1}$ (middle row). The energy and and $L^{2}$ errors are provided for both $\beta = 0$ and $\beta = 1$.}
		\label{fig:poisson rlambda3}
	\end{figure}
 
	\begin{example} \label{ex:poisson rlambda}
		We consider the case of the true solution $u(r,\theta) = r^{\lambda}\sin{\theta}$ for $\lambda = 1/4$ in the domain $\Omega = (-1,1)^{2} \backslash (-1,0)^{2}$. In this case, $f \notin L^{2}(\Omega)$ but $f \in L^{2}_{\beta}(\Omega)$ for $\beta > 1-\lambda$. We demonstrate that if $\beta = 0$ in \eqref{eq:poisson lsq}, then the approximation obtained is incorrect and does not converge to the true solution. Figure \ref{fig:poisson rlambda1} shows the approximation using 7 basis functions. All hyperparameters for this example are provided in Appendix \ref{app:examples}.
	
		We observe spurious oscillations of the approximation near $x=0$ when $\beta=0$ which are largely removed when $\beta=1$. Figures \ref{fig:poisson rlambda2}-\ref{fig:poisson rlambda3} show the true errors $u-u_{i-1}$ and basis functions $\varphi_{i}^{NN}$ for $\beta=0$ and $\beta=1$, respectively. We observe that in the $\beta=0$ case, the learned basis functions contain many high frequency modes not present in the corresponding error functions. In the $\beta=1$ case, the basis functions take on large values near $x=0$ and correctly identify the presence of large errors in the vicinity of the singularity, but convergence is still slow.
	\end{example}

	\begin{figure}[t!]
		\centering
		\includegraphics[width=1.3in]{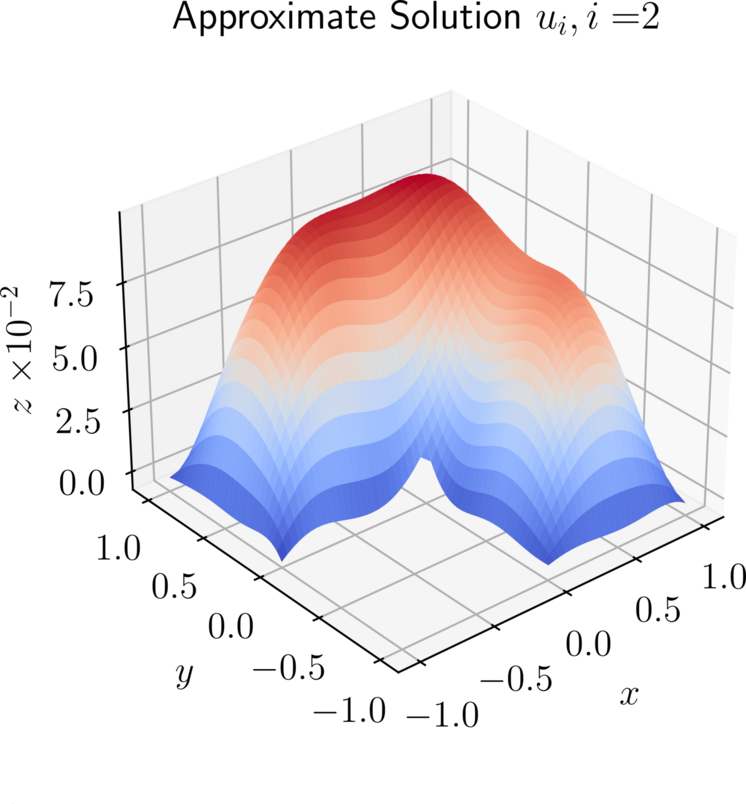}
		\quad
		\includegraphics[width=1.3in]{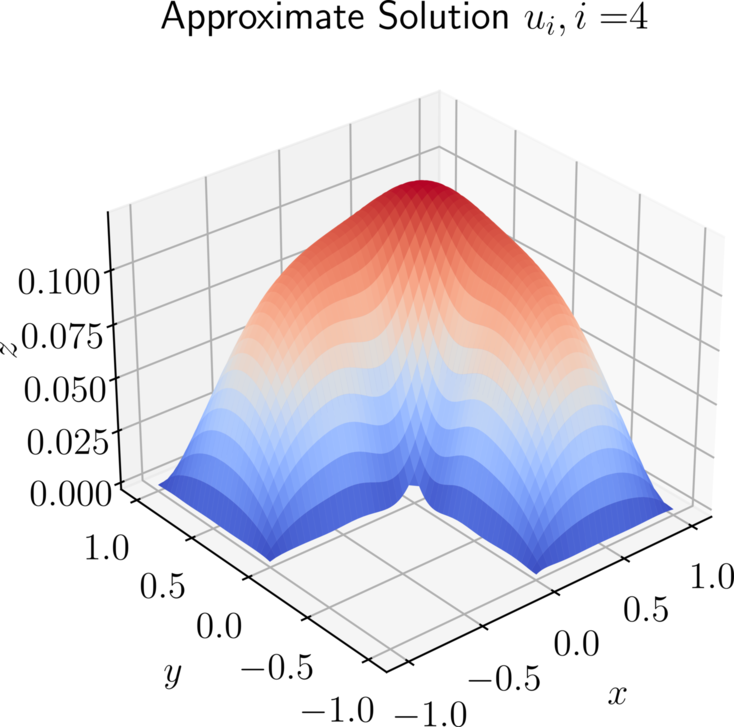}
		\quad
		\includegraphics[width=1.3in]{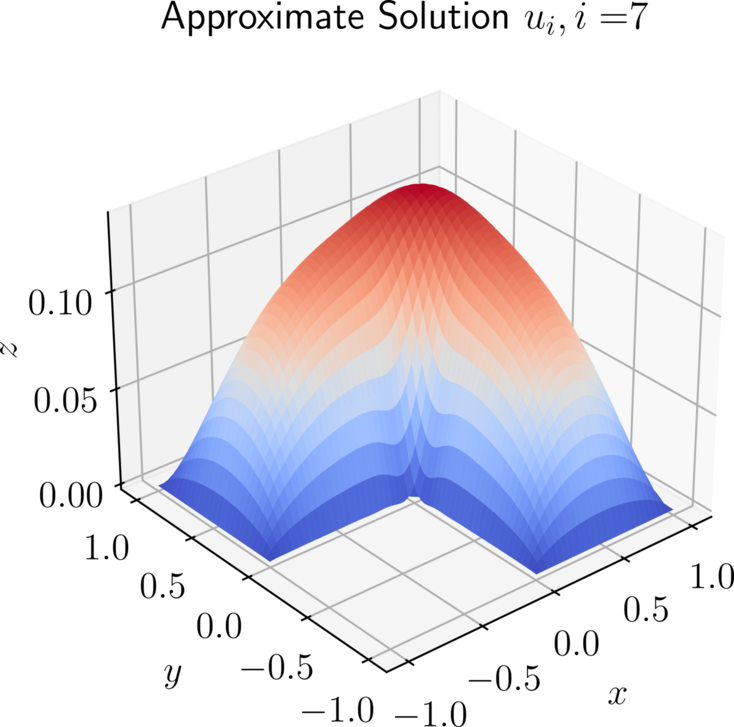}
		
		\includegraphics[width=1.3in]{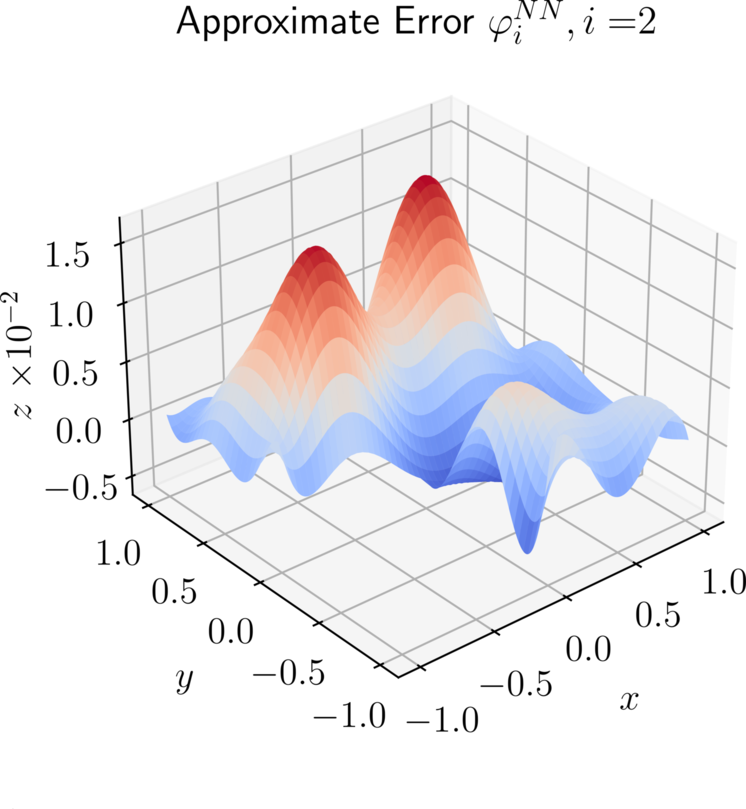}
		\quad
		\includegraphics[width=1.3in]{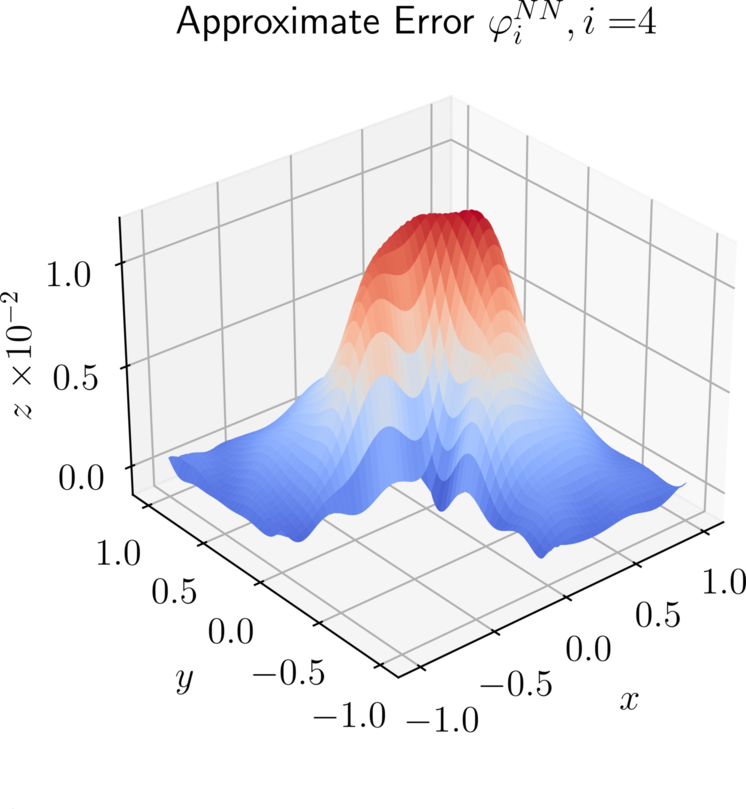}
		\quad
		\includegraphics[width=1.3in]{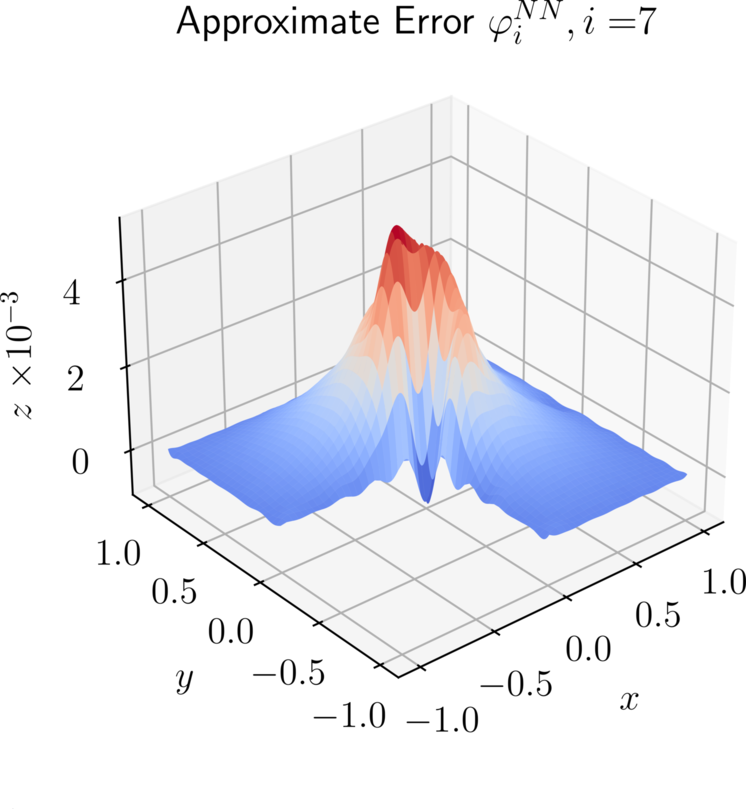}

		\caption{Example \ref{ex:poisson Lshaped}: for $\beta=4/3$, the approximate solutions $u_{i}$ and basis functions $\varphi_{i}^{NN}$ learned using the standard Galerkin neural network approach in Section \ref{sec:gnn}.}
		\label{fig:poisson Lshaped1}
	\end{figure}

	\subsubsection{Knowledge-based functions}
	
	Suppose $\Omega$ contains a single corner of interest $x^{(0)}$ with interior angle $\alpha^{(0)}$. For instance, the $L$-shaped domain $(-1,1)^{2}\backslash (-1,0)^{2}$ contains a non-convex corner with angle $\alpha^{(0)} = 3\pi/2$ at the origin. Solutions of \eqref{eq:poisson2} in such a domain have the well-known structure
	\begin{align}
		u(x,y) = \sum_{\lambda_{j}} c_{\lambda_{j}}r^{\lambda_{j}} \sin(\lambda_{j}\theta) + u^{*}(x,y),
	\end{align}
	
	\noindent where $u^{*} \in H^{2}(\Omega)$ and $\lambda_{n}$ are the eigenvalues of the operator pencil corresponding to the Laplacian operator \cite{mazya}. The eigenvalues are given by
	\begin{align} \label{eq:poisson eigenvalues}
		\lambda_{j} = \frac{j\pi}{\alpha^{(0)}}, \;\;\;j=1,2,3,\dots.
	\end{align}

	\begin{figure}[t!]
            \captionsetup{font={color=red}}
		\centering
		\includegraphics[width=1.3in]{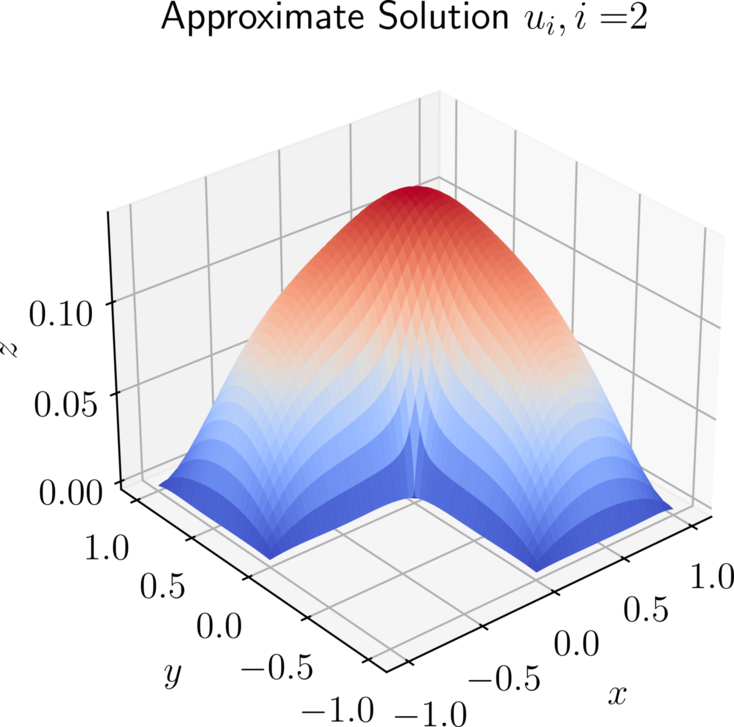}
		\quad
		\includegraphics[width=1.3in]{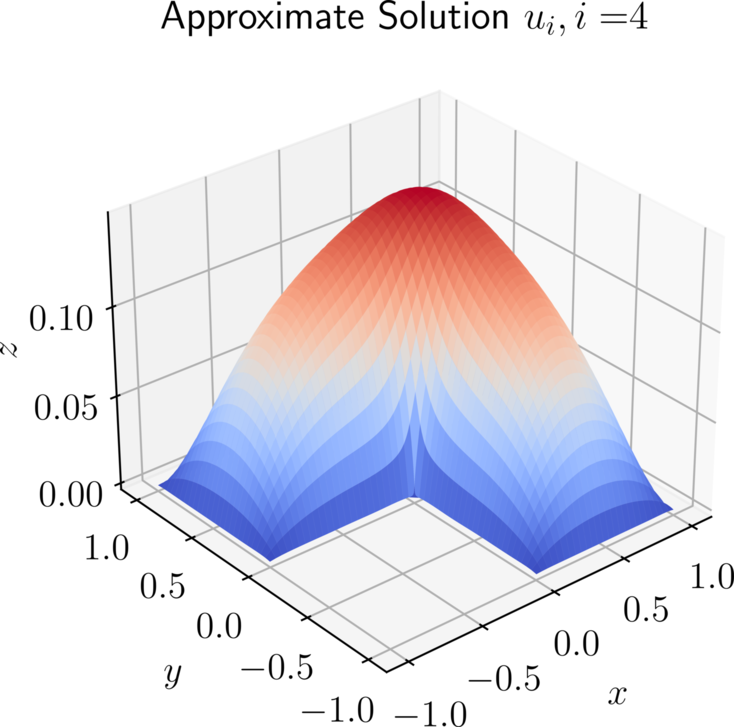}
		\quad
		\includegraphics[width=1.3in]{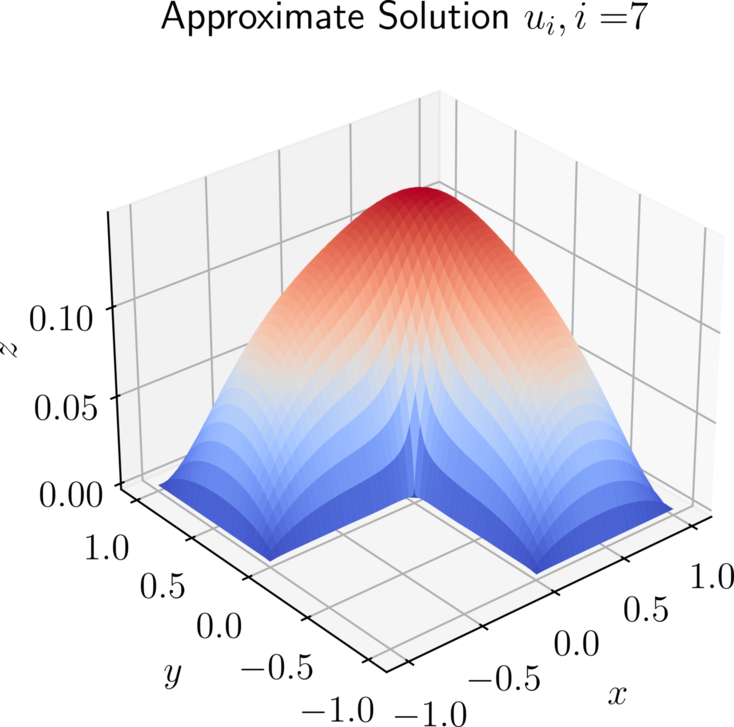}
		
		\includegraphics[width=1.3in]{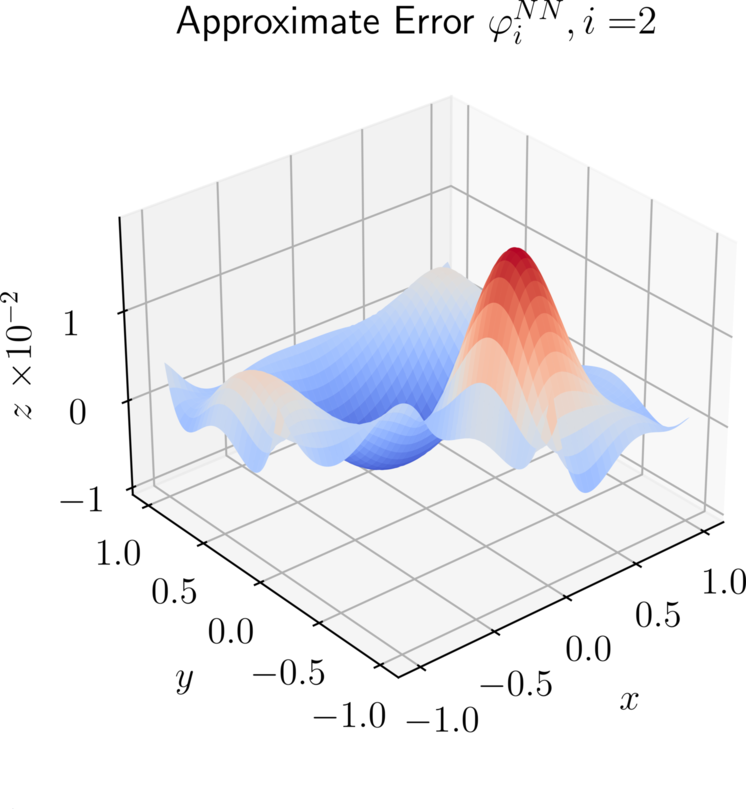}
		\quad
		\includegraphics[width=1.3in]{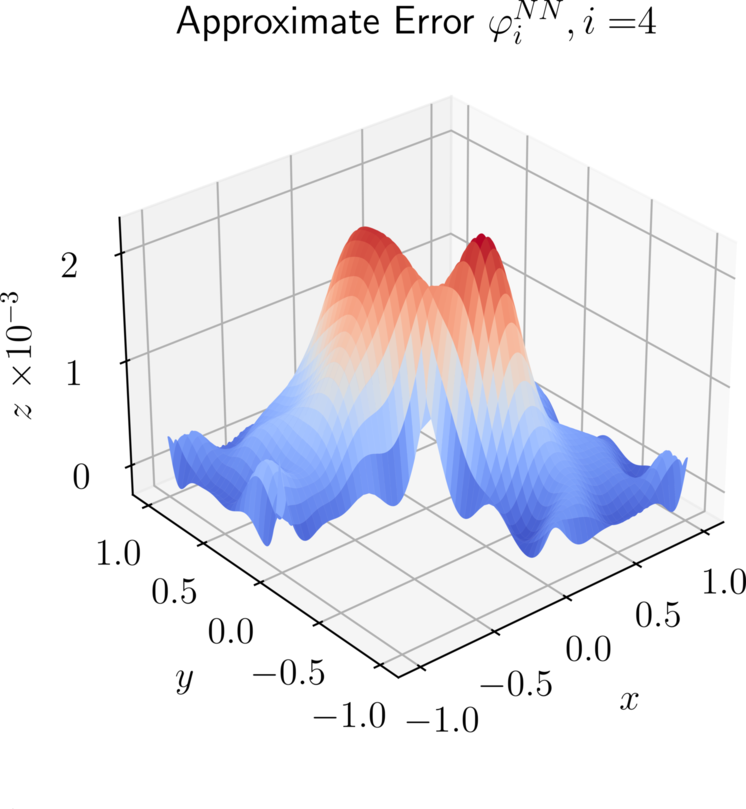}
		\quad
		\includegraphics[width=1.3in]{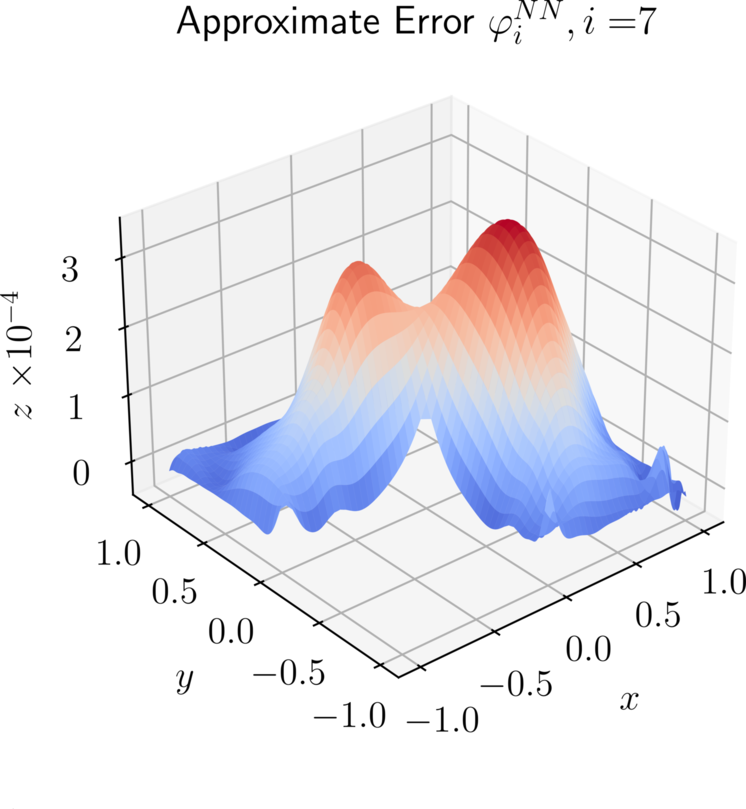}

          \includegraphics[width=1.9in]{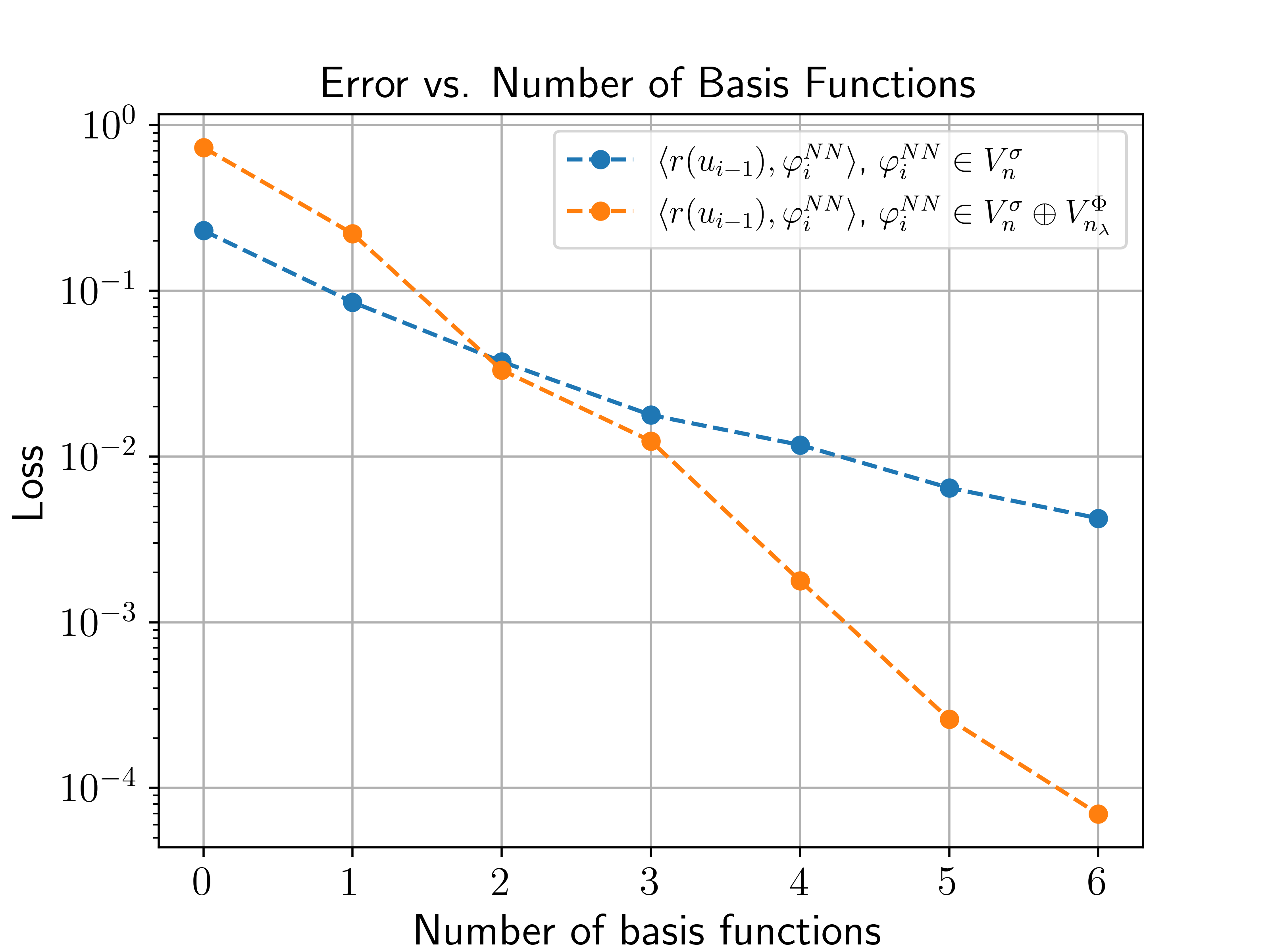}

		\caption{Example \ref{ex:poisson Lshaped}: for $\beta=4/3$, the approximate solutions $u_{i}$ (top row) and basis functions $\varphi_{i}^{NN}$ learned using the extended Galerkin neural network approach with $\Psi(x;\mu) = r(x)^{\mu}\sin(\mu\theta(x))$ (bottom row). A comparison of the loss function vs. number of learned basis functions is also provided for the basic feedforward network architecture vs. the extended architecture described in Figure \ref{fig:enriched GNN}.}
		\label{fig:poisson Lshaped2}
	\end{figure}
	
	\begin{example} \label{ex:poisson Lshaped}
		We consider the problem with homogeneous Dirichlet boundary condition in the domain $\Omega = (-1,1)^{2}\backslash (-1,0)^{2}$ and data $f = 1$. This problem was already considered in \cite{gnn1} for which the usual $H^{1}$ variational formulation of the Poisson equation was utilized. Nevertheless, this problem demonstrates many of the explemary singular features we wish to approximate and serves as an excellent example for how to utilize knowledge-based functions in the Galerkin-neural network framework.
        \end{example}

		 While there is no closed form solution to this problem, we provide the value of the \emph{a posteriori} error indicator $\langle r(u_{i-1}),\varphi_{i}^{NN} \rangle$ for each basis function learned as well as the basis functions $\varphi_{i}^{NN}$ themselves, which provide a function representation of the error $u-u_{i-1}$.  Figure \ref{fig:poisson Lshaped1} shows the first several basis functions using a standard feedforward neural network for learning each basis function. All hyperparameters for this example are provided in Appendix \ref{app:examples}. We observe that the approximation converges quite slowly near $x=0$ and the basis functions $\varphi_{i}^{NN}$ show that the magnitude of the pointwise error is large in the vicinity of the singularity. 
		
		Next, we repeat the simulation, this time with an augmented neural network architecture. Namely, we compute $\varphi_{i}^{NN}$ according to \eqref{eq:singular basis fn} with $\Psi(r,\theta;\mu) = r^{\mu}\sin(\mu\theta)$.  Since the eigenvalues are known, we simply take $\mu_{j} = 2j/3$ and $\mathcal{M} = 20$. However, we shall demonstrate in Section \ref{sec:eigenvalue training} how Galerkin neural networks may be used to approximate the eigenvalues $\lambda_{j}$ when their values are unknown. Figure \ref{fig:poisson Lshaped2} shows the analogous results to Figure \ref{fig:poisson Lshaped1}. We observe significantly faster convergence and the approximate pointwise errors decrease steadily from $\mathcal{O}(10^{-2})$ to $\mathcal{O}(10^{-4})$. Additionally, we note that the singularity is clearly captured correctly by the knowledge-based functions.

	\subsection{Stokes flow}
	
	The next application we consider is the incompressible steady Stokes flow \eqref{eq:stokes strong}. 
	
	\subsubsection{Regularity and variational formulation}
	Schwab and Guo proved the following result in \cite{schwabguo} for a polygonal domain $\Omega$.
	
	\begin{theorem} \label{thm:stokes a priori}
		\cite{schwabguo} Suppose $\mathbf{f} \in \mathbf{L^{2}_{\beta}}(\Omega)$, $g \in H^{1}_{\beta}(\Omega)$, $\mathbf{u}_{D} \in \mathbf{H^{3/2}_{\beta}}(\partial\Omega)$, and $\int_{\Omega} g\;d\mathbf{x} = -\int_{\partial\Omega} \mathbf{u}_{D}\cdot\mathbf{n}\;ds$, where $\beta = (\beta_{1},\dots,\beta_{M})$, $\beta_{i} > 1 - \kappa_{1}^{i}$ and $\kappa_{1}^{i}$ is the smallest real part of the eigenvalues of the operator pencil corresponding to \eqref{eq:stokes strong} with positive real part (see, e.g. \eqref{eq:eigenvalues}). Then \eqref{eq:stokes strong} admits a unique solution $(\mathbf{u},p) \in \mathbf{H^{2}_{\beta}}(\Omega) \times H^{1}_{\beta}(\Omega)/\mathbb{R}$ and 
		\begin{align} \label{eq:stokes estimate}
			||\mathbf{u}||_{\mathbf{H^{2}_{\beta}}(\Omega)} + ||p||_{H^{1}_{\beta}(\Omega)} \leqslant C(||\mathbf{f}||_{\mathbf{L^{2}_{\beta}}(\Omega)} + ||g||_{H^{1}_{\beta}(\Omega)} + ||\mathbf{u}_{D}||_{\mathbf{H^{3/2}_{\beta}}(\Omega)})
		\end{align}
		
		\noindent for some $C>0$ independent of $\mathbf{u}$, $p$, $\mathbf{f}$, $g$, and $\mathbf{u}_{D}$. 
	\end{theorem}
        
	Most notable is that \eqref{eq:stokes estimate} suggests that the residual of the incompressibility equation in \eqref{eq:stokes strong} should actually be taken in the $H^{1}_{\beta}$ norm, which leads to the continuous and coercive formulation: seek $(\mathbf{u},p) \in \mathbf{H^{2}_{\beta}}(\Omega) \times H^{1}_{\beta}(\Omega)/\mathbb{R}$ such that
	\begin{align}
            \begin{aligned}
		(-\Delta\mathbf{u} + \nabla p, -&\Delta\mathbf{v} + \nabla q)_{\beta,\Omega} + (\text{div}\;\mathbf{u}, \text{div}\;\mathbf{v})_{H^{1}_{\beta}(\Omega)} + \delta(\mathbf{u},\mathbf{v})_{H^{3/2}_{\beta}(\partial\Omega)}\\
		&= (\mathbf{f}, -\Delta\mathbf{v}+\nabla q)_{\beta,\Omega} + (g, \text{div}\;\mathbf{v})_{H^{1}_{\beta}(\Omega)} + \delta(\mathbf{u}_{D},\mathbf{v})_{\mathbf{H^{3/2}_{\beta}}(\Omega)}
            \end{aligned}
	\end{align}
	
	\noindent for all $(\mathbf{v},q) \in \mathbf{H^{2}_{\beta}}(\Omega) \times H^{1}_{\beta}(\Omega)/\mathbb{R}$. This problem takes the form \eqref{eq:general lsq} with $\mathcal{X} := \mathbf{H^{2}_{\beta}}(\Omega) \times H^{1}_{\beta}(\Omega)/\mathbb{R}$,
	\begin{align} \label{eq:stokes bilinear}
		a_{LS}((\mathbf{u},p), (\mathbf{v},q)) := (-\Delta\mathbf{u} &+ \nabla p, -\Delta\mathbf{v} + \nabla q)_{\beta,\Omega}\\
            &+ (\text{div}\;\mathbf{u}, \text{div}\;\mathbf{v})_{H^{1}_{\beta}(\Omega)} + \delta(\mathbf{u},\mathbf{v})_{H^{3/2}_{\beta}(\partial\Omega)},\notag
	\end{align}
        \noindent and
        \begin{align}
            F_{LS}((\mathbf{v},q)) := (\mathbf{f}, -\Delta\mathbf{v}+\nabla q)_{\beta,\Omega} + (g, \text{div}\;\mathbf{v})_{H^{1}_{\beta}(\Omega)} + \delta(\mathbf{u}_{D},\mathbf{v})_{\mathbf{H^{3/2}_{\beta}}(\Omega)}.
        \end{align}

	\subsubsection{Knowledge-based functions} \label{sec:stokes knowledge}
	
	Many features of interest in the Stokes flow develop in the corners of $\Omega$, including singularities at reentrant corners and eddies in convex corners that are driven by the far-field flow. In these cases, the solution structure is well-known. Noting that the velocity is the curl of the streamfunction $\psi$, $\mathbf{u} = \nabla^{\perp}\psi$, and that the streamfunction satisfies a biharmonic equation, the streamfunction takes the form
	\begin{align} \label{eq:singular structure}
		\psi(x,y) = \sum_{j=1}^{M} \chi(r^{(j)})\sum_{\lambda_{n}^{(j)}} c_{\lambda_{n}^{(j)}}r^{\lambda_{n}^{(j)}} \Psi_{\lambda_{n}^{(j)}}(\theta^{(j)}) + \psi^{*}(x,y),
	\end{align}
	
	\noindent where $(r^{(j)}, \theta^{(j)})$ are the local polar coordinates centered at $x^{(j)}$, $\{\lambda_{n}^{(j)}\}$ are the eigenvalues associated with the operator pencil of the biharmonic operator at vertex $x^{(j)}$ and $\Psi_{\lambda_{n}^{(j)}}$ are the corresponding generalized eigenfunctions, $\chi_{j} \in C^{s}(\bar{\Omega})$ is a cutoff function centered at $x^{(j)}$, and $\psi^{*} \in H^{4}(\Omega)$. The corresponding velocity and pressure solutions are derived from \eqref{eq:singular structure}. See Appendix \ref{app:stokes} for full details.
	
	For $\lambda \neq 0,1,2$, the eigenvalues of the biharmonic operator pencil at $x^{(j)}$ are given implicitly by \cite{blumrannacher}
	\begin{align} \label{eq:eigenvalues}
		\sin^{2}((\lambda-1)\alpha_{j}) - (\lambda-1)^{2}\sin^{2}(\alpha_{j}) = 0, 
	\end{align}
	
	\noindent where $\alpha_{j}$ is the interior angle of $x^{(j)}$. The eigenfunctions \cite{blumrannacher} are given by
	\begin{align} \label{eq:eigenfunctions}
		\Psi_{\lambda}(\theta) = A\sin(\lambda\theta) + B\cos(\lambda\theta) + C\sin((\lambda-2)\theta) + D\cos((\lambda-2)\theta),
	\end{align}

	\noindent with the constants $A,B,C,D$ chosen so that $\Psi$ satisfies suitable boundary conditions.
	\begin{figure}[t!]
		\centering
		\includegraphics[width=2.1in]{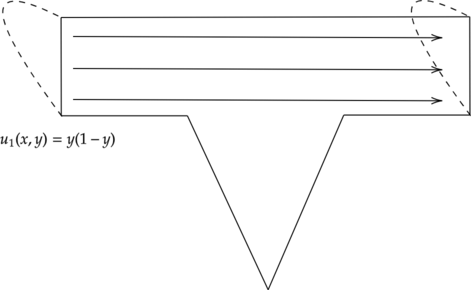}
		\caption{Channel with triangular cavity for Example \ref{ex:triangular}. The rectangle channel is described by $(-2,2) \times (0,2)$ while the triangular cavity has vertices at $(-1,0)$, $(1,0)$, and $(0,-3)$.}
            \label{fig:domain}
	\end{figure}
	 In the case when $\alpha_{j} \lessapprox 146^{\circ}$, solutions of \eqref{eq:eigenvalues} are complex and we replace $r^{\lambda_{n}}\Psi_{\lambda_{n}}(\theta)$ by $\mathfrak{Re}[r^{\lambda_{n}}\Psi_{\lambda_{n}}(\theta)]$ in \eqref{eq:singular structure}. Interestingly, complex eigenvalues are associated with infinite sequences of eddies in the flow due to disturbances in the far field. This behavior is described mathematically by, for instance
	\begin{align*}
		\mathfrak{Re}(r^{a+ib}\cos((a+ib)\theta)) = r^{a}\cos(b\log{r})\cos(a\theta)\cosh(b\theta).
	\end{align*}
	
	\noindent One can show, as in \cite{moffatt}, that each successive eddy decays exponentially in magnitude. Thus, Stokes flows in polygonal domains must be resolved to a very high degree of accuracy in order to adequately capture these eddies.

	\begin{example} \label{ex:triangular}
		Let $\Omega$ be the channel with triangular cavity as shown in Figure \ref{fig:domain}. The flow takes a parabolic profile prescribed at the inflow and outflow. Namely, the boundary conditions are 
        \begin{align*}
            u_{1}(x,y) = \begin{dcases}
                y(2-y), &x=-2 \;\text{or}\;2\\
                0, &\text{else},
            \end{dcases} \;\;\;\;\;\;u_{2}(x,y) = 0.
        \end{align*}
        We specify $\mathbf{f} = \mathbf{0}$ and $g = 0$. At the reentrant corners, we expect a singularity to manifest in both the velocity and pressure. At the bottom corner, we expect infinitely cascading eddies driven by the channel flow. As described in Section \ref{sec:stokes knowledge}, we shall supplement the standard feedforward network with $\Psi(r,\theta;\mu) = \chi(r)\mathfrak{Re}[r^{\mu}\Psi_{\mu}(\theta)]$ in the bottom corner and $\Psi(r,\theta;\mu) = \chi(r) r^{\mu}\Psi_{\mu}(\theta)$ at the reentrant corners. The exact values of the eigenvalue $\lambda$ are $\lambda=1.58223$ in the reentrant corners and $\lambda = 7.56813 + 3.37941i$ in the bottom corner with $\mathcal{M} = 1$. Precise details may be found in Appendix \ref{app:stokes}. As with Example \ref{ex:poisson Lshaped}, we shall demonstrate in Section \ref{sec:eigenvalue training} how these values of $\lambda$ may be approximated using Galerkin neural networks in the event that they are unknown.
		\begin{figure}[t!]
			\centering
			\includegraphics[width=1.1in]{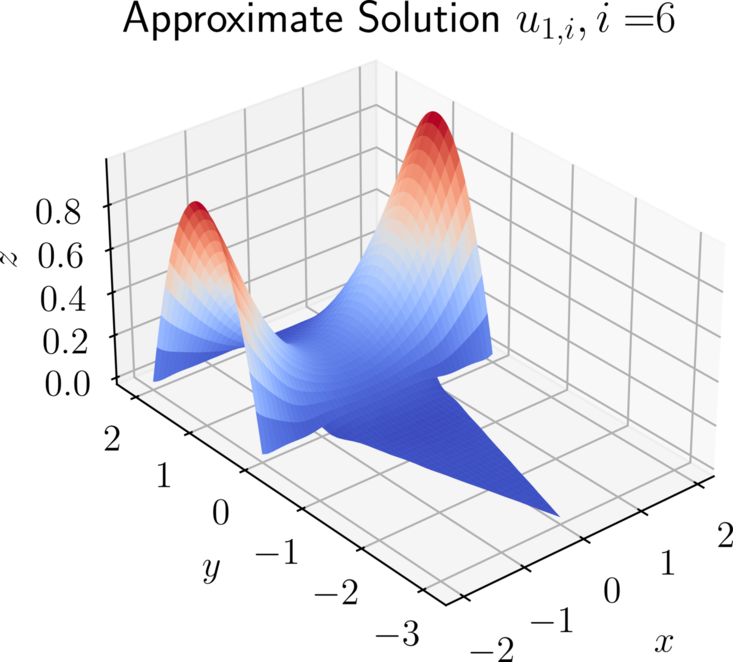}
			\quad
			\includegraphics[width=1.1in]{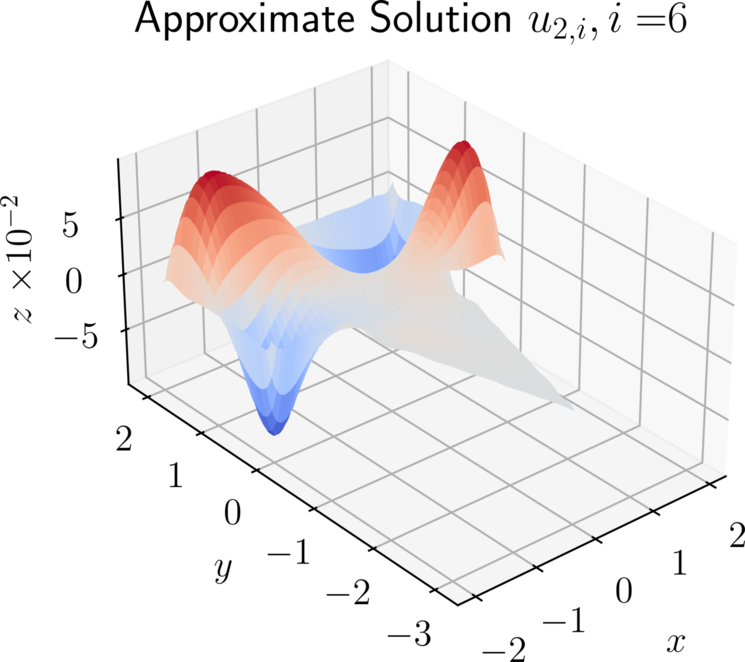}
			\quad
			\includegraphics[width=1.1in]{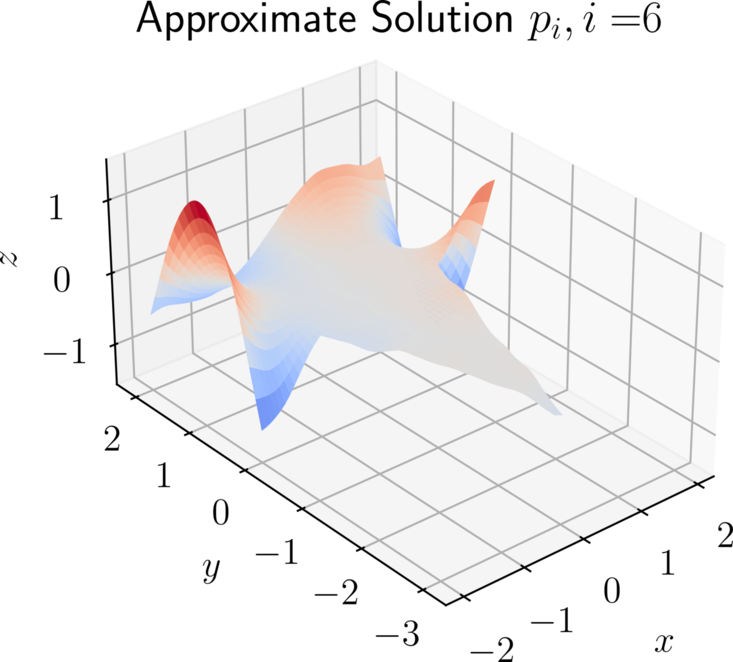}
                \quad
			\includegraphics[width=0.95in]{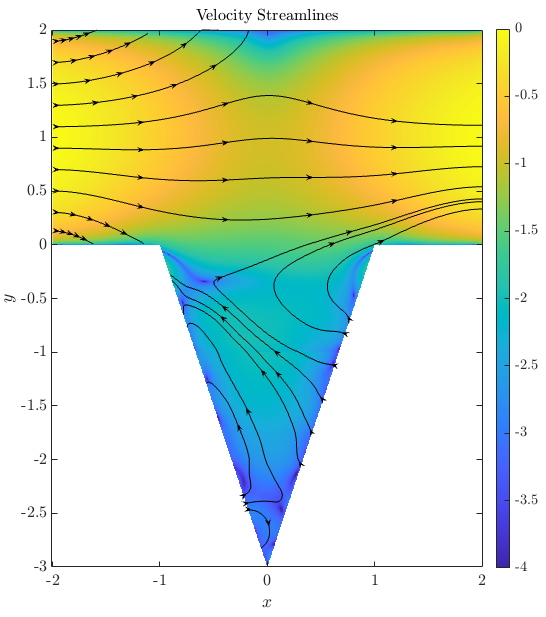}
			
			\includegraphics[width=1.1in]{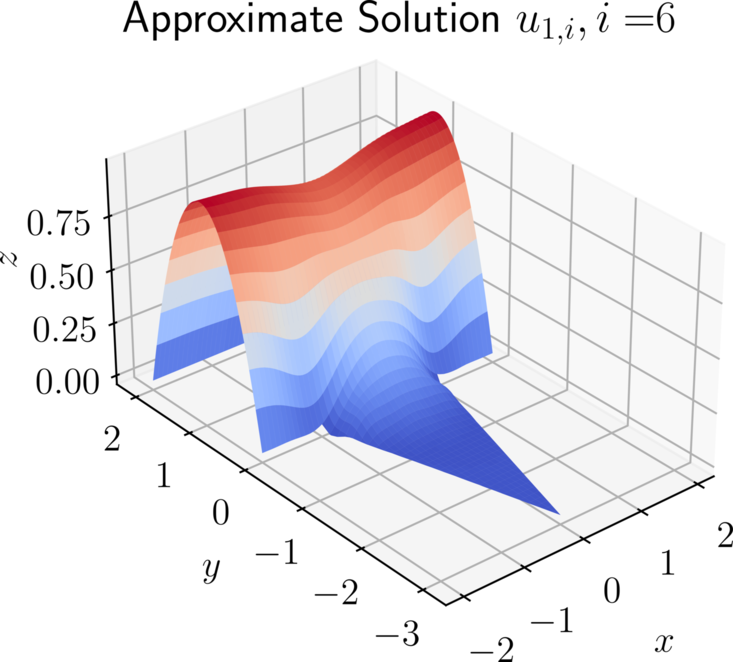}
			\quad
			\includegraphics[width=1.1in]{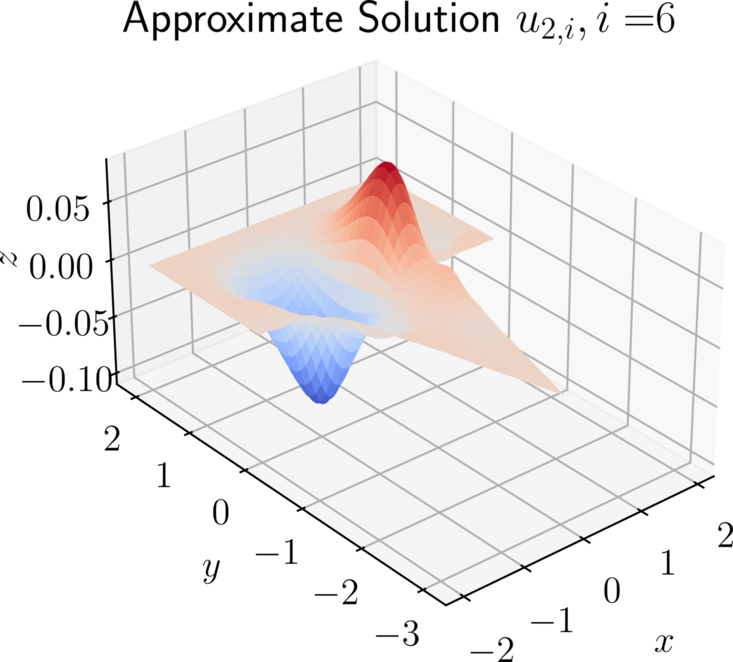}
			\quad
			\includegraphics[width=1.1in]{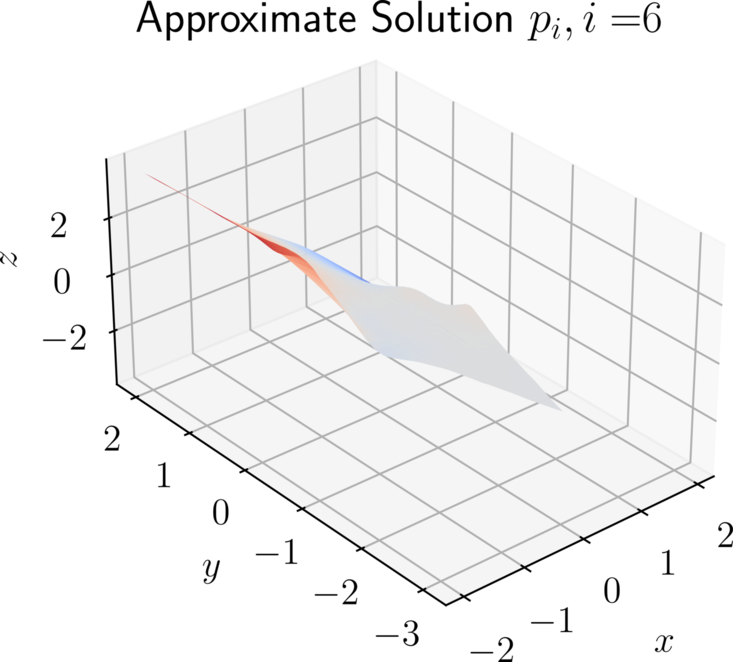}
                \quad
			\includegraphics[width=0.95in]{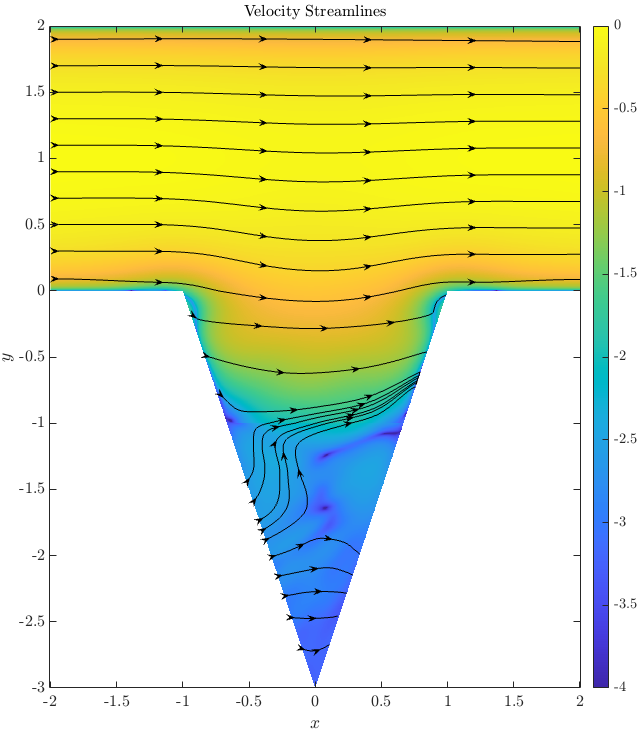}
			\caption{Example \ref{ex:triangular}: velocity and pressure with $\beta=0$ (top) and $\beta=5/3$ (bottom) in \eqref{eq:stokes bilinear} with basis functions learned using the standard Galerkin neural network approach described in Section \ref{sec:gnn}. The approximation is non-physical in the channel when the least squares variational problem is posed on unweighted Sobolev spaces. Without enriching the network with knowledge-based functions, the flow behavior at the non-convex corners and in the cavity remains unresolved.}
			\label{fig:stokes beta0+53}
		\end{figure}
		
		We first demonstrate the importance of choosing the Sobolev weight parameter $\beta$ appropriately. In contrast with Example \ref{ex:poisson rlambda}, this example contains smooth data $f$ and $g$ and continuous data $\mathbf{u}_{D}$. All hyperparameters for this example are provided in Appendix \ref{app:examples}. Figure \ref{fig:stokes beta0+53} shows the approximated velocities and pressures for $\beta = 0$. We shall use $\varphi_{1,i}^{NN}$, $\varphi_{2,i}^{NN}$, and $q_{i}^{NN}$ to denote the $i$th basis functions for the velocities and pressure, respectively. The resulting velocity and pressure fields in Figure \ref{fig:stokes beta0+53} show that even after learning six basis functions, the flow is nonphysical; we expect a near constant velocity $u_{1}$ along the $x$-direction in the channel. Figure \ref{fig:stokes beta0+53} also shows the analogous results for $\beta=5/3$. Despite not adequately capturing the singularities and the expected sequence of eddies, the approximate solutions do not exhibit spurious oscillations and the behavior in the channel is as expected.
		\begin{figure}[t!]
			\centering
			\includegraphics[width=1.3in]{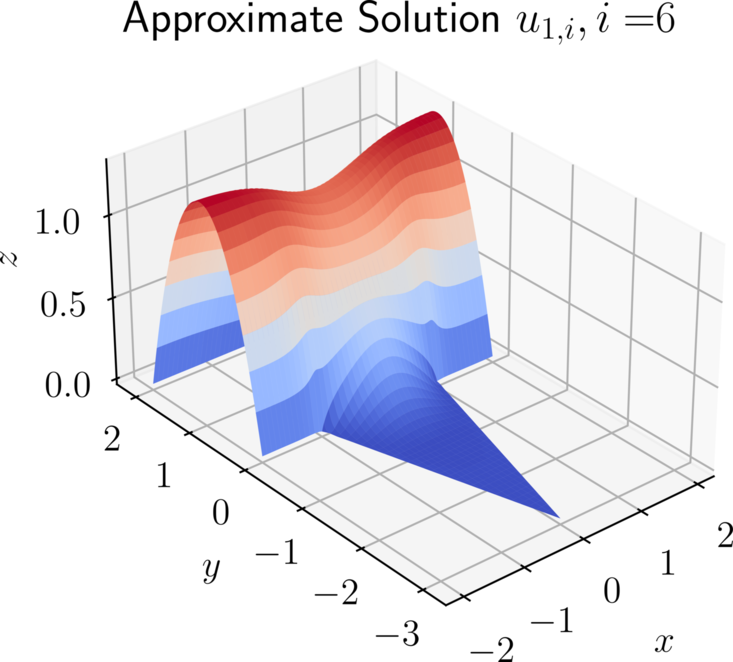}
			\quad
			\includegraphics[width=1.3in]{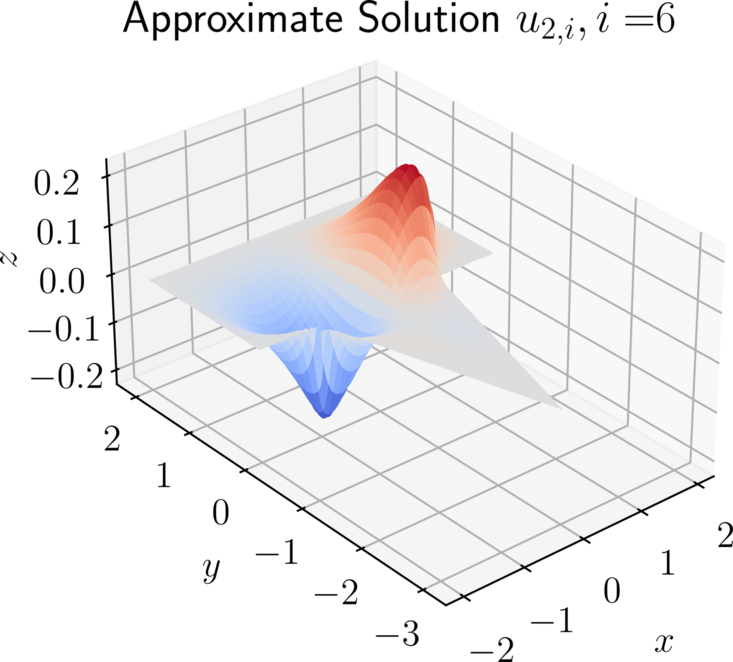}
			\quad
			\includegraphics[width=1.3in]{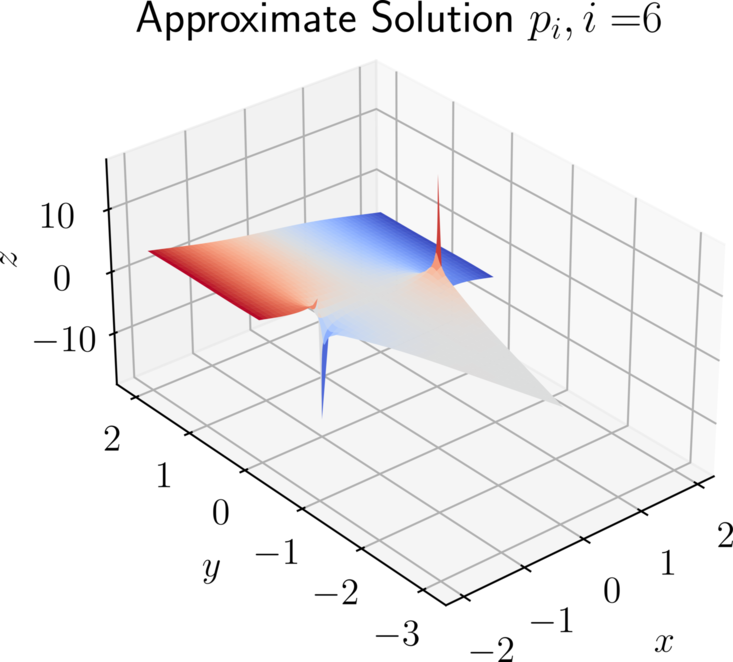}
			
			\includegraphics[width=1.3in]{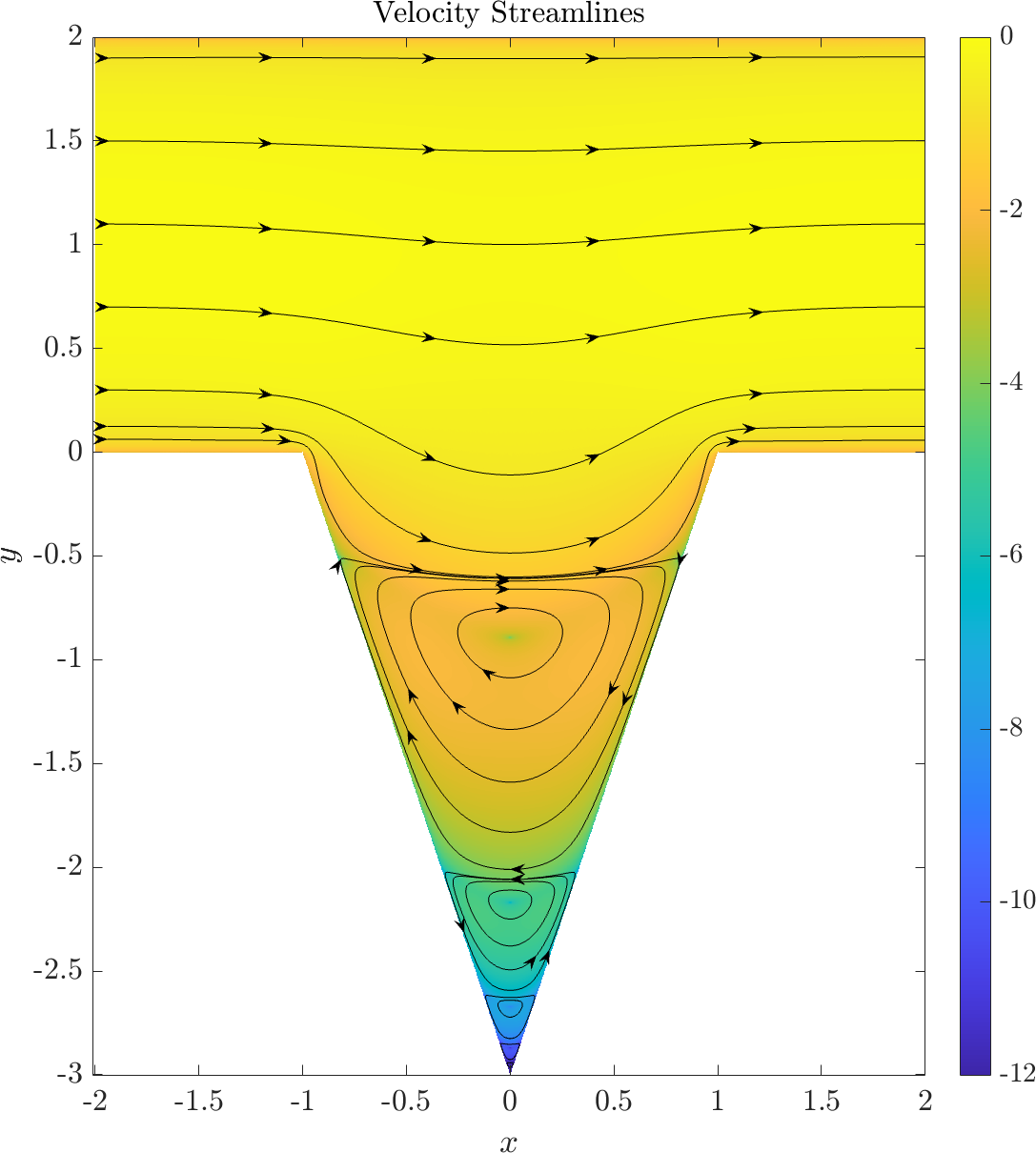}
                \quad
                \includegraphics[width=1.25in]{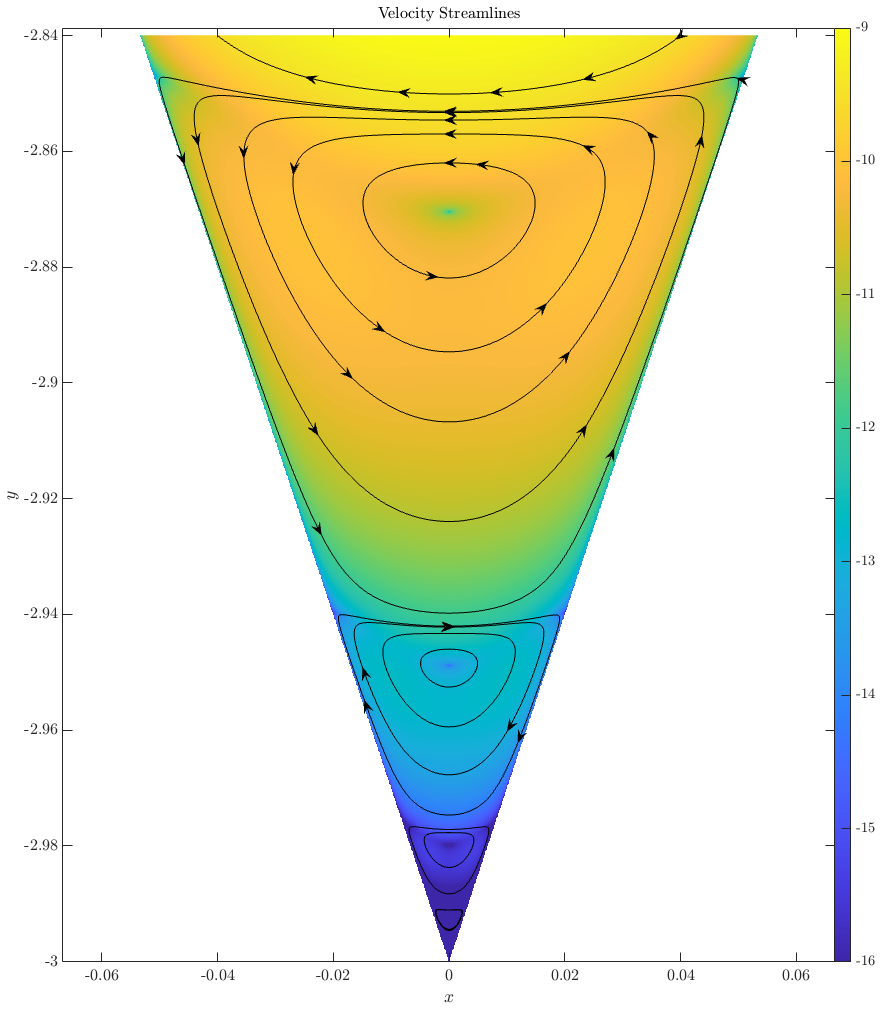}
                \quad
                \includegraphics[width=1.25in]{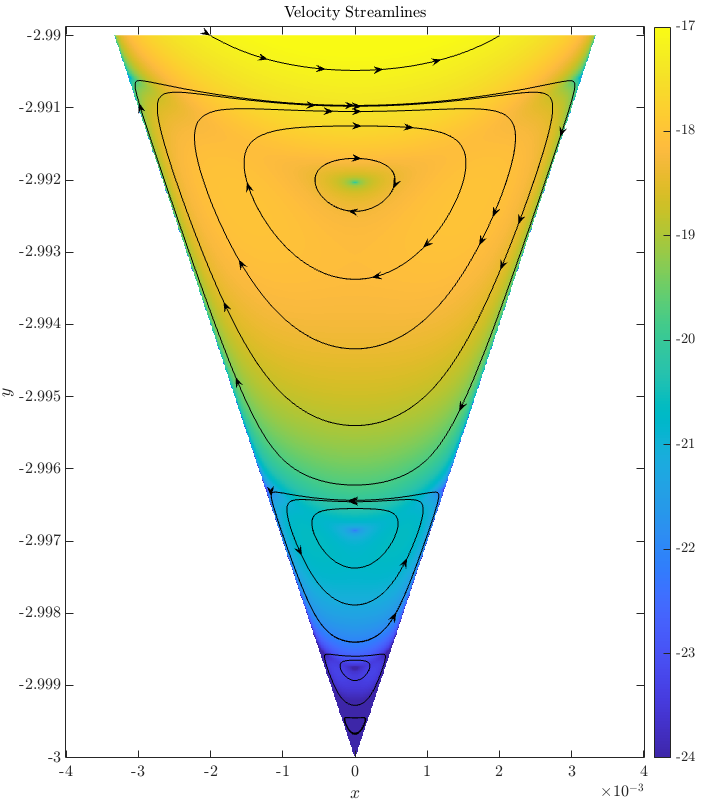}
			\caption{Example \ref{ex:triangular}: velocity and pressure with $\beta=5/3$ in \eqref{eq:stokes bilinear} and basis functions learned using the extended Galerkin neural network approach with extended neural network architecture (top row) and velocity streamlines with two progressive zooms showing resolution of the Moffatt eddies (bottom row).}
			\label{fig:stokes beta53}
		\end{figure}
		
		We fix $\beta=5/3$ next and approximate the solution using six basis functions learned by extended Galerkin neural networks. We remarkably observe that with the supplemental activation functions, i.e. $\varphi_{1,i}^{NN}, \varphi_{2,i}^{NN}, q_{i}^{NN} \in V_{\mathbf{n},L}^{\sigma} \oplus V_{\mathcal{M}}^{\Psi}$, both the singular features at the reentrant corners and the eddies in the cavity are captured accurately. As Figure \ref{fig:stokes beta53} shows, we are able to capture at least ten eddies. The basis functions used to approximate the velocities and pressures are shown in Figure \ref{fig:stokes beta53 basis}. 
		
		Finally, when $\beta$ is too large, the solution space contains functions which may have worse singular behavior than the true solution due to the strong damping effect of $r^{\beta}$ on the singularities. In particular, the bilinear operator $a_{LS}$ may admit as arguments functions with far stronger singularities than those of the true solution. We show the results for $\beta=4$ in Figure \ref{fig:stokes beta4}.
		\begin{figure}[t!]
			\centering
			\includegraphics[width=1.3in]{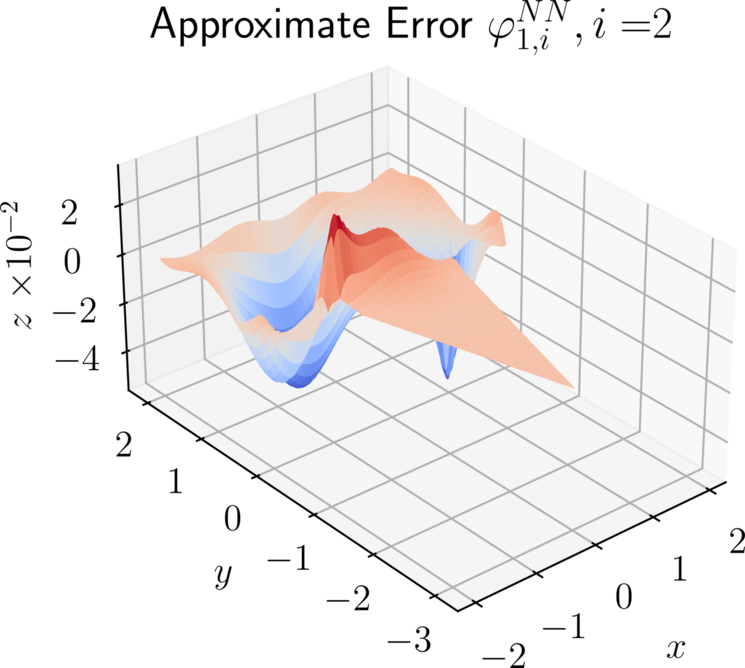}
			\quad
			\includegraphics[width=1.3in]{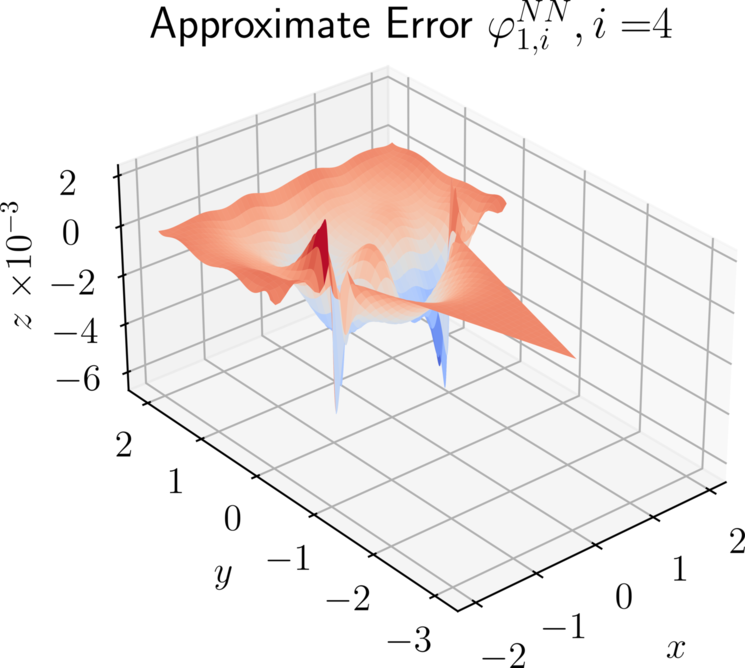}
			\quad
			\includegraphics[width=1.3in]{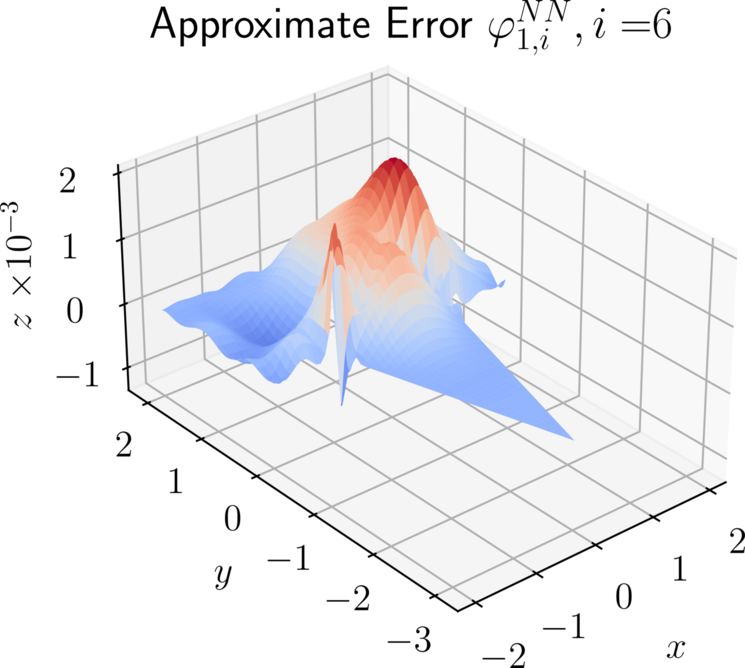}

			\includegraphics[width=1.3in]{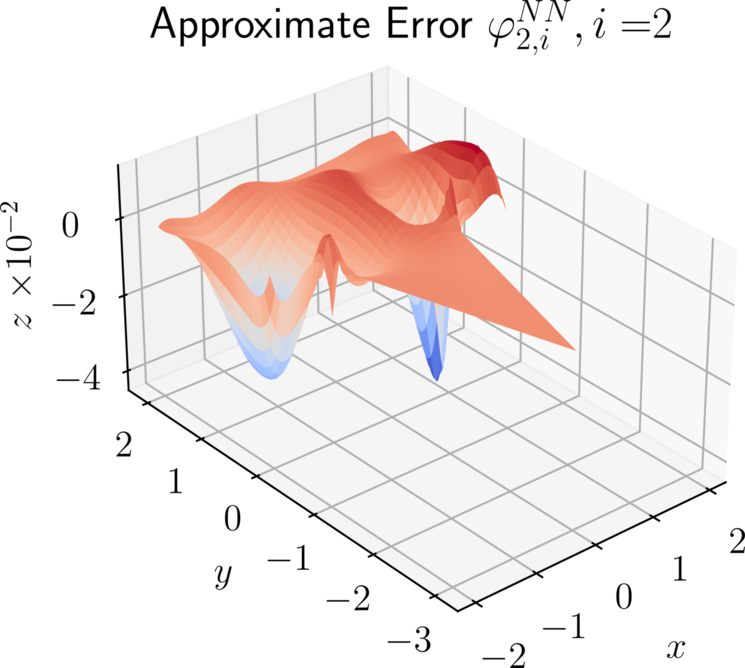}
			\quad
			\includegraphics[width=1.3in]{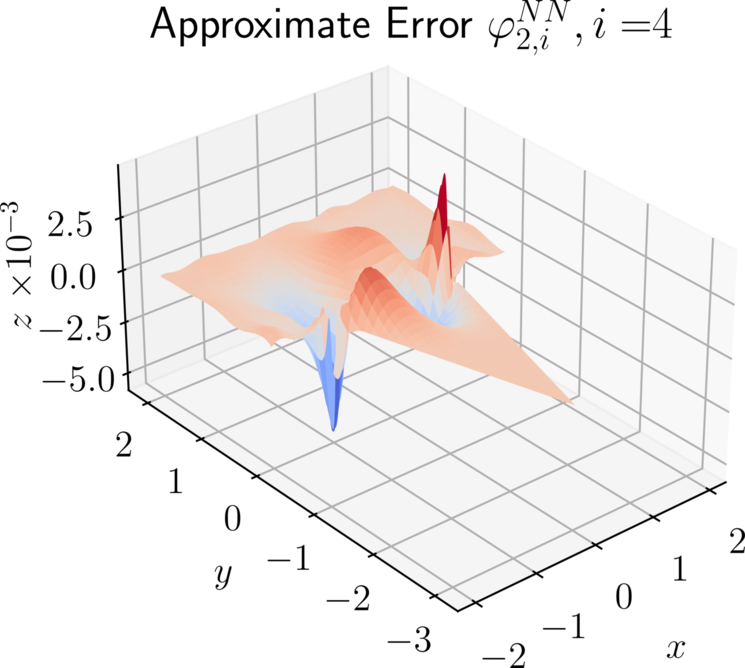}
			\quad
			\includegraphics[width=1.3in]{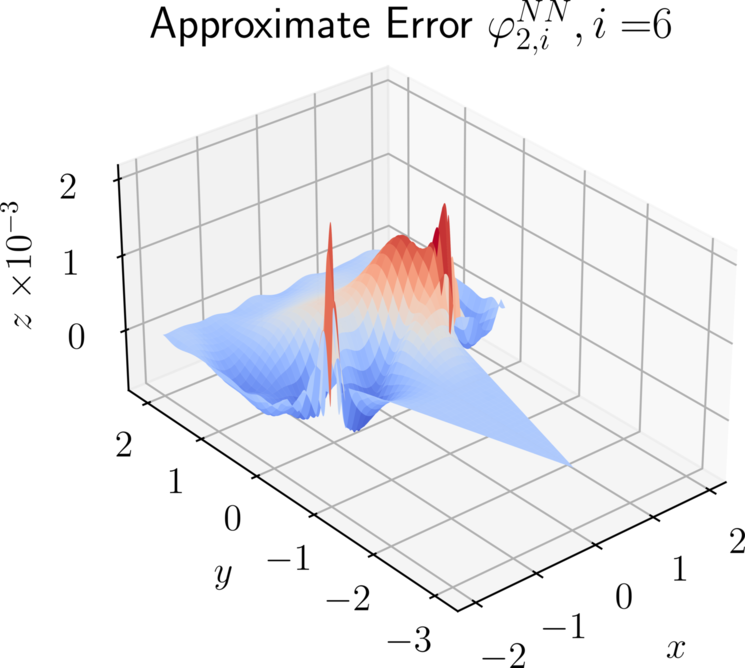}
			
			\includegraphics[width=1.3in]{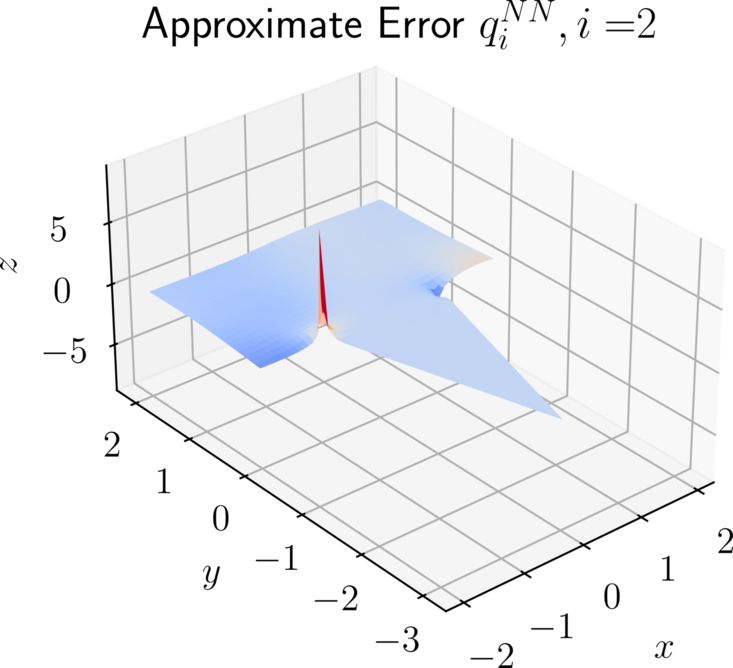}
			\quad
			\includegraphics[width=1.3in]{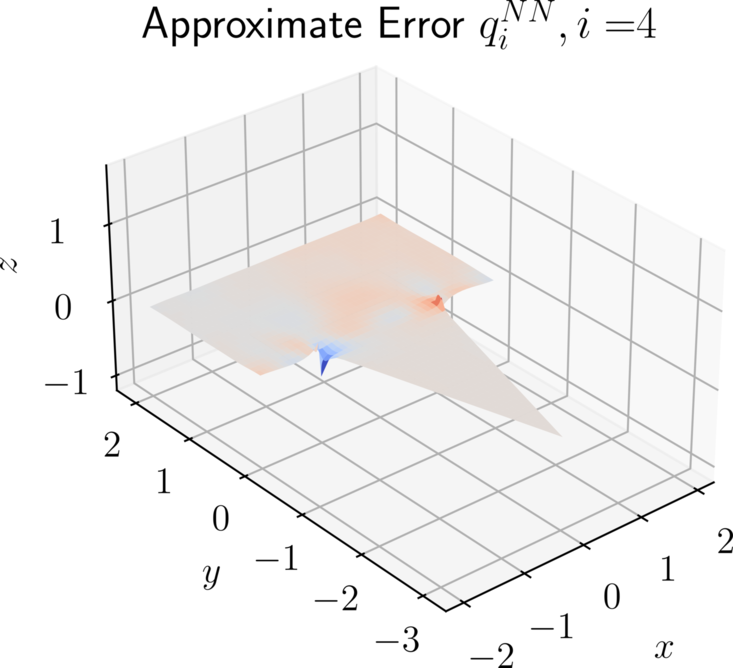}
			\quad
			\includegraphics[width=1.3in]{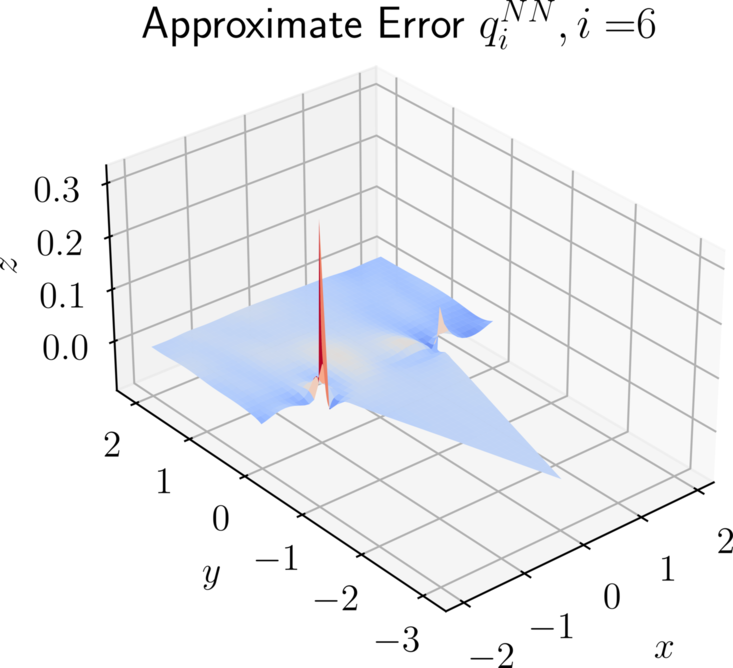}
			
			\caption{Example \ref{ex:triangular}: velocity and pressure basis functions with $\beta=5/3$ in \eqref{eq:stokes bilinear} and basis functions learned from extended neural network with knowledge-based functions.}
			\label{fig:stokes beta53 basis}
		\end{figure}
		

		\begin{figure}[t!]
			\centering
			\includegraphics[width=1.3in]{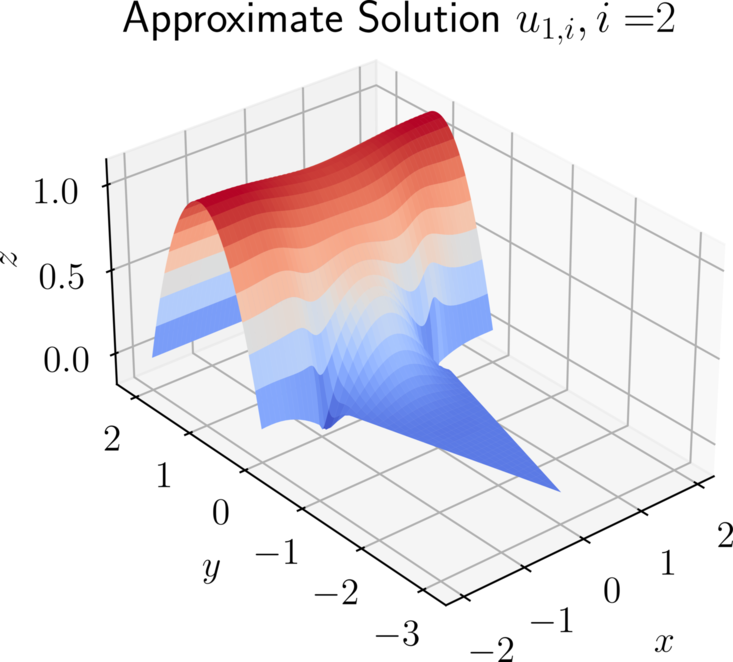}
			\quad
			\includegraphics[width=1.3in]{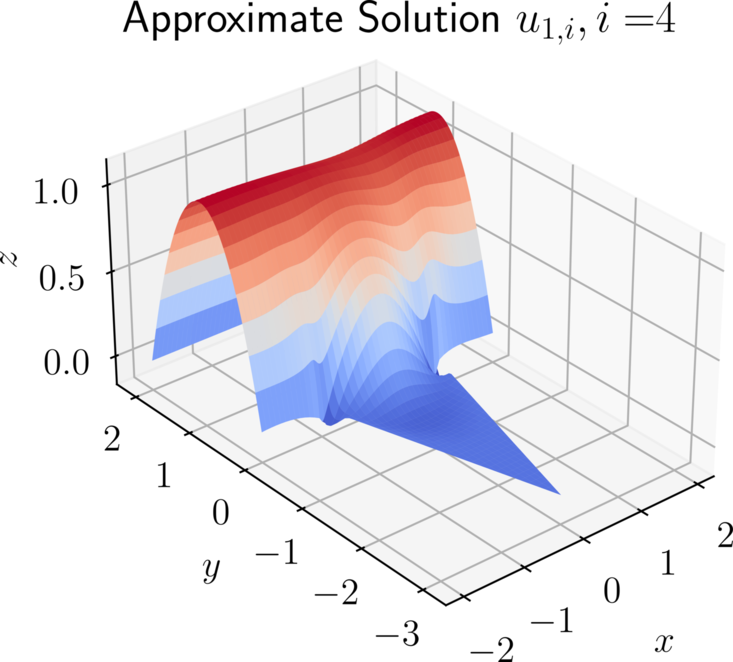}
			\quad
			\includegraphics[width=1.3in]{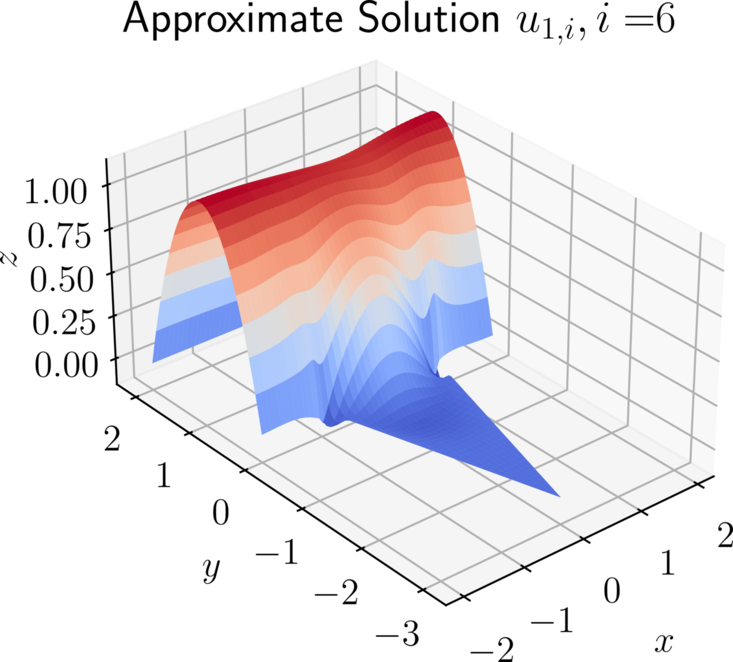}
			
			\includegraphics[width=1.3in]{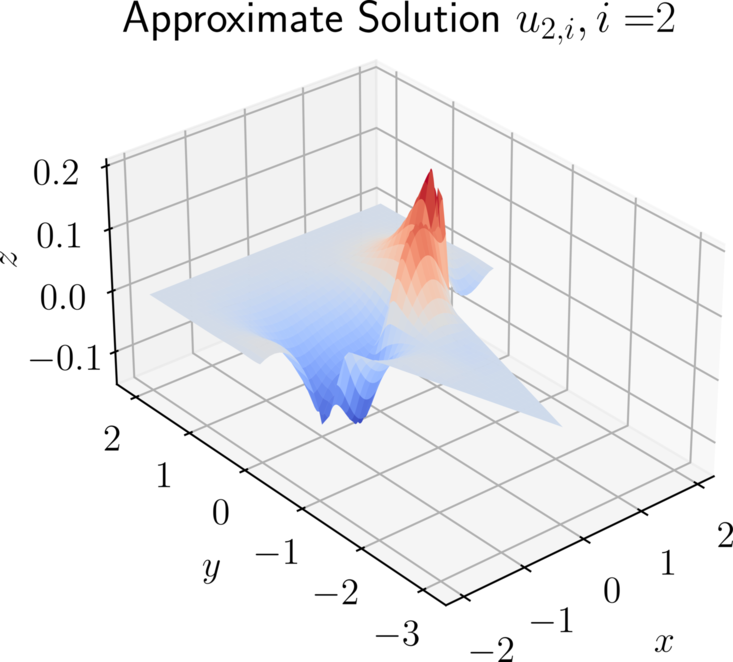}
			\quad
			\includegraphics[width=1.3in]{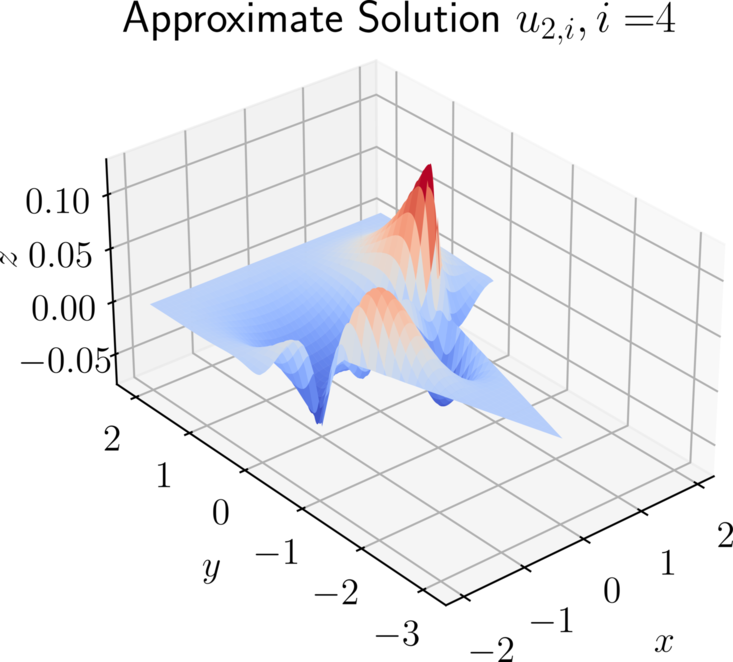}
			\quad
			\includegraphics[width=1.3in]{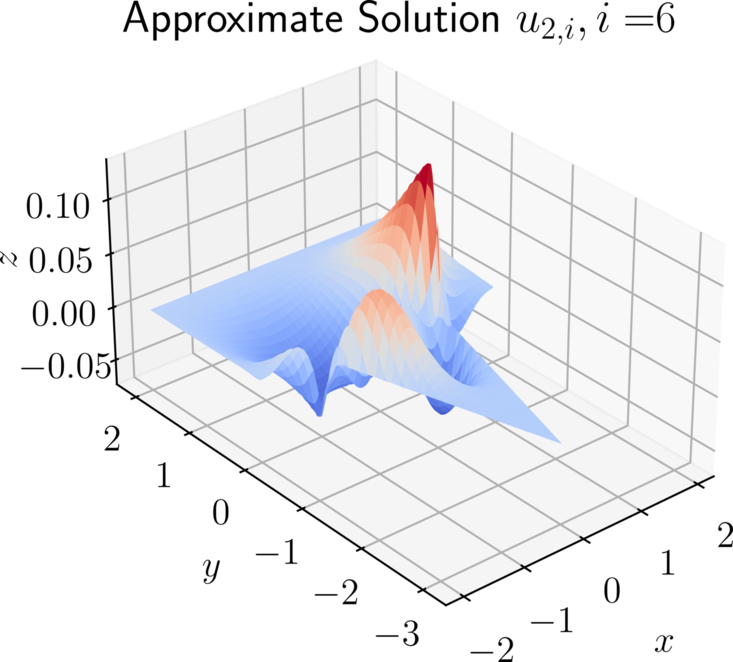}
			
			\includegraphics[width=1.3in]{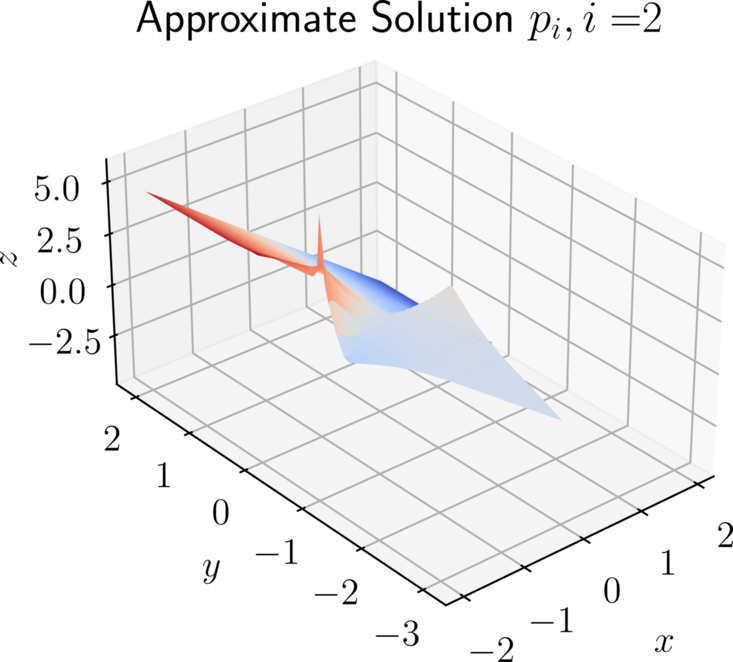}
			\quad
			\includegraphics[width=1.3in]{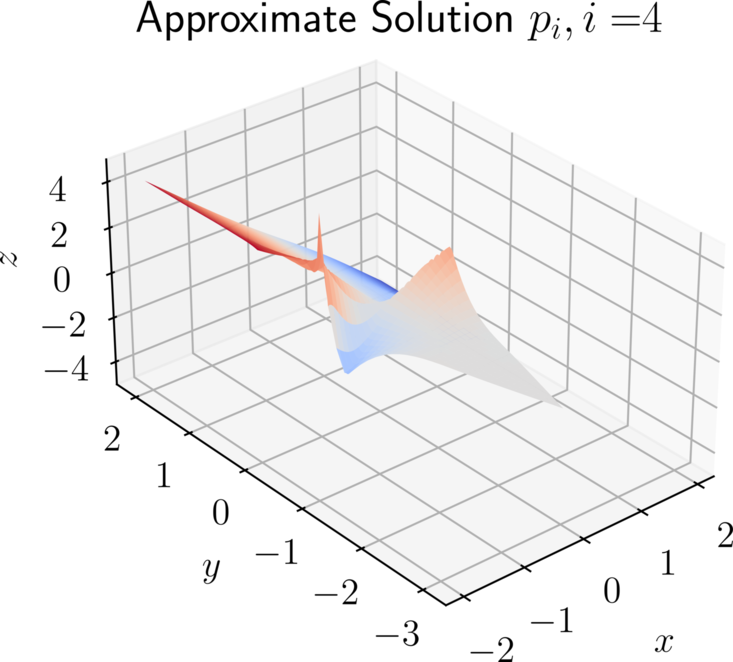}
			\quad
			\includegraphics[width=1.3in]{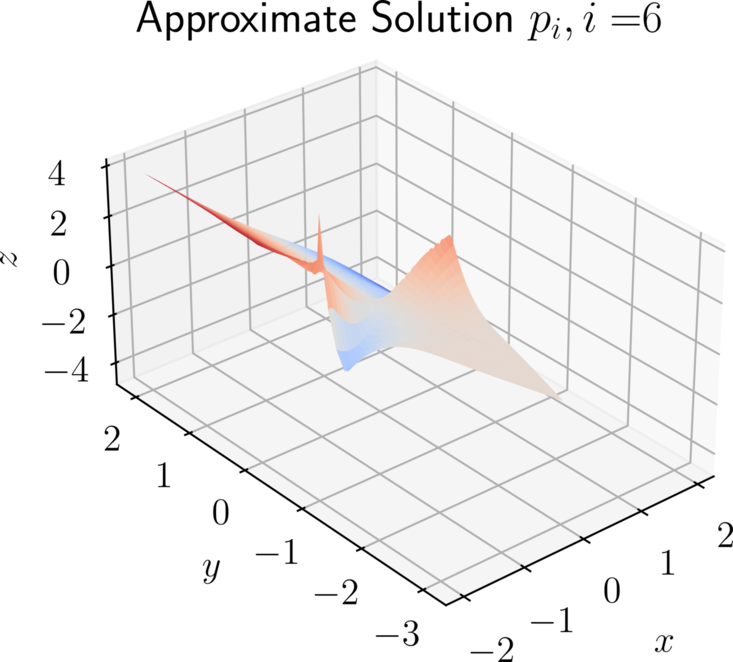}
			
			\caption{Example \ref{ex:triangular}: velocity and pressure basis functions with $\beta=4$ in \eqref{eq:stokes bilinear} and basis functions learned from extended neural network with knowledge-based functions.}
			\label{fig:stokes beta4}
		\end{figure}
	\end{example}

    \section{Extended Galerkin neural networks with learning of knowledge-based functions} \label{sec:eigenvalue training}

    In the numerical examples in Section \ref{sec:applications}, we approximated the solutions from the set $V_{\mathbf{n},L}^{\sigma} \oplus V_{\mathcal{M}}^{\Psi}(\lambda)$ by utilizing exact known values of $\lambda$ when including these functions in the approximation space. However, in many applications some information about the structure of $\Psi$ may be known while the parameters $\lambda$ may be unknown or difficult to compute exactly. In this section, we demonstrate how extended Galerkin neural networks may be applied to the previously studied examples in order to approximate $\lambda$.

    \subsection{Poisson equation}
    We revisit Example \ref{ex:poisson Lshaped} in which the knowledge-based functions took the form
    \begin{align*}
        \Psi(r,\theta;\mu) = r^{\mu}\sin(\mu\theta).
    \end{align*}

    \noindent The exact eigenvalues are given by $\lambda_{j}=2j/3$. We focus only on the smallest value of $\lambda$ which describes the asymptotic behavior of $\Psi$ and note that the feedforward part of the neural network, $V_{\mathbf{n},L}^{\sigma}$, is capable of approximating the smoother behavior of $\Phi$ for larger values of $\lambda$.
    \begin{example} \label{ex:poisson Lshaped training}
        We again consider the L-shaped domain $\Omega = (-1,1)^{2}\backslash(-1,0)^{2}$ with data $f=1$ and homogeneous Dirichlet boundary conditions. Each Galerkin neural network basis function is obtained according to \eqref{eq:singular basis fn} with $\mu$ being a trainable parameter. The initial value of $\mu$ is drawn from the uniform distribution $\mathcal{U}(0,1)$. Full details of the hyperparameters for this example may be found in Appendix \ref{app:examples}. Figure \ref{fig:poisson Lshaped3} shows how $\mu_{1}$ changes value during the backpropagation step of training over 20 independent runs. In each case $\mu_{1}$ moves toward the true value of $2/3$. 
    \end{example}

    \begin{figure}[t!]
        \centering
        \includegraphics[width=2.1in]{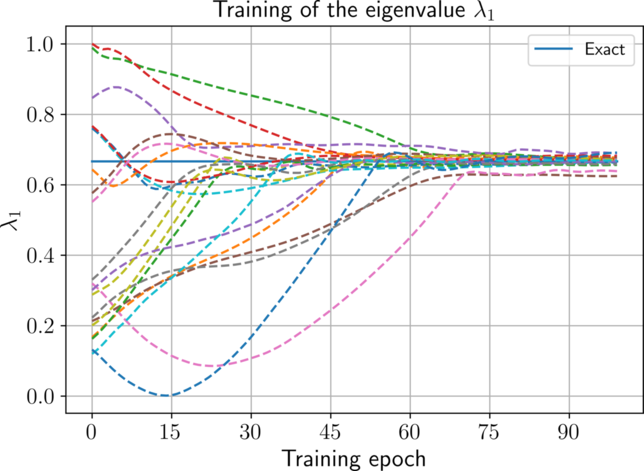}
        \caption{Example \ref{ex:poisson Lshaped training}:  learning of $\lambda_{1}$ and over several runs with initialization $\mu \sim \mathcal{U}(0,1)$.}
        \label{fig:poisson Lshaped3}
    \end{figure}

    \subsection{Stokes flow around a non-convex corner}

    \begin{example} \label{ex:stokes learning1}
        Analogous to Example \ref{ex:poisson Lshaped training}, we demonstrate how to learn the knowledge-based functions for the Stokes flow around a non-convex corner. We consider the circular sector domain $\Omega = \{(r,\theta) \;:\; 0 < r < 1, \;0 < \theta < \pi + \alpha^{(0)}\}$, $\alpha^{(0)} = \arccos(1/\sqrt{10})$ with $\mathbf{f} = \mathbf{0}$ and $g=0$. The boundary condition for the velocity is given by
        \begin{align*}
            \mathbf{u}(\theta) = \begin{bmatrix}
                -\frac{1}{2\alpha^{(0)}}\theta^{2} + \theta\\
                \frac{-2\alpha^{(0)}-1}{4(\alpha^{(0)})^{2}}\theta^{3} + \theta^{2} + \theta
            \end{bmatrix}.
        \end{align*}

        \noindent This example may be viewed as a localized BVP with artificial boundary conditions.

        \begin{figure}
            \centering
            \includegraphics[width=2.1in]{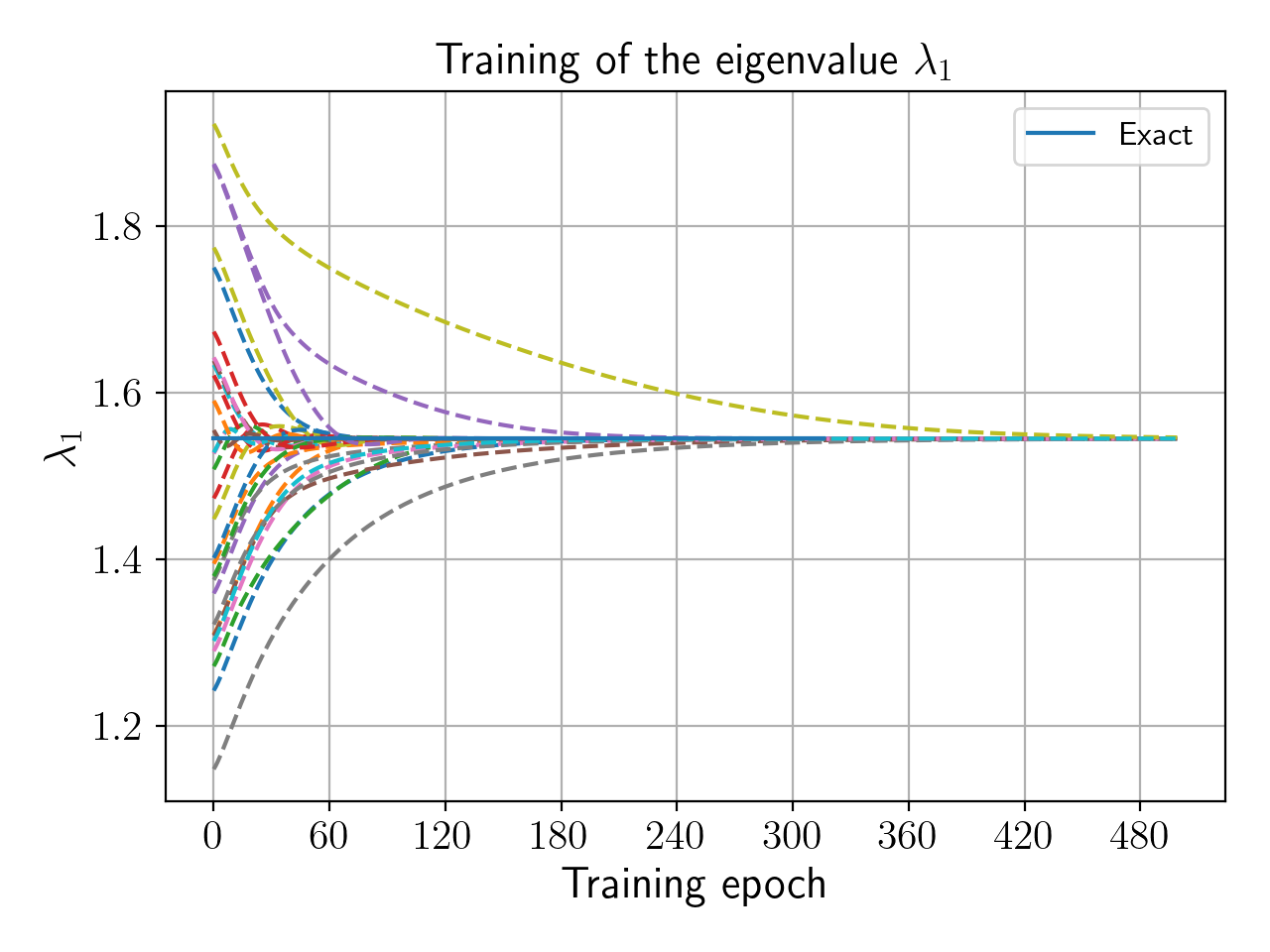}
            \caption{Eigenvalue training progress for Example \ref{ex:stokes learning1} with initialization $\mu \sim \mathcal{U}(\lambda_{\text{exact}}-1/3, \lambda_{\text{exact}} + 1/3)$.}
            \label{fig:stokes learning}
        \end{figure}

        We approximate the solution to this example using a feedforward network augmented with the knowledge-based function $\Psi(r,\theta;\mu) = r^{\mu}\Psi_{\mu}(\theta)$, where $\Psi_{\mu}$ is given in Appendix \ref{app:stokes} for the velocity and pressure. We initialize $\mu$ away from its true value as in Example \ref{ex:poisson Lshaped training} in order to demonstrate that the true eigenvalue of the singular part of the Stokes operator may also be learned by the neural network. Importantly, the singular-behavior of the solution at the origin should be localized with a cutoff function, namely in \eqref{eq:singular NN} we take
        \begin{align*}
            V_{m}^{\Phi} := \left\{ v \;:\; v(x) = \sum_{i=1}^{m} d_{i}\chi(r)\Phi(x;\mu_{i}), \;\mathbf{d} \in \mathbb{R}^{m}, \;\lambda_{j} \in \mathcal{I}_{\Phi} \right\},
        \end{align*}

        \noindent where $\chi$ is a $C^{2}$ cutoff function satisfying $\chi(r) = 1$ for $r < r_{0}$, $\chi(r) = 0$ for $r > r_{1}$. Full details of the hyperparameters for this example may be found in Appendix \ref{app:examples}.

        Figure \ref{fig:stokes learning} shows the training progress of $\lambda_{1}$ over 30 trials with $\mu$ initialized according to the uniform distribution $\mathcal{U}(\lambda_{\text{exact}}-1/3, \lambda_{\text{exact}} + 1/3)$. We observe that in each trial, $\mu$ converges to the true value. 
    \end{example}
    
    \subsection{Stokes flow in a convex corner induced by disturbance at a large distance}

    \begin{example} \label{ex:stokes learning2}
        We now demonstrate how to learn the knowledge-based functions for the Stokes flow in a convex corner induced by flow at a large distance. The localized domain $\Omega$ is the triangular wedge with vertices $(1,0)$, $(-1,0)$, and $(0,-3)$. The boundary condition for the velocity is a homogeneous Dirichlet condition along the diagonal edges of the wedge and $\mathbf{u}(x,y) = [(1-x)(1+x), 0]^{T}$ along the top edge. We again specify $\mathbf{f} = \mathbf{0}$ and $g=0$.
        \begin{figure}[t!]
            \centering
            \includegraphics[width=2.0in]{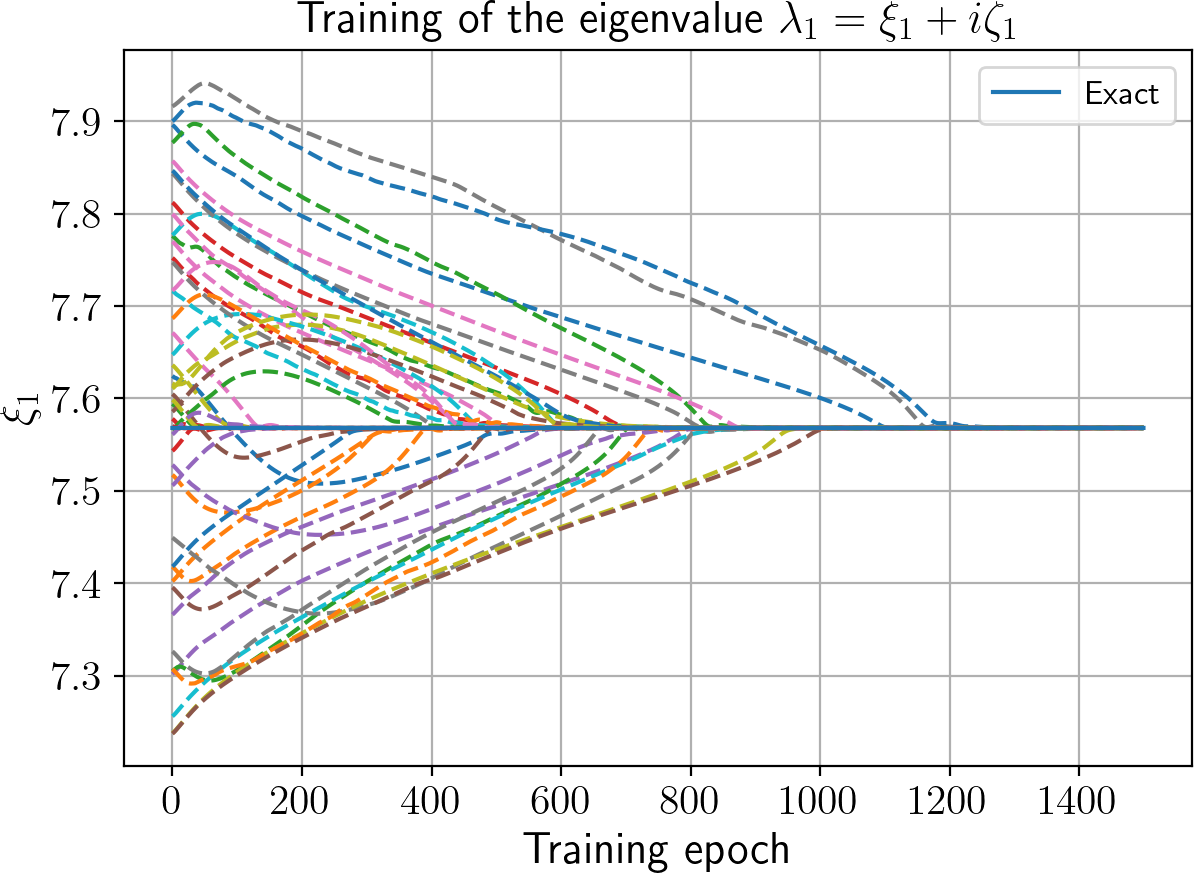}
            \quad
            \includegraphics[width=2.0in]{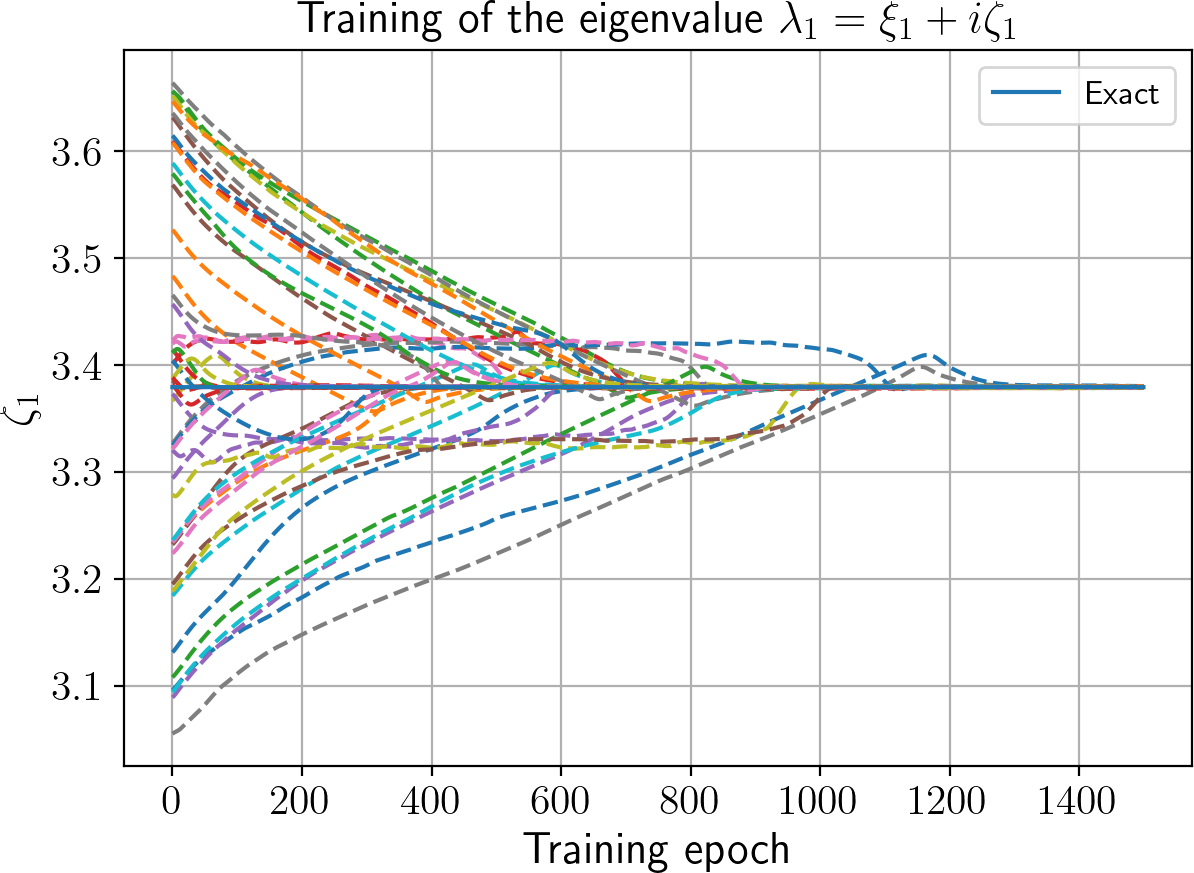}
            \caption{Eigenvalue training progress for Example \ref{ex:stokes learning2} with initialization $\xi \sim \mathcal{U}(\xi_{\text{exact}}-1/2, \xi_{\text{exact}} + 1/2)$ and $\zeta \sim \mathcal{U}(\zeta_{\text{exact}}-1/2, \zeta_{\text{exact}} + 1/2)$.}
            \label{fig:stokes learning eddy}
        \end{figure}
        
        We approximate the solution to this example using a feedforward network augmented with knowledge-based function $\Psi(r,\theta;\mu) = \mathfrak{Re}[r^{\mu}\Psi_{\mu}(\theta)]$. Here, $\mu$ is assumed to be an unknown complex number of the form $\mu = \xi + i\zeta$. A cutoff function is again used to isolate the behavior of $\Psi$ to a region around the vertex $(0,-3)$. We initialize $\xi$ and $\zeta$ away from their true values in order to demonstrate that the correct behavior may be learned by the augmented feedforward neural network. Full details of the hyperparameters for this example may be found in Appendix \ref{app:examples}.

        Figure \ref{fig:stokes learning eddy} shows the training progress of $\lambda$ over 50 trials with $\xi$ and $\zeta$ initialized according to the uniform distributions $\mathcal{U}(\xi_{\text{exact}}-1/2, \xi_{\text{exact}} + 1/2)$ and $\mathcal{U}(\zeta_{\text{exact}}-1/2, \zeta_{\text{exact}} + 1/2)$, respectively. Remarkably, we observe that the correct eddy behavior is obtained and the values of $\xi$ and $\zeta$ are learned to within about $0.005\%$.
    \end{example}

	\section{Conclusions} \label{sec:conclusions}

	We have presented the extended Galerkin neural network framework for approximating a general class of boundary value problems, including those with low-regularity solutions. The main features of this framework are:
    \begin{itemize}
        \item The presentation of well-posed least squares variational formulations for arbitrary boundary value problems that are well-suited to approximation by Galerkin neural networks and do not require conversion to first-order systems. We demonstrate that these formulations are dictated by the properties of the PDE under consideration and cannot be chosen arbitrarily by the practitioner.

        \item Extended feedforward neural network architectures that systematically incorporate knowledge-based functions into the approximation space and reduce the width of the network required to attain a given accuracy.

        \item A rigorous convergence theory which is borne out in the numerical experiments whereas other approaches lack sound theoretical footing and in practice tend to exhibit non-monotonic reduction of the error.
    \end{itemize}

	\section{Acknowledgments} \label{sec:acknowledgments}

	J.D. thanks the National Science Foundation Graduate Research Fellowship for its financial support under Grant No. 1644760. This work was performed under the auspices	of the U.S. Department of Energy by Lawrence Livermore National Laboratory under Contract DE-AC52-07NA27344; LLNL-JRNL-859911. M.A. gratefully acknowledges the support of National Science Foundation Grant No. DMS-2324364.

      \appendix	 
	 \section{Knowledge-based functions for Stokes flow} \label{app:stokes}
	 
	Given a vertex $x^{(i)}$ of $\Omega$ with interior angle $\alpha^{(i)}$, the solution of the Stokes flow in the vicinity of $x^{(i)}$ is may be written into terms of the eigenfunctions of the operator pencil of the Stokes operator \cite{mazya}. It is simpler to first consider the streamfunction $\psi$, which takes the form
	\begin{align*}
		\psi(r,\theta) = \sum_{\lambda_{n}} c_{\lambda_{n}} \mathfrak{Re}\left[r^{\lambda_{n}}\Psi_{\lambda_{n}}(\theta)\right],
	\end{align*}
	\noindent where for $\lambda_{n} \neq 0$, $\Psi_{\lambda_{n}}$ is given by the general form
        \begin{align*}
            \Psi_{\lambda_{n}}(\theta) := A_{1}\cos((\lambda_{n}-2)\theta) + A_{2}\cos(\lambda_{n}\theta) + A_{3}\sin((\lambda_{n}-2)\theta) + A_{4}\sin(\lambda_{n}\theta).
        \end{align*}

        \noindent The coefficients $A_{i}$ are chosen to satisfy the boundary conditions. The velocity may be obtained from the streamfunction with the relation $\mathbf{u} = \nabla^{\perp}\psi$ while the pressure may be obtained from the momentum equation $\nabla p = \Delta\mathbf{u}$.

        \subsection{No-slip boundary conditions}
        We shall first consider the simple case of no-slip boundary conditions, i.e. $\mathbf{u}\cdot\mathbf{n} = \partial\psi/\partial t = 0$. In this case, we have
        \begin{align} \label{eq:streamfunction psi}
		\Psi_{\lambda_{n}}(\theta) := \cos(\lambda_{n}\alpha)\cos((\lambda_{n}-2)\theta) - \cos((\lambda_{n}-2)\alpha)\cos(\lambda_{n}\theta).
	\end{align}

	\noindent In the case of complex eigenvalues, the streamfunction produces a sequence of infinitely cascading eddies, a qualitative discussion of which may be found in \cite{moffatt}. To the best of our knowledge, an explicit derivation of the full corresponding eigenfunctions of the eddies for the Stokes velocity and pressure does not exist in the literature. Thus, the purpose of this section is to state the velocities and pressure corresponding to the streamfunction $r^{\lambda_{n}}\Psi_{\lambda_{n}}(\theta)$ assuming complex eigenvalues $\lambda_{n} = \xi_{n} + i\zeta_{n}$. In this case, we have
     \begin{align*}
         \mathfrak{Re}\left[ r^{\lambda_{n}}\Psi_{\lambda_{n}}(\theta) \right] = \mathfrak{Re}[r^{\lambda_{n}}] \mathfrak{Re}[\Psi_{\lambda_{n}}(\theta)] - \mathfrak{Im}[r^{\lambda_{n}}] \mathfrak{Im}[\Psi_{\lambda_{n}}(\theta)]
     \end{align*}

     \noindent from which the velocity components may be derived:
	 \begin{align}
	 	\begin{bmatrix}
			u_{1}(r,\theta)\\
			u_{2}(r,\theta)
		\end{bmatrix}& = \sum_{\lambda_{n}} \begin{bmatrix}
			\frac{\partial}{\partial r}\mathfrak{Re}[r^{\lambda_{n}}] \mathfrak{Re}[\Psi_{\lambda_{n}}(\theta)]\sin{\theta} + \mathfrak{Re}[r^{\lambda_{n}}]\frac{\partial}{\partial\theta}\mathfrak{Re}[\Psi_{\lambda_{n}}(\theta)]\frac{\cos{\theta}}{r}\\
			-\frac{\partial}{\partial r}\mathfrak{Re}[r^{\lambda_{n}}] \mathfrak{Re}[\Psi_{\lambda_{n}}(\theta)]\cos{\theta} + \mathfrak{Re}[r^{\lambda_{n}}]\frac{\partial}{\partial\theta}\mathfrak{Re}[\Psi_{\lambda_{n}}(\theta)]\frac{\sin{\theta}}{r}
		\end{bmatrix}\notag\\
		&\;\;\;+\sum_{\lambda_{n}} \begin{bmatrix}
			- \frac{\partial}{\partial r}\mathfrak{Im}[r^{\lambda_{n}}]\mathfrak{Im}[\Psi_{\lambda_{n}}(\theta)]\sin{\theta} - \mathfrak{Im}[r^{\lambda_{n}}]\frac{\partial}{\partial\theta}\mathfrak{Im}[\Psi_{\lambda_{n}}(\theta)]\frac{\cos{\theta}}{r}\\
			\frac{\partial}{\partial r}\mathfrak{Im}[r^{\lambda_{n}}]\mathfrak{Im}[\Psi_{\lambda_{n}}(\theta)]\cos{\theta} - \mathfrak{Im}[r^{\lambda_{n}}]\frac{\partial}{\partial\theta}\mathfrak{Im}[\Psi_{\lambda_{n}}(\theta)]\frac{\sin{\theta}}{r}
		\end{bmatrix}\label{eq:moffatt complex}.
	 \end{align}

	 \noindent Here, 
	 \begin{align*}
	 	\mathfrak{Re}[r^{\lambda_{n}}] &= r^{\xi_{n}}\cos(\log(r)\zeta_{n}), \;\;\;\mathfrak{Im}[r^{\lambda_{n}}] = r^{\xi_{n}}\sin(\log(r)\zeta_{n})\\
		\mathfrak{Re}[\Psi_{\lambda_{n}}(\theta)] &= \cos(\xi_{n}\alpha^{(i)})\cos((\xi_{n}-2)\theta)\cosh(\zeta_{n}\alpha^{(i)})\cosh(\zeta_{n}\theta)\\
		&\;\;\;+\sin(\xi_{n}\alpha^{(i)})\sin((\xi_{n}-2)\theta)\sinh(\zeta_{n}\alpha^{(i)})\sinh(\zeta_{n}\theta)\\
		&\;\;\;-\cos((\xi_{n}-2)\alpha^{(i)})\cos((\xi_{n}-2)\theta)\cosh(\zeta_{n}\alpha^{(i)})\cosh(\zeta_{n}\theta)\\
		&\;\;\;+\sin((\xi_{n}-2)\alpha^{(i)})\sin(\xi_{n}\theta)\sinh(\zeta_{n}\alpha^{(i)})\sinh(\zeta_{n}\theta)\\
		\mathfrak{Im}[\Psi_{\lambda_{n}}(\theta)] &= -\cos(\xi_{n}\alpha^{(i)})\sin((\xi_{n}-2)\theta) \cosh(\zeta_{n}\alpha^{(i)})\sinh(\zeta_{n}\theta)\\
		&\;\;\;-\sin(\xi_{n}\alpha^{(i)})\cos((\xi_{n}-2)\theta) \sinh(\zeta_{n}\alpha^{(i)}) \cosh(\zeta_{n}\theta)\\
		&\;\;\;+\cos((\xi_{n}-2)\alpha^{(i)}) \sin(\xi_{n}\theta) \cosh(\zeta_{n}\alpha^{(i)}) \sinh(\zeta_{n}\theta)\\
		&\;\;\;+\sin((\xi_{n}-2)\alpha^{(i)}) \cos(\xi_{n}\theta) \sinh(\zeta_{n}\alpha^{(i)}) \cosh(\zeta_{n}\theta)).
	 \end{align*}
	 
    \noindent We note that \eqref{eq:moffatt complex} has been calculated to exactly satisfy a no-slip condition along $\theta = \pm\alpha^{(i)}$. The transformation $\theta \mapsto \theta - \phi$, $\phi \in \mathbb{R}$ should be taken if necessary depending on specific domain geometry. 

    Finally, to obtain the pressure corresponding to the streamfunction $\mathfrak{Re}[r^{\lambda_{n}}\Psi_{\lambda_{n}}(\theta)]$, it is straightforward to calculate $\partial p/\partial r$ and $\partial p/\partial\theta$ from the momentum equations in polar coordinates \cite{batchelor1967introduction}:
    \begin{align}
        \begin{dcases}
            -\left(\frac{\partial^{2}u_{r}}{\partial r^{2}} + \frac{1}{r}\frac{\partial u_{r}}{\partial r} + \frac{1}{r^{2}}\left( \frac{\partial^{2}u_{r}}{\partial\theta^{2}} - 2\frac{\partial u_{\theta}}{\partial\theta} - u_{r} \right)\right) + \frac{\partial p}{\partial r} = 0\\
            -\left(\frac{\partial^{2}u_{\theta}}{\partial r^{2}} + \frac{1}{r}\frac{\partial u_{\theta}}{\partial r} + \frac{1}{r^{2}}\left( \frac{\partial^{2}u_{\theta}}{\partial\theta^{2}} - 2\frac{\partial u_{r}}{\partial\theta} - u_{\theta} \right)\right) + \frac{1}{r}\frac{\partial p}{\partial \theta} = 0
        \end{dcases} \label{eq:momentum polar}
    \end{align}

    \noindent In \eqref{eq:momentum polar}, $u_{r}$ and $u_{\theta}$ are the radial and azimuthal components of the velocity given by $u_{r}(r,\theta) := r^{-1}\partial\psi/\partial\theta$ and $u_{\theta}(r,\theta) := -\partial\psi/\partial r$, respectively. Integrating $\partial p/\partial r$ and $\partial p/\partial\theta$ and equating constants yields
    \begin{align}
        p(r,\theta) &= 4r^{\xi_{n}-2}\biggl( \cos(\xi_{n}\alpha^{(i)})\cosh(\zeta_{n}\alpha^{(i)})\left( -\sin((\xi_{n}-2)\theta)\cosh(\zeta_{n}\theta)A(r) \right. \notag\\
        &\hspace{50mm}\left.\left.- \cos((\xi_{n}-2)\theta)\sinh(\zeta_{n}\theta)B(r) \right) + \right.\notag\\
        &\hspace{18mm}\left. \sin(\xi_{n}\alpha^{(i)})\sinh(\zeta_{n}\alpha^{(i)})\left( \sin((\xi_{n}-2)\theta)\cosh(\zeta_{n}\theta) B(r) \right.\right.\notag\\
        &\hspace{50mm}\left.- \cos((\xi_{n}-2)\theta)\sinh(\zeta_{n}\theta) A(r) \right) \biggr), \label{eq:moffatt pressure}\\
        A(r) &:= -(\xi_{n}-1)\cos(\log(r)\zeta_{n}) + \zeta_{n}\sin(\log(r)\zeta_{n})\notag\\
        B(r) &:= \zeta_{n}\cos(\log(r)\zeta_{n}) + (\xi_{n}-1)\sin(\log(r)\zeta_{n}).\notag
    \end{align}

    \noindent In the special case when $\lambda_{n} \in \mathbb{R}$, we simply take $\zeta_{n} = 0$ to arrive at 
    \begin{align} \label{eq:no slip eigenfunctions}
	\begin{bmatrix}
			u_{1}(r,\theta)\\
			u_{2}(r,\theta)\\
			p(r,\theta) 
		\end{bmatrix} &= \sum_{\lambda_{n}} \begin{bmatrix} \lambda_{n}r^{\lambda_{n}-1}\mathfrak{Re}[\Psi_{\lambda_{n}}(\theta)]\sin{\theta} - r^{\lambda_{n}}\frac{\partial}{\partial\theta}\mathfrak{Re}[\Psi_{\lambda_{n}}(\theta)]\cos{\theta}\\
        -\lambda_{n}r^{\lambda_{n}-1}\mathfrak{Re}[\Psi_{\lambda_{n}}(\theta)]\cos{\theta} + r^{\lambda_{n}-1}\frac{\partial}{\partial\theta}\mathfrak{Re}[\Psi_{\lambda_{n}}(\theta)]\sin{\theta}\\
        4(\lambda_{n}-1)r^{\lambda_{n}-2} \cos(\lambda_{n}\alpha^{(i)}) \sin((\lambda_{n}-2)\theta)
        \end{bmatrix}.
    \end{align}


    \subsection{Homogeneous Dirichlet boundary conditions}
    \noindent For the boundary homogeneous Dirichlet boundary condition $\mathbf{u} = \mathbf{0}$ used in the examples in Section \ref{sec:applications}, we leave the constants $A_{i}$ as unknowns to be determined during the least squares step \eqref{eq:lstsq solve}. In this case, the velocity corresponding to the streamfunction $\mathfrak{Re}[r^{\lambda_{n}}\Psi_{\lambda_{n}}(\theta)]$ consists of four terms:
	 \begin{align} \label{eq:velocity dirichlet}
	 	\begin{bmatrix}
			u_{1}(r,\theta)\\
			u_{2}(r,\theta)
		\end{bmatrix}& = \sum_{i=1}^{4}A_{i}\left\{\sum_{\lambda_{n}} \begin{bmatrix}
			\frac{\partial}{\partial r}\mathfrak{Re}[r^{\lambda_{n}}] \mathfrak{Re}[\Psi_{i,\lambda_{n}}(\theta)]\sin{\theta} + \mathfrak{Re}[r^{\lambda_{n}}]\frac{\partial}{\partial\theta}\mathfrak{Re}[\Psi_{i,\lambda_{n}}(\theta)]\frac{\cos{\theta}}{r}\\
			-\frac{\partial}{\partial r}\mathfrak{Re}[r^{\lambda_{n}}] \mathfrak{Re}[\Psi_{i,\lambda_{n}}(\theta)]\cos{\theta} + \mathfrak{Re}[r^{\lambda_{n}}]\frac{\partial}{\partial\theta}\mathfrak{Re}[\Psi_{i,\lambda_{n}}(\theta)]\frac{\sin{\theta}}{r}
		\end{bmatrix}\right.\notag\\
		&\left.\;\;\;+\sum_{\lambda_{n}} \begin{bmatrix}
			- \frac{\partial}{\partial r}\mathfrak{Im}[r^{\lambda_{n}}]\mathfrak{Im}[\Psi_{i,\lambda_{n}}(\theta)]\sin{\theta} - \mathfrak{Im}[r^{\lambda_{n}}]\frac{\partial}{\partial\theta}\mathfrak{Im}[\Psi_{i,\lambda_{n}}(\theta)]\frac{\cos{\theta}}{r}\\
			\frac{\partial}{\partial r}\mathfrak{Im}[r^{\lambda_{n}}]\mathfrak{Im}[\Psi_{i,\lambda_{n}}(\theta)]\cos{\theta} - \mathfrak{Im}[r^{\lambda_{n}}]\frac{\partial}{\partial\theta}\mathfrak{Im}[\Psi_{i,\lambda_{n}}(\theta)]\frac{\sin{\theta}}{r}
		\end{bmatrix}\right\}.
	 \end{align}

    \noindent For ease of notation, we use $\Psi_{i,\lambda_{n}}$ to denote the four component eigenfunctions
    \begin{align*}
        \Psi_{1,\lambda_{n}}(\theta) &= \cos((\lambda_{n}-2)\theta), \;\;\;\Psi_{2,\lambda_{n}}(\theta) = \cos(\lambda_{n}\theta),\\
        \Psi_{3,\lambda_{n}}(\theta) &= \sin((\lambda_{n}-2)\theta), \;\;\;\;\Psi_{4,\lambda_{n}}(\theta) = \sin(\lambda_{n}\theta).
    \end{align*}

    \noindent The real and complex components of $\Psi_{i,\lambda_{n}}$ are given by
    \begin{align*}
        \mathfrak{Re}[\Psi_{1,\lambda_{n}}(\theta)] &= \cos((\xi_{n}-2)\theta) \cosh(\zeta_{n}\theta), \;\;\;\mathfrak{Im}[\Psi_{1,\lambda_{n}}(\theta)] = -\sin((\xi_{n}-2)\theta) \sinh(\zeta_{n}\theta),\\
        \mathfrak{Re}[\Psi_{2,\lambda_{n}}(\theta)] &= \cos(\xi_{n}\theta) \cosh(\zeta_{n}\theta), \;\;\;\;\;\;\;\;\;\;\;\;\mathfrak{Im}[\Psi_{2,\lambda_{n}}(\theta)] = -\sin(\xi_{n}\theta) \sinh(\zeta_{n}\theta),\\
        \mathfrak{Re}[\Psi_{3,\lambda_{n}}(\theta)] &= \sin((\xi_{n}-2)\theta) \cosh(\zeta_{n}\theta), \;\;\;\mathfrak{Im}[\Psi_{3,\lambda_{n}}(\theta)] = \cos((\xi_{n}-2)\theta) \sinh(\zeta_{n}\theta),\\
        \mathfrak{Re}[\Psi_{4,\lambda_{n}}(\theta)] &= \sin(\xi_{n}\theta) \cosh(\zeta_{n}\theta), \;\;\;\;\;\;\;\;\;\;\;\;\mathfrak{Im}[\Psi_{4,\lambda_{n}}(\theta)] = \cos(\xi_{n}\theta) \sinh(\zeta_{n}\theta).
    \end{align*}

    The pressure corresponding to $\mathfrak{Re}[r^{\lambda_{n}}\Psi_{\lambda_{n}}(\theta)]$ is given by
    \begin{align} \label{eq:pressure dirichlet}
        p(r,\theta) = \sum_{i=1}^{4} A_{i} P_{i,\lambda_{n}}(r,\theta),
    \end{align}
    
    \noindent where the four component eigenfunctions $P_{i,\lambda_{n}}$ are 
    \begin{align*} 
        P_{1,\lambda_{n}}(r,\theta) &= 4r^{\xi_{n}-2} \left( \sin((\xi_{n}-2)\theta) \cosh(\zeta_{n}\theta) A(r) + \cos((\xi_{n}-2)\theta) \sinh(\zeta_{n}\theta) B(r) \right)\notag\\
        P_{3,\lambda_{n}}(r,\theta) &= 4r^{\xi_{n}-2} \left( -\cos((\xi_{n}-2)\theta) \cosh(\zeta_{n}\theta)A(r) + \sin((\xi_{n}-2)\theta) \sinh(\zeta_{n}\theta) B(r) \right)\notag
    \end{align*}

    \noindent with $P_{2,\lambda_{n}}(r,\theta) = P_{4,\lambda_{n}}(r,\theta) = 0$ and $A(r)$, $B(r)$ as in \eqref{eq:moffatt pressure}.

	 
	 

    \section{Hyperparameters for numerical examples} \label{app:examples}

    \subsection{Example \ref{ex:cont coercivity}} \label{app:first example}
    \hphantom{4mm}
    \vspace{-3mm}
    \begin{table}[H]
        \small
        \centering
        \begin{tabular}{||c||c||}
            \hline
            \textbf{Parameter} & \\
            \hline \hline
            $u_{0}(x)$ & 0\\
            \hline
            $\mathbf{n}^{(i)}$ & $20\cdot 2^{i-1}$\\
            \hline
            $L_{i}$ & 1\\
            \hline
            $\sigma_{i}$ & $\tanh((1 + \frac{i}{4})t)$\\
            \hline
            training data & \begin{tabular}{@{}c@{}}tensor-product Gauss-Legendre $128 \times 128$ in $\Omega$;\\ left Riemann sum 256 on $\partial\Omega$\end{tabular}\\
            \hline
            learning rate & $10^{-3}/1.1^{i-1}$\\
            \hline
            $\delta$ & Variable. See Example \ref{ex:cont coercivity} for details. \\
            \hline
        \end{tabular}
        \label{tab:ex1}
    \end{table}

    \subsection{Examples \ref{ex:poisson rlambda}-\ref{ex:poisson Lshaped}}
    \hphantom{0mm}
    \vspace{-3mm}
    \begin{table}[H]
        \small
        \centering
        \begin{tabular}{||c||c||}
            \hline 
            \textbf{Parameter} & \\
            \hline \hline
            $u_{0}(x)$ & 0\\
            \hline
            $\mathbf{n}^{(i)}$ & $20\cdot 2^{i-1}$\\
            \hline
            $L_{i}$ & 1\\
            \hline
            $\sigma_{i}$ & $\tanh((1 + \frac{i}{4})t)$\\
            \hline
            $\mathcal{M}$\tablefootnote{\label{foot:exception} Only for Example \ref{ex:poisson Lshaped}.} & 20\\
            \hline
            $\Psi$ & $r^{\mu}\sin(\mu\theta)$\\
            \hline
            $\mu_{j}$\footref{foot:exception} & $2j/3$\\
            \hline
            training data & \begin{tabular}{@{}c@{}}tensor-product Gauss-Legendre $128 \times 128$ in each square\\ quadrant of $\Omega$; Gauss-Legendre 128 on each edge of $\partial\Omega$ \end{tabular}\\
            \hline
            learning rate & $4\cdot 10^{-3}/1.1^{i-1}$\\
            \hline
            $\delta$ & $10^{3}$\\
            \hline
        \end{tabular}
        \label{tab:ex2-3}
    \end{table}

    \vspace{-3mm}
    \subsection{Example \ref{ex:triangular}}
    \hphantom{0mm}
    \vspace{-3mm}
    \begin{table}[H]
        \small
        \centering
        \begin{tabular}{||c||c||}
            \hline 
            \textbf{Parameter} & \\
            \hline \hline
            $u_{0}(x)$ & 0\\
            \hline
            $\mathbf{n}^{(i)}$ & $20\cdot 2^{i-1}$\\
            \hline
            $L_{i}$ & 1\\
            \hline
            $\sigma_{i}$ & $\tanh((1 + \frac{i}{4})t)$\\
            \hline
            $\mathcal{M}$ & 1 in each non-convex corner and at $(0,-3)$\\
            \hline
            $\Psi$ & Equations \eqref{eq:velocity dirichlet} and \eqref{eq:pressure dirichlet}\\
            \hline
            $\mu_{j}$ & Equation \eqref{eq:eigenvalues}\\
            \hline
            training data & \begin{tabular}{@{}c@{}}tensor-product Gauss-Legendre $128 \times 128$ in each square\\ quadrant of channel and rectangle enclosing triangular cavity\\ of $\Omega$\tablefootnote{The square quadrants are given by $(-2,0)\times (0,2)$ and $(0,2)\times (0,2)$ in the channel. For the triangular cavity, a $128\times 128$ Gauss-Legendre quadrature rule is generated in the rectangle $(-1,1)\times (0,-3)$ and any node falling outside of the cavity has its weight set to 0.}; Gauss-Legendre 128 on each edge of $\partial\Omega$\end{tabular}\\
            \hline
            learning rate & $3\cdot 10^{-3}/1.1^{i-1}$\\
            \hline
            $\delta$ & $10^{3}$\\
            \hline
        \end{tabular}
        \label{tab:ex4}
    \end{table}

    \vspace{-3mm}
    \subsection{Example \ref{ex:poisson Lshaped training}}
    \hphantom{0mm}
    \vspace{-3mm}
    \begin{table}[H]
        \small
        \centering
        \begin{tabular}{||c||c||}
            \hline 
            \textbf{Parameter} & \\
            \hline \hline
            $u_{0}(x)$ & 0\\
            \hline
            $\mathbf{n}^{(i)}$ & $40\cdot 2^{i-1}$\\
            \hline
            $L_{i}$ & 1\\
            \hline
            $\sigma$ & $\tanh((1 + \frac{i}{4})t)$\\
            \hline
            $m_{*}$ & 1\\
            \hline
            $\Psi$ & $r^{\mu}\sin(\mu\theta)$\\
            \hline
            $\mu$ & Initialized in $\mathcal{U}(0,1)$\\
            \hline
            training data & \begin{tabular}{@{}c@{}}tensor-product Gauss-Legendre $128 \times 128$ in each square\\ quadrant of $\Omega$; Gauss-Legendre 128 on each edge of $\partial\Omega$ \end{tabular}\\
            \hline
            learning rate & $4\cdot \times 10^{-3}/1.1^{i-1}$\\
            \hline
            $\delta$ & $10^{3}$\\
            \hline
        \end{tabular}
        \label{tab:ex5}
    \end{table}

    \subsection{Examples \ref{ex:stokes learning1}-\ref{ex:stokes learning2}} \label{app:last example}
    \hphantom{0mm}
    \vspace{-3mm}
    \begin{table}[H]
        \small
        \centering
        \begin{tabular}{||c||c||}
            \hline 
            \textbf{Parameter} & \\
            \hline \hline
            $u_{0}(x)$ & 0\\
            \hline
            $\mathbf{n}^{(i)}$ & $40\cdot 1.9^{i-1}$\\
            \hline
            $L_{i}$ & 1\\
            \hline
            $\sigma$ & $\tanh((1 + \frac{i}{4})t)$\\
            \hline
            $m_{*}$ & 1\\
            \hline
            $\Psi$ & Equation \eqref{eq:no slip eigenfunctions} in \ref{ex:stokes learning1}; \eqref{eq:moffatt complex} and \eqref{eq:moffatt pressure} in \ref{ex:stokes learning2} \\
            \hline
            $\mu$ & \begin{tabular}{@{}c@{}}Initialized in $\mathcal{U}(\lambda_{\text{exact}}-1/3, \lambda_{\text{exact}}+1/3)$ for \ref{ex:stokes learning1}; initialized\\ in $\mathcal{U}((\xi_{\text{exact}}-1/2, \xi_{\text{exact}}+1/2) \times (\zeta_{\text{exact}}-1/2, \zeta_{\text{exact}}+1/2))$\\
            for \ref{ex:stokes learning2}\end{tabular}\\
            \hline
            training data & \begin{tabular}{@{}c@{}}tensor-product polar coordinate Gauss-Legendre $128 \times 128$ in \\ $\Omega$ and 512 Gauss-Legendre nodes along $\partial\Omega$ for \ref{ex:stokes learning1};\\ Gauss-Legendre $128\times 128$ in $\Omega$\tablefootnote{In Example \ref{ex:stokes learning2}, quadrature rule is generated in the rectangle $(-1,1)\times (0,-3)$ and any points falling outside of $\Omega$ have their weights set to 0.} and Gauss-Legendre 256\\ on each edge of $\partial\Omega$ in \ref{ex:stokes learning2}\end{tabular}\\
            \hline
            learning rate & $3\cdot 10^{-3} / 1.1^{i-1}$\\
            \hline
            $\delta$ & $10^{3}$\\
            \hline
        \end{tabular}
        \label{tab:ex6}
    \end{table}

    \vspace{-7mm}
    \bibliographystyle{siamplain}
    \bibliography{references}

\begin{thebibliography}{10}

\bibitem{gnn1}
{\sc M.~Ainsworth and J.~Dong}, {\em Galerkin neural networks: A framework for
  approximating variational equations with error control}, SIAM Journal on
  Scientific Computing, 43 (2021), pp.~A2474--A2501.

\bibitem{gnn2}
{\sc M.~Ainsworth and J.~Dong}, {\em Galerkin neural network approximation of
  singularly-perturbed elliptic systems}, Computer Methods in Applied Mechanics
  and Engineering, 402 (2022), p.~115169.

\bibitem{arzani2023theory}
{\sc A.~Arzani, K.~W. Cassel, and R.~M. D'Souza}, {\em Theory-guided
  physics-informed neural networks for boundary layer problems with singular
  perturbation}, Journal of Computational Physics, 473 (2023), p.~111768.

\bibitem{bar2019learning}
{\sc Y.~Bar-Sinai, S.~Hoyer, J.~Hickey, and M.~P. Brenner}, {\em Learning
  data-driven discretizations for partial differential equations}, Proceedings
  of the National Academy of Sciences, 116 (2019), pp.~15344--15349.

\bibitem{batchelor1967introduction}
{\sc G.~K. Batchelor}, {\em An introduction to fluid dynamics}, Cambridge
  university press, 1967.

\bibitem{blumrannacher}
{\sc H.~Blum, R.~Rannacher, and R.~Leis}, {\em On the boundary value problem of
  the biharmonic operator on domains with angular corners}, Mathematical
  Methods in the Applied Sciences, 2 (1980), pp.~556--581.

\bibitem{bochev}
{\sc P.~B. Bochev and M.~D. Gunzburger}, {\em Least-squares finite element
  methods}, vol.~166, Springer Science \& Business Media, 2009.

\bibitem{dissanayake}
{\sc M.~Dissanayake and N.~Phan-Thien}, {\em Neural-network-based
  approximations for solving partial differential equations}, Communications in
  Numerical Methods in Engineering, 10 (1994), pp.~195--201.

\bibitem{schwabguo}
{\sc B.~Guo and C.~Schwab}, {\em Analytic regularity of stokes flow on
  polygonal domains in countably weighted sobolev spaces}, Journal of
  Computational and Applied Mathematics, 190 (2006), pp.~487--519.

\bibitem{hornik1990universal}
{\sc K.~Hornik, M.~Stinchcombe, and H.~White}, {\em Universal approximation of
  an unknown mapping and its derivatives using multilayer feedforward
  networks}, Neural networks, 3 (1990), pp.~551--560.

\bibitem{kidger2020universal}
{\sc P.~Kidger and T.~Lyons}, {\em Universal approximation with deep narrow
  networks}, in Conference on learning theory, PMLR, 2020, pp.~2306--2327.

\bibitem{kondratev}
{\sc V.~A. Kondrat'ev}, {\em Boundary value problems for elliptic equations in
  domains with conical or angular points}, Trudy Moskovskogo Matematicheskogo
  Obshchestva, 16 (1967), pp.~209--292.

\bibitem{mazya}
{\sc V.~Kozlov, V.~G. Maz'ya, and J.~Rossmann}, {\em Spectral problems
  associated with corner singularities of solutions to elliptic equations},
  American Mathematical Soc., 2001.

\bibitem{lagaris}
{\sc I.~E. Lagaris, A.~Likas, and D.~I. Fotiadis}, {\em Artificial neural
  networks for solving ordinary and partial differential equations}, IEEE
  transactions on neural networks, 9 (1998), pp.~987--1000.

\bibitem{fourieronet}
{\sc Z.~Li, N.~B. Kovachki, K.~Azizzadenesheli, K.~Bhattacharya, A.~Stuart,
  A.~Anandkumar, et~al.}, {\em Fourier neural operator for parametric partial
  differential equations}, in International Conference on Learning
  Representations, 2020.

\bibitem{deeponet}
{\sc L.~Lu, P.~Jin, G.~Pang, Z.~Zhang, and G.~E. Karniadakis}, {\em Learning
  nonlinear operators via deeponet based on the universal approximation theorem
  of operators}, Nature Machine Intelligence, 3 (2021), pp.~218--229,
  \url{https://doi.org/10.1038/s42256-021-00302-5}.

\bibitem{moffatt}
{\sc H.~K. Moffatt}, {\em Viscous and resistive eddies near a sharp corner},
  Journal of Fluid Mechanics, 18 (1964), pp.~1--18.

\bibitem{pakravan2021solving}
{\sc S.~Pakravan, P.~A. Mistani, M.~A. Aragon-Calvo, and F.~Gibou}, {\em
  Solving inverse-pde problems with physics-aware neural networks}, Journal of
  Computational Physics, 440 (2021), p.~110414.

\bibitem{petersen2021topological}
{\sc P.~Petersen, M.~Raslan, and F.~Voigtlaender}, {\em Topological properties
  of the set of functions generated by neural networks of fixed size},
  Foundations of computational mathematics, 21 (2021), pp.~375--444.

\bibitem{pinns}
{\sc M.~Raissi, P.~Perdikaris, and G.~E. Karniadakis}, {\em Physics-informed
  neural networks: A deep learning framework for solving forward and inverse
  problems involving nonlinear partial differential equations}, Journal of
  Computational physics, 378 (2019), pp.~686--707.

\bibitem{boundarypenalty}
{\sc S.~A. Sauter and C.~Schwab}, {\em Boundary element methods}, Springer,
  2011.

\bibitem{schwab1998p}
{\sc C.~Schwab}, {\em p-and hp-finite element methods: Theory and applications
  in solid and fluid mechanics}, Clarendon Press, 1998.

\bibitem{genfem}
{\sc T.~Strouboulis, K.~Copps, and I.~Babu{\v{s}}ka}, {\em The generalized
  finite element method}, Computer methods in applied mechanics and
  engineering, 190 (2001), pp.~4081--4193.

\bibitem{xfem}
{\sc N.~Sukumar, N.~Mo{\"e}s, B.~Moran, and T.~Belytschko}, {\em Extended
  finite element method for three-dimensional crack modelling}, International
  journal for numerical methods in engineering, 48 (2000), pp.~1549--1570.

\bibitem{yang2022multi}
{\sc J.~Yang, K.~Mittal, T.~Dzanic, S.~Petrides, B.~Keith, B.~Petersen,
  D.~Faissol, and R.~Anderson}, {\em Multi-agent reinforcement learning for
  adaptive mesh refinement}, arXiv preprint arXiv:2211.00801,  (2022).

\bibitem{wan}
{\sc Y.~Zang, G.~Bao, X.~Ye, and H.~Zhou}, {\em Weak adversarial networks for
  high-dimensional partial differential equations}, Journal of Computational
  Physics, 411 (2020), p.~109409.

\end{thebibliography}

    \clearpage
	
\end{document}